\documentclass{amsart}
\usepackage[utf8]{inputenc}

\usepackage{fullpage}
\usepackage[utf8]{inputenc}
\usepackage[T1]{fontenc}
\usepackage{graphicx,subfigure} 
\usepackage{amsmath,amsfonts,amssymb,mathtools}
\usepackage{stmaryrd}
\usepackage{mathrsfs,dsfont}
\usepackage{amsthm}
\usepackage{mathabx}
\usepackage{tabularx}
\usepackage{float}
\usepackage[usenames]{color}
\allowdisplaybreaks
\usepackage{amscd}
\usepackage{ulem}

\numberwithin{equation}{section}

\newtheorem{theo}{Theorem}
\newtheorem{lemma}{Lemma}[section]
\newtheorem{assum}{Assumption}[section]
\newtheorem{rem}{Remark}[section]
\newtheorem{propo}{Proposition}[section]
\newtheorem{exam}{\it Example}[section]
\newtheorem{cor}{Corollary}[section]

\newcommand{\N}{\mathbb{N}}
\newcommand{\R}{\mathbb{R}}
\newcommand{\E}{\mathbb{E}}

\newcommand{\D}{\mathbb{D}}

\newcommand{\IC}{\mathcal{C}}
\newcommand{\IL}{\mathcal{L}}

\newcommand{\IW}{\mathcal{W}}

\newcommand{\IB}{\mathcal{B}}

\newcommand{\ID}{\mathcal{D}}
\newcommand{\IS}{\mathcal{S}}
\newcommand{\GP}{\mathfrak{p}}

\title{Weak convergence orders for a fully discrete approximation of semilinear stochastic PDEs driven by multiplicative Gaussian and Poisson noise}

\author{Chuchu Chen}
\address{State Key Laboratory of Mathematical Sciences, Academy of Mathematics and Systems Science, Chinese Academy of Sciences, Beijing 100190, China; School of Mathematical Sciences, University of Chinese Academy of Sciences, Beijing 100049, China}
\email{chenchuchu@lsec.cc.ac.cn}

\author{Jialin Hong}
\address{State Key Laboratory of Mathematical Sciences, Academy of Mathematics and Systems Science, Chinese Academy of Sciences, Beijing 100190, China; School of Mathematical Sciences, University of Chinese Academy of Sciences, Beijing 100049, China}
\email{hjl@lsec.cc.ac.cn}

\author{Ziyi Lei}
\address{State Key Laboratory of Mathematical Sciences, Academy of Mathematics and Systems Science, Chinese Academy of Sciences, Beijing 100190, China; Department of Mathematics, College of Sciences, China University of Mining and Technology, Beijing 100083, China}
\email{ziyilei@lsec.cc.ac.cn}

\author{Xu Wang}

\address{State Key Laboratory of Mathematical Sciences, Academy of Mathematics and Systems Science, Chinese Academy of Sciences, Beijing 100190, China; School of Mathematical Sciences, University of Chinese Academy of Sciences, Beijing 100049, China}

\email{wangxu@lsec.cc.ac.cn}

\thanks{This work is funded by the National Key R\&D Program of China (No. 2024YFA1015900 and No. 2024YFA1012300), the National Natural Science Foundation of China (No. 12031020, No. 12461160278, No. 12471386,  and No. 12288201),  the Youth Innovation Promotion Association CAS, China, and by the CAS Project for Young Scientists in Basic Research (YSBR-087)}

\keywords {Fully discrete scheme $\cdot$ Weak convergence order $\cdot$ Multiplicative L\'evy noise $\cdot$ Stochastic partial differential equation}

\begin{document}

\maketitle

\begin{abstract}

 This paper investigates the weak convergence of a fully discrete approximation for a class of semilinear  stochastic partial differential equations (SPDEs) driven by multiplicative Gaussian and Poisson noise. The scheme combines a spectral Galerkin discretization in space with an exponential Euler discretization in time. Compared with single-noise settings, the simultaneous presence of continuous and jump perturbations introduces additional challenges arising from discontinuous trajectories, nonlocal jump operators, and their interaction with state-dependent multiplicative coefficients. To overcome these difficulties, we establish the exponential integrability of the numerical solution processes and refined moment estimates in negative Sobolev norms for the associated first- and second-order variational processes. Combining these estimates with a discrete-time approximation argument and Wiener--Poisson Malliavin calculus, we derive regularity estimates for the solutions of the Kolmogorov equations associated with the spatial Galerkin approximations. Based on these estimates, we obtain the temporal and spatial weak convergence orders arbitrarily close to one with respect to discretization parameters $\tau$ and $\lambda_N^{-1}$, respectively. To the best of our knowledge, this provides the first weak convergence result for a fully discrete approximation of semilinear SPDEs with simultaneous multiplicative Gaussian and Poisson noise.

\end{abstract}

\section{Introduction}\label{sec.1}

In this work, we consider the numerical analysis for SPDEs of the form:
\begin{equation}\label{SPDE:MultiplicativeNoise}
\!\!\!\left\{
 \begin{aligned}
&dX(t)=(AX(t)+F(X(t)))dt+B(X(t))dW^Q(t)+\int_\chi G(X(t-), z)\tilde{\mathbf{N}}(dz,dt),\quad t\in(0,T],~~T>0,\\
&X(0)=x_0,
\end{aligned}   
\right.
\end{equation}
where $A:=\Delta:\mbox{Dom}(A)\subset H\to H$ is the Laplacian with homogeneous Dirichlet boundary condition, and $H:=L^2(\mathcal O)$ with an open bounded domain $\mathcal{O}\subset\mathbb R^d, d=1,2,3$ and smooth  boundary $\partial \mathcal{O}$, is endowed with the usual inner product $\langle\cdot,\cdot\rangle$ and norm $\|\cdot\|$. The mark set is denoted by $\chi:=H\backslash\{0\}$. Here the driving noises $\{W^Q(t)\}_{t\in[0,T]}$ and $\{\tilde{\mathbf{N}}(\cdot,t)\}_{t\in[0,T]}$ are, respectively, a Gaussian process and a compensated Poisson random measure. These types of noises model different sources of randomness: Gaussian noise captures continuous stochastic fluctuations, while the compensated Poisson random measure describes sudden, discontinuous perturbations. This framework covers a large class of L\'evy-driven SPDEs, in view of the L\'evy--It\^o decomposition theorem. The precise description of these driving noises, as well as the assumptions on coefficients $F, B$ and $G$, will be given in Section~\ref{sec2}. Since the explicit expression for the solution of \eqref{SPDE:MultiplicativeNoise} is generally unavailable, numerical methods play a crucial role in approximating the solution. In addition to the strong error, which measures convergence in a specified $Lp$ sense (cf. \cite{JA-2019,CTJZ-2025}), one is often interested in approximating the expected values of certain functionals of the solution $X$, which leads to the study of weak convergence for the proposed numerical method.

 The literature on the weak convergence for numerical approximations of SPDEs  of the form~\eqref{SPDE:MultiplicativeNoise} primarily focuses on equations driven solely by Gaussian noise ($G\equiv0$) or by pure-jump L\'evy noise ($B\equiv0$). For Gaussian-driven SPDEs, a substantial body of  weak convergence results is available; see, e.g., \cite{CH-2019,CHS-2021,AS-2016,DAR-2018,AD-2011,CE-2018,Kruse-2014,AJ-2009,MMS-2009,J-2003,Zhang2026,Brehier2026a,Wang2016,Cai2021} and the references therein. The analyses in these works typically employ several key tools, including the Kolmogorov equation, Malliavin calculus, and mild It\^o's formula. In the purely Gaussian setting, these approaches benefit particularly from the regularity induced by the continuous perturbations and the local diffusion structure.  By contrast, when the Poisson noise is taken into account (i.e., $G\neq0$), one has to deal instead with discontinuous perturbations and non-local jump operators. Consequently, within analytical frameworks based on the Kolmogorov equation and Malliavin calculus, both the weak error representation and the regularity propagation mechanism become more involved,  so that the weak convergence analysis developed for the purely Gaussian noise case cannot be directly extended to the jump setting or to the more general Gaussian–Poisson framework.

For SPDEs driven by pure-jump L\'evy noise, available weak convergence results are more limited and largely focus on linear or semilinear equations with additive noise; see, e.g., \cite{AT2018,FR-2013,MFR-2015,AF2019}. These works rely on the special structure of the additive noise, which is essentially different from that of the setting of \eqref{SPDE:MultiplicativeNoise} with the multiplicative Gaussian and Poisson noises. In our setting, the main difficulty lies not only in the jumps themselves, but also in their coupling with the state-dependent multiplicative coefficients. This coupling may propagate and accumulate jump-induced irregularities along the stochastic evolution, thereby deteriorating the temporal regularity required in the weak error analysis. To the best of our knowledge, no weak convergence result is currently available for numerical approximations of \eqref{SPDE:MultiplicativeNoise} driven simultaneously by the multiplicative Gaussian noise and Poisson random measure.

The aim of this work is to fill this gap by developing a unified numerical and analytical framework that incorporates both types of noise and establishes weak convergence results for fully discrete approximations of \eqref{SPDE:MultiplicativeNoise}. We propose the fully discrete approximation $\{\bar X_m^N\}_{m\in\{0,1,\ldots,M\}}$, which combines the spectral Galerkin method in space with the exponential Euler method in time and is defined by
 \begin{align}\label{eq:one-step-I}
	\bar{X}^{N}_{m+1}&=E(\tau)\Big(\bar{X}^{N}_m+\tau F_N(\bar{X}^{N}_m)+P_NB(\bar{X}^{N}_m)\Delta W^Q_{m}+\int_{t_m}^{t_{m+1}}\int_\chi G_N(\bar{X}^{N}_m,z)\tilde{\mathbf{N}}(dz,ds)\Big),\quad
\end{align}
where $ \bar{X}_0^N=P_Nx_0$ and coefficients are defined by $F_N := P_N F$ and $G_N := P_N G$ with $P_N$ being the orthogonal projection from $H$ onto its finite-dimensional subspace, the time step size $\tau:=T/M$, and $\Delta W^Q_{m}:=W^Q(t_{m+1})-W^Q(t_m)$. Our main result establishes the weak convergence orders of the fully discrete approximation~\eqref{eq:one-step-I}, which are arbitrarily close to one with respect to both the time step size $\tau$ and the inverse eigenvalue $\lambda_N^{-1}$.

\begin{theo}\label{the:weakconvergence-full}
Suppose that Assumptions \ref{assum:nonlinearity-F}--\ref{assum:nonlinearity-G} hold and let  $\varphi\in\IC^3_{\mathfrak{p}}(H;\R)$. Then for any sufficiently small  $\epsilon>0$ and $x_0\in \dot{H}^{\frac{1-\epsilon/2}{2}}$, there exists a constant $C_{\epsilon,T}>0$ such that 
\begin{align*}
	&\big|\E[\varphi(X(T))]-\E[\varphi(\bar{X}^N_M)]\big|\le C_{\epsilon,T}\big(\lambda_N^{-1+3\epsilon}+ \tau^{1-3\epsilon}\big)
\end{align*}
for all $N\in\N_+$ and $\tau\in(0, \tau_0)$ satisfying the CFL condition $\tau\lambda_N<\frac{\sqrt{\ln 2}}{4}$  with $\tau_0$ given in \eqref{def:tau_0}.
\end{theo}

The key to deriving the above weak convergence orders lies in a set of new regularity estimates. Based on the Kolmogorov equations and Wiener--Poisson Malliavin calculus, we use a discrete-time analytical framework, which transforms the analysis of continuous-time stochastic integrals into recursive estimates driven by independent increments. This enables us to fully exploit the  independence and martingale structure of both Gaussian and jump noises to provide sharper estimates. Within this framework, we first establish the exponential integrability of the solution processes in the presence of Poisson jumps, and then derive refined moment estimates in  negative Sobolev norms for the first- and second-order variational processes associated with the numerical scheme. Further, by combining these estimates with Wiener and Poisson Malliavin integration-by-parts formulas, we are able to recover additional regularity of the solution to the associated Kolmogorov equations, thereby effectively overcoming the difficulties posed by the low regularity induced by multiplicative Gaussian and Poisson noises. These estimates then provide the key analytical ingredients in the rigorous derivation of weak convergence orders for the fully discrete approximation~\eqref{eq:one-step-I} of \eqref{SPDE:MultiplicativeNoise}.

The rest of the paper is organized as follows. The preliminaries, including notations, assumptions, the well-posedness of the exact solution, and a brief introduction to the Malliavin calculus are given in Section \ref{sec2}. Section \ref{sec:3} is devoted to the introduction of the spatio-temporal full-discretization, the properties of the corresponding solution processes, as well as the estimates of Malliavin derivatives. In Section \ref{sec:4}, we present the statements of the associated Kolmogorov equations and the regularity estimates on the solution of the Kolmogorov equations. Finally, the proof of Theorem \ref{the:weakconvergence-full} is provided in Section~\ref{sec.5}.

\section{Preliminaries}\label{sec2}

This section fixes the basic settings and assumptions, and introduces the definitions and properties of Malliavin calculus that furnish the basis and key analytical tools for our subsequent analysis.

\subsection{Settings and assumptions}\label{sec:2.1}
Let $\mathbb{N}$ be the set of nonnegative integers and $\mathbb{N}_+:=\mathbb{N}\backslash\{0\}$. Set $a\vee b:=\max\{a,b\}$ and $a\wedge b:=\min\{a,b\}$. For an arbitrary $k\in\N$, denote by $\IC^k_b(H;\R)$ the space of bounded continuous mappings from $H$ to $\R$ that are $k$ times continuously Fr\'echet differentiable with bounded derivatives up to order $k$.  Let $\IB_b(H;\R)$ be the space of bounded and Borel measurable functions from $H$ to $\R$.

Denote by $\IC^2_{\mathfrak{p}}(H;\R)$ the space of twice continuously Fr\'echet differentiable mappings $h:H\to \R$ with $h$ and its derivatives satisfying polynomial growth, i.e., there exist constants $\mathfrak{p}\in\N_+$ and  $C_{\GP}>0$ such that 
\begin{align}
   |h(x)|\vee \sup_{y_1\in H\setminus\{0\}}\frac{|Dh(x).y_1|}{\|y_1\|}\vee \sup_{y_1,y_2\in H\setminus\{0\}}\frac{|D^{2}h(x).(y_1,y_2)|}{\|y_1\|\|y_2\|}\le C_{\GP}(1+\|x\|)^{\mathfrak{p}},\quad x\in H, 
\end{align}
where $D^{j}h$ is the $j$th order Fr\'echet derivative of $h$ identified as
\[D^jh(x).(y_1,\ldots,y_j)=\langle D^{j}h(x).(y_1,\ldots,y_{j-1}),y_j\rangle,\quad \mbox{for}~x,y_j\in H,~ j\in\mathbb{N}_+\]
 with $Dh(x).y_1=\langle Dh(x),y_1\rangle$. 
 Define further a subspace $\IC^3_{\GP}(H;\R)\subset \IC^2_{\mathfrak{p}}(H;\R)$ consisting of three times continuously Fr\'echet differentiable mappings from $H$ to $\mathbb{R}$ whose third derivative satisfies
 \begin{align}\label{setting-varrho1}
  |D^{3}h(x).(h_1,h_2,h_3)|\le C_{\GP}(1+\|x\|)^{\GP}\|h_1\|\|h_2\|\|h_3\|_{\dot{H}^{\varrho_1}}, \quad   x,h_1,h_2\in H, ~h_3\in \dot{H}^{\varrho_1}
\end{align}
for some $\varrho_1\in (\frac12,1)$. Here $\dot{H}^{\varrho_1}$ will be introduced in Subsection \ref{subsec:2.1.1}. 

Let $\IL(U;H)$ and $\IL_2(U;H)$ be spaces of bounded linear operators and Hilbert--Schmidt operators, respectively, from separable Hilbert spaces $U$ to $H$. The associated norms are defined as
\begin{align*}
  \|T\|_{\IL(U;H)}:=\sup_{x\in U,\|x\|_U\le 1}\|Tx\|,\quad   \|S\|_{\IL_2(U;H)}:=\big(\sum_{i\in\N_+}\|Sf_i\|^2\big)^{\frac12}
\end{align*}
for $T\in \IL(U;H)$ and $S\in \IL_2(U;H)$,
where $\{f_i\}_{i\in\N_+}$ is an orthonormal basis of $U$ and $\|\cdot\|_U$ denotes the norm on $U$. In particular, we use the notations $\IL(H):=\IL(H;H)$, $\IL_2(H):=\IL_2(H;H)$, and $\IL_2^0(H):=\IL_2(U_0;H)$ with $U_0:=Q^{\frac12}(H)$ for simplicity. Here the covariance operator $Q\in\IL(H)$ is assumed to be symmetric and non-negative definite. The Hilbert space $U_0$ is endowed with the inner product $\langle u,v\rangle_{U_0}:=\langle Q^{-\frac12}u,Q^{-\frac12}v\rangle$ for all $u,v\in U_0$, where $Q^{-\frac12}$ is the pseudoinverse of $Q^{\frac12}$. If $T\in\IL(H)$ and $S\in \IL_2(H)$, then $TS\in \IL_2(H)$ and $\|TS\|_{\IL_2(H)}\le \|T\|_{\IL(H)}\|S\|_{\IL_2(H)}$.  Let $\IL_1(H):=\IL_1(H;H)$ be the space of nuclear operators from $H$ to $H$. The trace of operator $S\in\IL_1(H)$ is defined as $\mbox{Tr}(S):=\sum_{i\in\N_+}\langle Se_i,e_i\rangle$ with $\{e_i\}_{i\in\N_+}$ being an orthonormal basis of $H$. 

Throughout the paper, $C$ denotes a generic positive  constant which may vary from one occurrence to another, and numbered subscripts in $C_{a, b, \ldots}$ indicate that the constant depends on parameters $a, b, \ldots$  

In the following subsections, we will introduce definitions of the linear operator $A$, and the driving noises $W^Q$ and $\tilde{\mathbf{N}}$, as well as assumptions imposed on the coefficients $F$, $B$, and $G$.

\subsubsection{Linear operator}\label{subsec:2.1.1}
Let $A=\Delta$ with domain $\mbox{Dom}(A)=H^2(\mathcal{O})\bigcap H^1_0(\mathcal{O})\subset H$ denote the Laplacian with homogeneous Dirichlet boundary condition, where $H^2(\mathcal{O})$ and $H^1_0(\mathcal{O})$ are the standard Sobolev spaces equipped with their usual norms. It possesses a non-decreasing sequence $\{\lambda_i\}_{i\in\N_+}$ of positive eigenvalues satisfying $-Ae_i=\lambda_ie_i$ for all $i\in\N_+$. Moreover, Weyl's law yields $\lambda_i\sim c_Di^{2/d}
~\text{as }i\to\infty,$ for some constant \(c_D>0\).

For $s\ge 0$, let $\dot{H}^s$ be the Sobolev space generated by $(-A)^{\frac{s}{2}}$, endowed with the inner product $\langle u,v\rangle_{\dot{H}^s}:=\langle(-A)^{\frac{s}{2}}u,(-A)^{\frac{s}{2}}v\rangle$ and the norm $\|u\|_{\dot{H}^s}:=\sqrt{\langle u,u\rangle_{\dot{H}^s}}$ for $u,v\in \dot{H}^s$. For \(s<0\), we define $\dot H^s:=(\dot H^{-s})^*,$ where the duality is taken with respect to the pivot space \(H\). Equivalently, \(\dot H^s\) is the completion of \(H\) under the norm
\[
\|u\|_{\dot H^s}^2=\sum_{i\in\mathbb N_+}\lambda_i^s
|\langle u,e_i\rangle_H|^2.
\]

The linear operator $A$ is the infinitesimal generator of the semigroup $\{E(t):=e^{tA}\}_{t\geq0}$ on $H$, defined by $E(t)x=e^{tA}x=\sum_{i\in\N_+}e^{-t\lambda_i}\langle x, e_i\rangle e_i$ for $t\ge 0$ and $x\in H$. The operator $E(t)$ is  bounded and self-adjoint on $H$ with $\|E(t)\|_{\IL(H)}\leq e^{-t\lambda_1}\leq 1$. 
Recall also the smoothing properties for $\{E(t)\}_{t\ge0}$ (cf. \cite{Kruse-2014}): 
\begin{align}
\|(-A)^{a}E(t)\|_{\IL(H)}&\leq Ct^{-a},\quad t\in(0,T],~a\ge 0,~T>0,\label{prop:semigroup1}\\
\|(-A)^{-b}(E(t)-E(s))\|_{\IL(H)}&\leq C(t-s)^{b},\quad 0\leq s<t\leq T, ~T>0,~b\in[0,1].\label{prop:semigroup2}
\end{align}

\subsubsection{Driving noises}
For \(T>0\), let $W^Q$ and $\tilde{\mathbf{N}}$ denote, respectively, the Gaussian and jump driving noises in \eqref{SPDE:MultiplicativeNoise}, defined on a filtered probability space \((\Omega,\mathcal F,\{\mathcal F_t\}_{t\in[0,T]},\mathbb P)\), and assume that they are independent.

More precisely, $\{W^Q(t)\}_{t\in[0,T]}$ is a $Q$-Wiener process with the covariance operator $Q$ being a linear, self-adjoint and positive semi-definite nuclear operator on $H$ that commutes with $A$. Then the Karhunen--Lo\`eve expansion yields
\begin{align*}
W^Q(t)=\sum_{i\in\N_+}\sqrt{q_i}\beta_i(t)e_i,\quad t\in[0,T],
\end{align*}
where $q_i$, $i\in\mathbb{N}_+$ are eigenvalues of $Q$ satisfying $Qe_i=q_ie_i$ and $\{\beta_i\}_{i\in\N_+}$ is a family of independent real-valued standard Wiener processes. 

Based on the mark space $\chi:=H\setminus\{0\},$
let \(\nu\) be a L\'evy measure on \((H,\mathcal B(H))\) satisfying $\nu(\{0\})=0$ and $\int_H (1\wedge \|z\|^2)\,\nu(dz)<\infty$.
We identify \(\nu\) with its restriction to \((\chi,\mathcal B(\chi))\). Set
\[
\chi_T:=[0,T]\times\chi,
\qquad
\Xi_T:=\mathcal B([0,T])\otimes\mathcal B(\chi),
\]
and equip \(\chi_T\) with the product measure \(dt\otimes\nu(dz)\). Let \(\mathbf{N}\) be a Poisson random measure on \(\chi_T\) with intensity measure \(\nu(dz)\,dt\), and define the associated compensated Poisson random measure by
\[
\tilde{\mathbf{N}}(dz,dt):=\mathbf{N}(dz,dt)-\nu(dz)\,dt.
\]
For $\hat\Xi_T:=\bigl\{B\in\Xi_T:(dt\otimes\nu)(B)<\infty\bigr\},$
the random variable \(\mathbf{N}(B)\) is finite almost surely for every \(B\in\hat\Xi_T\). If \(\nu(\Lambda)=\infty\) for some \(\Lambda\in\mathcal B(\chi)\), then $\mathbf{N}((0,t]\times\Lambda)=\infty~\text{a.s.}$ for every \(t>0\). We work with the completed, right-continuous natural filtration generated jointly by \(W^Q\) and \(\mathbf{N}\). More precisely, let
\[
\mathcal G_t^0
:=
\sigma\bigl(W^Q(s):0\le s\le t\bigr)
\vee
\sigma\bigl(\mathbf{N}((0,s]\times\Lambda):
0\le s\le t,\,
\Lambda\in\mathcal B(\chi),\,
\nu(\Lambda)<\infty
\bigr),
\]
and let \(\{\mathcal F_t\}_{t\in[0,T]}\) be the usual augmentation of
\(\{\mathcal G_t^0\}_{t\in[0,T]}\). In particular, for \(0\le s<t\le T\), the Wiener increment $W^Q(t)-W^Q(s)$ and the restriction of \(\mathbf{N}\) to \((s,t]\times\chi\) are independent of \(\mathcal F_s\).

Throughout the paper, all stochastic integrands are understood in their predictable versions. In particular, whenever a c\`adl\`ag process appears in an integral with respect to the compensated Poisson random measure, its predictable left-limit version is used. Given predictable processes \(\Phi\) and \(\Psi\) satisfying
\[
\mathbb E\Big[\int_0^T
\|\Phi(s)\|_{\mathcal L_2^0(H)}^2\,ds\Big]<\infty\quad \mbox{and}\quad
\mathbb E\Big[\int_0^T\int_{\chi}
\|\Psi(s,z)\|_H^2\,\nu(dz)\,ds\Big]<\infty,
\]
integrals  $\int_0^T\Phi(s)dW^Q(s)$ and $\int_0^T\int_{\chi}\Psi(s,z)\,\tilde{\mathbf{N}}(dz,ds)$ are well-defined $H$-valued stochastic integrals. We refer to \cite[Section 2]{SJ_2013}, \cite[Section 2]{TE-2016} and references therein for more details on properties of $\mathbf{N}$.

\subsubsection{Coefficients}\label{subsec:2.1.3.}

We now impose the following assumptions on the coefficients $F$, $B$ and $G$, together with some valid examples.

\begin{assum}\label{assum:nonlinearity-F}
Let $F\in \IC^3(H;H)$  satisfy $\|F(0)\|\leq L_0$ and the following conditions.
\begin{itemize}
\item[(i)] For all $r_1\in(0,\frac14)$, and all $\varrho\in(0,\frac12-r_1]$, there exists a constant  $C_{r_1,\varrho}>0$ such that 
\begin{align*}
	\|(-A)^{r_1}F(x)\|\le C_{r_1,\varrho}(1+\|x\|_{\dot{H}^{2(r_1+\varrho)}}),\quad x\in \dot{H}^{2(r_1+\varrho)}.
\end{align*}
\item[(ii)] There exists a constant $L_F>0$ such that $\sup_{x\in H}\|DF(x).y\|\le L_F\|y\|$ for all $y\in H$. In addition,  there exist constants $\alpha_0\in(0,1)$ and $\varrho_0\in(\frac12,1)$ such that for any $\varsigma_0\in(0,1)$ and any $\varsigma\in(0,1)$, there exist positive constants $C_{\varsigma, \varsigma_0}$, $C_{\varsigma_0}$, and $C_{\varsigma_0,\varrho_0}$ such that 
 \begin{align}
&\|(-A)^{-\frac{1+\varsigma_0}{4}}DF(x).y\|\leq C_{\varsigma, \varsigma_0}(1+\|x\|_{\dot{H}^{\frac{1-\varsigma}{2}}})^{\alpha_0}\|y\|_{\dot{H}^{-\frac{1-\varsigma}{2}}},\quad x\in\dot{H}^{\frac{1-\varsigma}{2}}, ~y\in H,\nonumber\\
&\|(-A)^{-\frac{1+\varsigma_0}{4}}D^2F(x).(y,h)\|\leq C_{\varsigma_0}\|y\|\|h\|,\quad x, y,h\in H,\nonumber\\
&\|(-A)^{-\frac{1+\varsigma_0}{4}}D^3F(x).(y,h,k)\|\leq C_{\varsigma_0,\varrho_0}\|y\|\|h\|\|k\|_{\dot{H}^{\varrho_0}},\quad x,y,h\in H, ~k\in \dot{H}^{\varrho_0}.\nonumber
\end{align} 
\end{itemize}
\end{assum}

\begin{exam}
We give two examples of drift coefficients satisfying Assumption 2.1.

(i) \textbf{Linear drift.} Let $F(x)=ax,~ x\in H$ for some \(a\in\mathbb R\). Then \(F\in C^\infty(H;H)\), \(F(0)=0\), and 
$DF(x)\cdot y=ay,$ $D^2F(x)=D^3F(x)=0.$ For all \(r_1\in(0,1/4)\), \(\varrho\in(0,1/2-r_1]\), and
\(x\in\dot H^{2(r_1+\varrho)}\),
\[
\|(-A)^{r_1}F(x)\|
=
|a|\|x\|_{\dot H^{2r_1}}
\le
C_{r_1,\varrho}|a|
\|x\|_{\dot H^{2(r_1+\varrho)}}.
\]
Moreover, $\sup_{x\in H}\|DF(x)\cdot y\|\le |a|\|y\|,~ y\in H$, and for \(s_0,s\in(0,1)\),
\[
\|(-A)^{-\frac{1+s_0}{4}}DF(x).y\|=
|a|\|y\|_{\dot H^{-\frac{1+s_0}{2}}}\le C_{s,s_0}|a|\|y\|_{\dot H^{-\frac{1-s}{2}}}.
\]
Since \(D^2F=D^3F=0\), all the remaining conditions in Assumption~\ref{assum:nonlinearity-F} are immediate.

(ii) \textbf{A nonlinear diagonal drift.} Define
\[
F(x)
:=
\sum_{i\in\mathbb N_+}
\sin(\langle x,e_i\rangle_H)e_i,
\qquad x\in H.
\]
Then we derive that \(F\in C^3(H;H)\), \(F(0)=0\), and \begin{align*}
	DF(x).y=\sum_{i\in\N_+}\cos(\langle x, e_i\rangle)\langle y, e_i\rangle e_i, & \ \ \quad ~~
	D^2F(x).(y,h)=-\sum_{i\in\N_+}\sin(\langle x, e_i\rangle)\langle h, e_i\rangle \langle y, e_i\rangle e_i,\\
	D^3F(x).(y,h,k)&=-\sum_{i\in\N_+}\cos(\langle x, e_i\rangle)\langle k, e_i\rangle\langle h, e_i\rangle \langle y, e_i\rangle e_i.
\end{align*}
 Since $|\sin z|\le |z|$, for every
\(x\in\dot H^{2(r_1+\varrho)}\), we obtain
\begin{align*}
\|(-A)^{r_1}F(x)\|^2&=
\sum_{i\in\mathbb N_+}
\lambda_i^{2r_1}
\left|
\sin(\langle x,e_i\rangle_H)
\right|^2\le
\|x\|_{\dot H^{2r_1}}^2\le
C_{r_1,\varrho}
\|x\|_{\dot H^{2(r_1+\varrho)}}^2,
\end{align*}
which verifies Assumption~\ref{assum:nonlinearity-F} (i). Furthermore, $
\sup_{x\in H}\|DF(x).y\|
\le \|y\|$ for $y\in H.$ For any $\varsigma_0\in(0,1)$ and $\varsigma\in(0,1)$,  it holds that
\begin{align*}
	\|(-A)^{-\frac{1+\varsigma_0}{4}}DF(x).y\|^2&= \sum_{i\in\N_+}\lambda_i^{-\frac{1+\varsigma_0}{2}} \left( \cos(\langle x, e_i\rangle) \langle y, e_i\rangle\right)^2 \le C_{\varsigma,\varsigma_0} \sum_{i\in\N_+}\lambda_i^{-\frac{1-\varsigma}{2}} \langle y, e_i\rangle^2 = C_{\varsigma,\varsigma_0}  \|y\|_{\dot{H}^{-\frac{1-\varsigma}{2}}}^2,
\end{align*}
and all the conditions in Assumption~\ref{assum:nonlinearity-F} are satisfied.  Other standard constructions satisfying all  the conditions in Assumption~\ref{assum:nonlinearity-F} include the nonlinearity defined by $F(x):=(1+\|x\|_{\dot{H}^{-\frac{1}{2}}}^2)^{-1}(-A)^{-\frac{1}{4}}x$.
\end{exam}

\begin{assum}\label{assum:nonlinearity-B}
Let $B\in\IC^3(H;\IL_2^0(H))$  satisfy  $\|B(0)\|_{\IL_2^0(H)}\leq L_1$ and the following conditions.
\begin{itemize} 
\item[(i)] There exists some $\gamma\in[0,\frac12)$ such that for any $r\in[0,\frac{1}{4})$, there exists  a constant $C_r>0$ such that
\begin{align*}
\|(-A)^{r}B(x)\|_{\IL_2^0(H)}\le C_{r}(1+\|(-A)^{r}x\|^{\gamma}),\quad x\in \dot{H}^{2r}.
\end{align*} 
\item[(ii)] For any $\varsigma_1\in[0,\frac14]$, there exists a positive constant $C_{\varsigma_1}$ such that
\begin{align*}
&\sup_{x\in H}\|(-A)^{-\varsigma_1}DB(x).y\|_{\IL_2^0(H)}\le C_{\varsigma_1}\|(-A)^{-\varsigma_1}y\|,\quad y\in H,\\
&\sup_{x\in H}\|(-A)^{-\varsigma_1}D^2B(x).(y,h)\|_{\IL_2^0(H)}\le C_{\varsigma_1}\|(-A)^{-\varsigma_1}y\|\|(-A)^{-\varsigma_1}h\|,\quad y,h\in H.
\end{align*}
\item[(iii)]
There exists a positive constant $C_B$ such that
\begin{align*}
\sup_{x\in H}\|D^3B(x).(y,h,k)\|_{\IL_2^0(H)} \le C_{B}\|y\|\|h\|\|k\|,\quad y,h,k\in H.
\end{align*}
\end{itemize}
\end{assum}

\begin{exam} 

Set \(d=1\) and assume \(Q^{1/2}\in\mathcal L_2(H)\). For \(x\in \dot H^{2r}\), define the diagonal operator
\[
\widehat B(x):H\to H,
\qquad
\widehat B(x)e_i:=\tilde\sigma_i(x)e_i,
\qquad i\in\mathbb N_+,
\]
and let
\[
B(x):=\widehat B(x)\big|_{U_0}:U_0\to H,
\qquad
U_0:=Q^{1/2}(H).
\]
Let \(V\in C^\infty([0,\infty);\mathbb R)\) satisfy \(V(s)=1\) for \(s\in[0,1]\), \(V(s)\in[0,1]\) for \(s\in[1,2]\), and \(V(s)=0\) for \(s\ge2\). For \(r\in[0,1/4)\), define $
\tilde\sigma_i:\dot H^{2r}\to\mathbb R$
by
\[
\tilde\sigma_i(x):=\frac{1}{i+1}\Bigl[V\bigl(\|x\|_{\dot{H}^{-1/2}}^2\bigr)\bigl(\sqrt{1+\|x\|_{\dot{H}^{-1/2}}^2}-1\bigr)+\bigl(1-V(\|x\|_{\dot{H}^{-1/2}}^2)\bigr)\bigl(\sqrt{1+\|x\|_{\dot{H}^{-1/2}}^{2\gamma}}-1\bigr)\Bigr]
\] 
for some \(\gamma\in[0,1/2)\). By the definition of \(\mathcal L_2^0(H)\) and the ideal property of Hilbert--Schmidt operators,
\begin{align*}
\|(-A)^rB(x)\|_{\mathcal L_2^0(H)}^2&=
\|(-A)^r\widehat B(x)Q^{1/2}\|_{\mathcal L_2(H)}^2\le
\|(-A)^r\widehat B(x)\|_{\mathcal L(H)}^2
\|Q^{1/2}\|_{\mathcal L_2(H)}^2\\
&\le C\sup_{i\in\mathbb N_+}
\bigl|\tilde\sigma_i(x)\lambda_i^r\bigr|^2\le
C\bigl(1+\|x\|_{\dot H^{-1/2}}^{2\gamma}\bigr)
\sup_{i\in\mathbb N_+}(i^{2r-1})^2\le
C_r\bigl(1+\|x\|_{\dot H^{2r}}^{2\gamma}\bigr),
\end{align*}
which verifies Assumption \ref{assum:nonlinearity-B} (i). Moreover, for \(x,y\in H\),
\[
D\widehat B(x).y\, e_i
=
D\tilde\sigma_i(x).y\, e_i,
\qquad i\in\mathbb N_+,
\]
where
\begin{align*}
&D\tilde\sigma_i(x). y
=\frac{1}{i+1}\Bigg[
2V'(\|x\|_{\dot H^{-1/2}}^2)
\big\langle (-A)^{-1/4}x,(-A)^{-1/4}y\big\rangle_H\times
\left(
\sqrt{1+\|x\|_{\dot H^{-1/2}}^2}
-
\sqrt{1+\|x\|_{\dot H^{-1/2}}^{2\gamma}}
\right)\\
&+
V(\|x\|_{\dot H^{-1/2}}^2)
\frac{\big\langle (-A)^{-1/4}x,(-A)^{-1/4}y\big\rangle_H}
{\sqrt{1+\|x\|_{\dot H^{-1/2}}^2}}+
\bigl(1-V(\|x\|_{\dot H^{-1/2}}^2)\bigr)
\frac{\gamma\|x\|_{\dot H^{-1/2}}^{2(\gamma-1)}
\big\langle (-A)^{-1/4}x,(-A)^{-1/4}y\big\rangle_H}
{\sqrt{1+\|x\|_{\dot H^{-1/2}}^{2\gamma}}}
\Bigg].
\end{align*}
Using the boundedness of \(V\) and \(V'\), together with the fact that \(V'\) is supported in \([1,2]\), we obtain
\[
\sup_{x\in H}
|D\tilde\sigma_i(x).y|
\le\frac{C}{i+1}\|(-A)^{-1/4}y\|_H.
\]
Consequently, for any \(s_1\in[0,1/4]\),
\begin{align*}
\sup_{x\in H}
\|(-A)^{-s_1}DB(x).y\|_{\mathcal L_2^0(H)}^2
&=
\sup_{x\in H}
\|(-A)^{-s_1}D\widehat B(x).y\,Q^{1/2}\|_{\mathcal L_2(H)}^2\le
\|Q^{1/2}\|_{\mathcal L_2(H)}^2
\sup_{x\in H}
\|(-A)^{-s_1}D\widehat B(x).y\|_{\mathcal L(H)}^2\\
&\le
C_{s_1}
\sup_{i\in\mathbb N_+}
\lambda_i^{-2s_1}
\sup_{x\in H}|D\tilde\sigma_i(x).y|^2\le
C_{s_1}
\sup_{i\in\mathbb N_+}
\frac{\lambda_i^{-2s_1}}{(i+1)^2}
\|(-A)^{-1/4}y\|_H^2\\
&\le
C_{s_1}\|(-A)^{-s_1}y\|_H^2.
\end{align*}
The estimates required for the second- and third-order derivatives in Assumption \ref{assum:nonlinearity-B} (ii)--(iii) follow analogously.
\end{exam}

\begin{assum}\label{assum:nonlinearity-G}
Let \(G:H\times\chi\to H\) be
\(\mathcal B(H)\otimes\mathcal B(\chi)/\mathcal B(H)\)-measurable. Assume that, for \(\nu\)-a.e. \(z\in\chi\), the mapping \(x\mapsto G(x,z)\) is three times continuously Fr\'echet differentiable.  Moreover, for \(j=1,2,3\), assume that $
(x,z,y_1,\ldots,y_j)
\mapsto
D_x^jG(x,z).(y_1,\ldots,y_j)$ is jointly measurable. The following conditions hold.
\begin{itemize}
	\item[(i)] There exists a measurable function $\bar{c}\colon \chi\to \R_+$ such that $\sup_{x\in H}\|G(x,z)\|_{\dot{H}^{\frac12}}\le \bar{c}(z)$ for all $z\in\chi$, and for all  $L\ge 0$ there exists a positive constant $C_L$ such that 	
		$\int_{\chi}e^{L\bar{c}(z)} \big(\bar{c}(z)+(\bar{c}(z))^2\big)\,\nu(dz)\le C_{L}.$
	\item[(ii)]  For every $\varsigma_2\in[0,\frac14]$ and for all $p\ge 1$, there exists a positive constant $C_{\varsigma_2,p}$ such that
	\begin{align*}
	& \Big(\sup_{x\in H}\int_{\chi} \big\|(-A)^{-\varsigma_2}DG(x,z).y\big\|^p\,\nu(dz)\Big)^{\frac1p}\le C_{\varsigma_2,p}\|(-A)^{-\varsigma_2}y\|,\quad y\in H,\\
	&\Bigl(\sup_{x\in H}\int_{\chi} \big\|(-A)^{-\varsigma_2}D^2G(x,z).(y,h)\big\|^p\,\nu(dz)\Big)^{\frac1p}\le C_{\varsigma_2,p}\|(-A)^{-\varsigma_2}y\|\|(-A)^{-\varsigma_2}h\|,\quad y,h\in H.
	\end{align*}
   \item[(iii)]  For all $p\ge 1$, there exists a positive constant $C_{g,p}$ such that
\begin{align*}
\Bigl(\sup_{x\in H}\int_{\chi}\|D^3G(x,z).(y,h,k)\|^p\,\nu(dz)\Big)^{\frac1p} \le C_{g,p}\|y\|\|h\|\|k\|,\quad y,h,k\in H.
\end{align*}
\end{itemize}
\end{assum}

\begin{rem}
The sublinear growth condition on coefficient $B$ in Assumption \ref{assum:nonlinearity-B} (i) and the restriction imposed on the coefficient $G$ in 
Assumption \ref{assum:nonlinearity-G} (i) are required to obtain the exponential integrability of the solution process in the subsequent analysis; see, for example, Lemmas \ref{lem:momentbounds-fulldiscre} and \ref{lem:exponential-semidiscrete}. 

Notice that Assumption \ref{assum:nonlinearity-G}  (i) does not require $\nu(\chi)<\infty$ and therefore allows infinite-activity L\'evy noise. For example, let $\chi=H\setminus\{0\}$ and, for some $\alpha\in(0,1)$, define a L\'evy measure supported on the truncated ray $\{re_1:0<r<1\}$, where $e_1$ is the first eigenvector of $-A$, by  $\nu(B):=\int_0^1\mathbf 1_B(re_1)\,r^{-1-\alpha}\,dr$, $B\in\mathcal B(\chi)$. Then we derive $\nu(\chi)=\infty$ and $\int_{\chi}(1\wedge\|z\|^2)\,\nu(dz)<\infty.$ Moreover, the coefficient $G(x,z):=
\langle z,e_1\rangle_H
\tanh(\langle x,e_1\rangle_H)e_1 $ satisfies Assumption \ref{assum:nonlinearity-G}  with $\bar c(z)=\lambda_1^{1/4}\|z\|_H$. Indeed, for every $L\ge0$,
\[
\int_{\chi}
e^{L\bar c(z)}
\bigl(\bar c(z)+\bar c(z)^2\bigr)\nu(dz)
=
\int_0^1e^{L\lambda_1^{1/4}r}\bigl(\lambda_1^{1/4}r^{-\alpha}+\lambda_1^{1/2}r^{1-\alpha}\bigr)\,dr
<\infty,
\]
and the required bounds on the first three derivatives of $G$ follow from the boundedness of the first three derivatives of $\tanh$.

\end{rem}

The well-posedness of the problem \eqref{SPDE:MultiplicativeNoise}, along with related properties of the solution process, can be established under the above conditions on the coefficients.

\begin{propo}[{\cite[Theorem 2.4, Remark 2.5]{CCM-2010}}]\label{prop:wellposedness}
Let Assumptions \ref{assum:nonlinearity-F}--\ref{assum:nonlinearity-G} hold and let $X(0)\in L^2(\Omega,\mathcal F_0,\mathbb P;H)$. Then there exists a unique mild solution $\{X(t)\}_{t\in[0,T]}$ to \eqref{SPDE:MultiplicativeNoise} such that $\E\big[\sup_{t\in[0,T]}\|X(t)\|^2\big]<\infty$. Moreover, the solution $X$ has a c\`{a}dl\`{a}g modification and is a time-homogeneous Markov process.
\end{propo}

\subsection{Malliavin calculus}\label{subsec2}

In this subsection, we gather the concepts and results of Wiener--Poisson Malliavin calculus in Hilbert spaces that will be used in the weak convergence analysis, and refer to \cite{DE2018,GBF2009} and references therein for a comprehensive treatment.

For \(\mathbf B_k=(B_1,\ldots,B_k)\) with \(k\in\mathbb N_+\) and \(B_1,\ldots,B_k\in\hat\Xi_T\), we set $\mathbf N(\mathbf B_k)
:=
\bigl(N(B_1),\ldots,N(B_k)\bigr),$
where the sets \(B_1,\ldots,B_k\) are always assumed to be pairwise disjoint.

Let $I_Q:=\{i\in\mathbb N_+:q_i>0\}$ and $g_i:=Q^{1/2}e_i=\sqrt{q_i}\,e_i,~ i\in I_Q$.  \(\{g_i\}_{i\in I_Q}\) is an orthonormal basis of \(U_0=Q^{1/2}(H)\). For
\(\psi\in L^2([0,T];U_0)\), define the Hilbert--Schmidt functional
\[
\ell_{\psi(t)}:U_0\to\mathbb R,
\qquad
\ell_{\psi(t)}u:=\langle\psi(t),u\rangle_{U_0},
\]
and set
\[
W^Q(\psi)
:=
\int_0^T \ell_{\psi(t)}\,dW^Q(t)
:=
\sum_{i\in I_Q}
\int_0^T
\langle\psi(t),g_i\rangle_{U_0}\,d\beta_i(t).
\]
The series converges in \(L^2(\Omega;\mathbb R)\) and
$
\mathbb E[|W^Q(\psi)|^2]
=
\int_0^T\|\psi(t)\|_{U_0}^2\,dt.
$

We denote by \(\mathcal S_H\) the class of \(H\)-valued smooth cylindrical random variables of the form
\begin{align*}
\IS_{H}&:=\Big\{F_{H}=\sum_{r=1}^Rf_r(W^Q(\psi_1),W^Q(\psi_2),\ldots,W^Q(\psi_{\ell}),\mathbf{N}(\mathbf B_k))h_r,~ R,\ell,k\in \N_+,\\
&~~\quad\quad f_r\in \IC^\infty_{c,\mathrm{fin}}(\R^{\ell}\times\N^k;\R), ~\psi_i\in L^2([0,T];U_0), ~1\le i\le \ell, ~B_1,\ldots,B_k\in\hat{\Xi}{_T}, ~h_1,\ldots,h_R\in H\Big\}.
\end{align*}
 Here \(C_{c,\mathrm{fin}}^\infty(\mathbb R^\ell\times\mathbb N^k;\mathbb R)\) denotes the class of functions that are smooth with compact support in the \(\mathbb R^\ell\)-variables and have finite support in the discrete variable \(\mathbf n\in\mathbb N^k\). More precisely, for each \(\mathbf n\in\mathbb N^k\), $x\mapsto f(x,\mathbf n)$ belongs to \(C_c^\infty(\mathbb R^\ell;\mathbb R)\), and \(f(\cdot,\mathbf n)\equiv0\) for all but finitely many \(\mathbf n\in\mathbb N^k\).  Let $L^2_{\mathcal F_T}(\Omega;H)$ be the space of $\mathcal F_T$-measurable square integrable random variables.  It is known that $\IS_{H}$ is dense in $L^2_{\mathcal F_T}(\Omega;H)$ (cf. \cite[Lemma 3.1]{JJ-2008}).

 For \(F_H\in\mathcal S_H\) of the form
$F_H
=
\sum_{r=1}^R
f_r\bigl(
W^Q(\psi_1),\ldots,W^Q(\psi_\ell),
\mathbf N(\mathbf B_k)
\bigr)h_r,
$
the Wiener--Malliavin derivative is defined by
\[
\ID_t^W F_H
:=
\sum_{r=1}^R\sum_{i=1}^{\ell}
\partial_i f_r\bigl(
W^Q(\psi_1),\ldots,W^Q(\psi_\ell),
\mathbf N(\mathbf B_k)
\bigr)
\bigl(\psi_i(t)\otimes h_r\bigr),
\qquad t\in[0,T],
\]
where, for \(\psi_i(t)\in U_0\) and \(h_r\in H\), $(\psi_i(t)\otimes h_r)u
:=
\langle\psi_i(t),u\rangle_{U_0}h_r$ for $u\in U_0$. Thus $\ID_t^W F_H\in\mathcal L_2(U_0;H)
=\mathcal L_2^0(H)$. We write $\ID_t^{W,u}F_H:=\ID_t^W F_H.u$ for \(u\in U_0\). For \(dt\otimes\nu(dz)\)-a.e. \((t,z)\in\chi_T\), the Poisson--Malliavin derivative is given by
\[
\begin{aligned}
\ID_{(t,z)}^\mathbf{N} F_H
:=
\sum_{r=1}^R\sum_{i=1}^k
&\Big[
f_r\bigl(
W^Q(\psi_1),\ldots,W^Q(\psi_\ell),
\mathbf N(\mathbf B_k)+\mathbf e_i^{(k)}
\bigr)\\
&\quad-
f_r\bigl(
W^Q(\psi_1),\ldots,W^Q(\psi_\ell),
\mathbf N(\mathbf B_k)
\bigr)
\Big]
h_r\,\mathbf 1_{B_i}(t,z),
\end{aligned}
\]
where \(\mathbf e_i^{(k)}\) denotes the \(i\)-th canonical basis vector of \(\mathbb R^k\).

The operators $\ID^W$ and $\ID^\mathbf{N}$ are closable from $\IS_{H}$ into $L^2(\Omega\times[0,T];\IL_2^0(H))$ and into $L^2(\Omega\times\chi_T;H)$, respectively (cf. \cite[Proposition 4.4]{Kruse-2014} and \cite[Theorem 5.6]{SJ_2013}). Thus their domain can be extended to the Hilbert space $\D^{1,2}:=\overline{\mathcal S}_H^{\|\cdot\|_{\mathbb{D}^{1,2}}}$ with norm 
\begin{align*}
\|F_H\|_{\mathbb D^{1,2}}^2
:={}&
\mathbb E[\|F_H\|_H^2]
+
\mathbb E\Big[\int_0^T
\|D_t^W F_H\|_{\mathcal L_2^0(H)}^2\,dt\Big]+
\mathbb E\Big[\int_0^T\int_{\chi}
\|D_{(t,z)}^N F_H\|_H^2
\,\nu(dz)\,dt\Big].
\end{align*} 
By closability,  \(\ID^W\) and \(\ID^\mathbf{N}\) extend uniquely to $\ID^W\colon \mathbb D^{1,2}\to L^2(\Omega\times[0,T];\IL_2^0(H))$ and $\ID^\mathbf{N}\colon \mathbb D^{1,2}\to L^2(\Omega\times\chi_T;H)$, respectively.

According to \cite[Lemma 4.7]{Kruse-2014} and \cite[Proposition 2.6]{AT2018}, together with the independence of the Wiener process and the compensated Poisson process, the chain rule is valid: for $F_H\in\D^{1,2}$ and $\varphi\in\IC_b^1(H;\tilde{H})$ with another separable Hilbert space $\tilde{H}$, it holds $\varphi(F_H)\in \D^{1,2}(\tilde{H})$ and 
\begin{align}\label{eq:chainrule}
 \ID^{W}\varphi(F_H)=D\varphi(F_H).\ID^{W} F_H,  \quad \ID^\mathbf{N}\varphi(F_H)=\varphi(F_H+\ID^\mathbf{N}F_H)-\varphi(F_H).
\end{align}
For any \(A\in\hat{\Xi}_T\), let $\mathcal F_A^{\mathbf{N}}
:=
\sigma\bigl\{
\mathbf{N}(B): B\in\hat{\Xi}_T,\ B\subset A
\bigr\}\vee\mathcal N,$ where \(\mathcal N\) denotes the collection of \(\mathbb P\)-null sets.
Thus, \(\mathcal F_A^{\mathbf{N}}\) contains the information generated by the Poisson random measure restricted to \(A\). If \(F_H\) is \(\mathcal F_t\)-measurable and \(G_H\) is \(\mathcal F_A^{\mathbf{N}}\)-measurable for \(A\in\hat{\Xi}_T\), then
\begin{align}
 \ID^W_sF_H=0,&\quad  \mathbb P\otimes ds\text{-a.s.} ~\mbox{on}~ \Omega\times (t,T],\label{eq:zero-W}\\
 \ID^\mathbf{N}_{(s,z)}G_H=0, &\quad\mathbb P\otimes ds\otimes\nu(dz)\text{-a.s.} ~\mbox{on}~ \Omega\times A^c.\label{eq:zero-N}
\end{align}
We refer to \cite{GBF2009,SJ_2013,AF2019} for more details, including the following product rule. 

\begin{propo}[\textbf{product rule}]\label{prop:product}
For all $F_H, G_H\in \IS_H$ and $\tilde{u}\in U_0$, it holds that $\langle F_H,G_H\rangle\in\IS_\R$ and  
\begin{align*}
\ID^{W,\tilde{u}}(\langle F_H,G_H\rangle)&=\langle\ID^{W,\tilde{u}} F_H,G_H\rangle+\langle F_H, \ID^{W,\tilde{u}} G_H\rangle,\\
\ID^{\mathbf{N}}(\langle F_H,G_H\rangle)&=\langle\ID^{\mathbf{N}} F_H,G_H\rangle+\langle F_H, \ID^{\mathbf{N}} G_H\rangle+\langle\ID^{\mathbf{N}}F_H,\ID^{\mathbf{N}}G_H\rangle.
\end{align*}
\end{propo}

In addition, owing to the independence of the Wiener process and the Poisson random measure, the corresponding Malliavin duality formulas extend componentwise to the joint Wiener--Poisson probability space. Indeed, one may regard either noise component as an independent random parameter and apply the duality formula with respect to the other one. Consequently, by \cite[Proposition 4.3]{Kruse-2014} and  \cite[Proposition 5.5]{SJ_2013}, we conclude the following duality formula.

\begin{propo}\label{prop:dualityformula}
Let $F_H\in \D^{1,2}$, and let $\Phi\in L^2(\Omega\times[0,T];\IL_2^0(H))$ and $\phi\in L^2(\Omega\times\chi_T;H)$ be  predictable processes. Then 
\begin{align}
\E\Big[\Big\langle F_H, \int_0^T\Phi(t)dW^Q(t)\Big\rangle_H\Big]&=\E\big[\langle \ID^WF_H,\Phi\rangle_{L^2([0,T];\IL_2^0(H))}\big], \label{eq:duality-W}\\
\E\Big[\Big\langle F_H, \int_0^T\int_\chi\phi(t,z)\tilde{\mathbf{N}}(dz,dt)\Big\rangle_H\Big]&=\E\big[\langle \ID^\mathbf{N}F_H,\phi\rangle_{L^2(\chi_T;H)}\big].	 \label{eq:duality-N}
\end{align}
\end{propo}

As a consequence of the duality formula, we obtain the following conclusion.
\begin{propo}[{\cite[Lemma 2.1]{AD-2011}}]\label{prop:Malliavinderivative-W}
Let $F_H\in \D^{1,2}$, $\varphi\in\IC^2_b(H;\R)$, and $\Phi\in L^2(\Omega\times[0,T];\IL_2^0(H))$ be an adapted process. Then 
\begin{align*}
\E\Big[\Big\langle D\varphi(F_H), \int_0^T\Phi(t)dW^Q(t)\Big\rangle_H\Big]&=\sum_{k\in\N_+}\E\Big[\int_0^T D^2\varphi(F_H).\big(\mathcal D^{W, Q^{\frac12}e_k}_tF_H,  \Phi(t)Q^{\frac12}e_k\big)\,dt\Big].
\end{align*}
\end{propo}

As mentioned in \cite[Lemma 2.1]{AD-2011}, the above result remains valid for $\varphi\in C_{\mathfrak{p}}^2(H;\mathbb R)$, provided that the required integrability conditions are satisfied. This follows by approximating $\varphi$ with bounded smooth functions and passing to the limit by standard dominated convergence arguments.


\section{Numerical methods and properties of numerical solutions}\label{sec:3}

This section introduces the numerical methods studied in this paper and establishes their key properties, including the moment boundedness, exponential integrability, and estimates of the Malliavin derivatives of the numerical solutions. 

\subsection{Numerical methods}\label{subsec:3.1}

In the spatial direction, the spectral Galerkin method is applied to obtain the spatial semi-discrete solution $\{X^N(t)\}_{t\in[0,T]}$ with discretization parameter $N\in \N_+$. In the temporal direction, the exponential Euler method is further employed to obtain the spatio-temporal fully discrete solution $\{\bar{X}^N_m\}_{m\in\{0,1,\ldots,M\}}$ with some $M\in \N_+$ and time step size $\tau=T/M$.

\subsubsection{Spatial semi-discretization}
For arbitrary $N\in \N_+$,  let $P_N$ be the orthogonal projection operator of $H$ onto $H_N:=\mbox{span}\{e_1,e_2,...,e_N\}$ defined by 
$P_Nx:=\sum_{j=1}^N\langle x,e_j\rangle e_j$ for all $x\in H$. The spectral Galerkin method applied to  \eqref{SPDE:MultiplicativeNoise} satisfies 
\begin{align}\label{eq:spectralgalerkin}
dX^N(t)&=AX^N(t)dt+F_N(X^N(t))dt+P_NB(X^N(t))dW^Q(t)\nonumber\\
&\quad+\int_{\chi}G_N(X^N(t),z)\,\tilde{\mathbf{N}}(dz,dt),\quad t\in(0,T]
\end{align}
with initial value $X^N(0)=P_Nx_0$.  The unique mild solution to \eqref{eq:spectralgalerkin} is given by
\begin{align*}
X^N(t,P_Nx_0)&=E(t)P_Nx_0+\int_0^tE(t-s)F_N(X^N(s))ds+\int_0^tE(t-s)P_NB(X^N(s))dW^Q(s)\\
&~~~+\int_0^t\int_\chi  E(t-s)G_N(X^N(s),z)\tilde{\mathbf{N}}(dz,ds).
\end{align*}

\subsubsection{Spatio-temporal full-discretization}
Applying further the exponential Euler method to \eqref{eq:spectralgalerkin}, we obtain the spatio-temporal full-discretization  $\{\bar{X}^{N}_m\}_{m\in\{0,1,\ldots,M\}}$ satisfying $\bar{X}^N_0=X^N(0)$ and 
\begin{align}\label{eq:one-step}
\bar{X}^{N}_{m+1}&=E(\tau)\Big(\bar{X}^{N}_m+\tau F_N(\bar{X}^{N}_m)+P_NB(\bar{X}^{N}_m)\Delta W^Q_{m}+\int_{t_m}^{t_{m+1}}\int_\chi G_N(\bar{X}^{N}_m,z)\tilde{\mathbf{N}}(dz,ds)\Big)
\end{align}
 with $\Delta W^Q_{m}=W^Q(t_{m+1})-W^Q(t_m)$ and $m=0,1,\ldots,M-1$.  Define the continuous-time interpolation $\{\tilde{X}^{N}(t)\}_{t\in[0,T]}$ by
\begin{align}\label{eq:continuoustime}
\tilde{X}^{N}(t)&=E(t-t_m)\Big(\bar{X}^{N}_m+(t-t_m)F_N(\bar{X}^{N}_m)+P_NB(\bar{X}^{N}_m)(W^Q(t)\!-\!W^Q(t_m))\nonumber\\
&~~~+\!\int_{t_m}^{t}\int_\chi G_N(\bar{X}^{N}_m,z)\tilde{\mathbf{N}}(dz,ds)\Big),\quad t\in[t_m,t_{m+1})
\end{align}
for $m\in\{0,1,\ldots,M-1\}$, and set
$\tilde X^N(T):=\bar X_M^N$. Then it holds that $\tilde{X}^{N}(t_i)=\bar{X}^{N }_i$ for $i\in \{0,1,\ldots,M\}$ and  
\begin{align*}
\!\!\!d\tilde{X}^{N}(t)=A\tilde{X}^{N }(t)dt+E(t-\lfloor t\rfloor)\Big(F_N(\bar{X}^{N}_{\lfloor t\rfloor/\tau})dt+P_NB(\bar{X}^{N }_{\lfloor t\rfloor/\tau})dW^Q(t)+\int_\chi G_N(\bar{X}^{N}_{\lfloor t\rfloor/\tau},z)\tilde{\mathbf{N}}(dz,dt)\Big)
\end{align*}
for $t\in[t_m,t_{m+1})$ with $m\in\{0,1,\ldots,M-1\}$, where $\lfloor t\rfloor:=\max\{t_m\colon t_m\le t, ~m\in\N\}$.

In the sequel, we fix $\epsilon \in (0, 1]$ and occasionally impose the following restriction on the time step size $\tau$ for technical reasons:  \begin{align*}
\tau<\Big(\frac{e^{-\frac{\sqrt{\ln 2}}{4}}}{4(1\vee M_F)}\Big)^4\wedge \big(\lambda_1^{\frac{1-\epsilon}{2}} e^{-\frac{\sqrt{\ln 2}}{2}}\big),
\end{align*}
together with the CFL condition $\tau\lambda_N<\frac{\sqrt{\ln 2}}{4}$, or equivalently, $\tau<\frac{\sqrt{\ln 2}}{4\lambda_N}$. For simplicity, we denote 
\begin{align}\label{def:tau_0}
\tau_0:=\Big(\frac{e^{-\frac{\sqrt{\ln 2}}{4}}}{4(1\vee M_F)}\Big)^4\wedge \big(\lambda_1^{\frac{1-\epsilon}{2}} e^{-\frac{\sqrt{\ln 2}}{2}}\big), \quad \tau_N:=\frac{\sqrt{\ln 2}}{4\lambda_N}\wedge\tau_0, \quad  M_F:=C_E(L_F+L_0),
\end{align}
where $C_E>0$ is a constant satisfying $\|E(t)(-A)^{\frac{1-\epsilon}{4}}\|_{\IL(H)}\le C_Et^{-\frac{1-\epsilon}{4}}$ for $t>0$, and $L_F$ and $L_0$ are given in  Assumption \ref{assum:nonlinearity-F}. This domain of the time step size $\tau$ will be specified for each result in what follows.

\subsection{ Exponential integrability and moment estimates}

We establish in this subsection the exponential integrability and moment estimates of numerical solutions given in \eqref{eq:spectralgalerkin} and \eqref{eq:continuoustime}.   The exponential integrability of the numerical solutions will be used to  establish weak convergence results for test functions $\varphi \in \IC^3_{\mathfrak{p}}(H;\R)$, whose derivatives  
 exhibit polynomial growth.
We also note that this type of exponential integrability admits many further applications (cf. \cite{JJ-2018,CD-2023,JJZ-2017,CTJF-2025}).

\begin{lemma}\label{lem:momentbounds-fulldiscre}
Let Assumptions \ref{assum:nonlinearity-F}--\ref{assum:nonlinearity-G} hold. For any fixed $\epsilon\in(0,1]$, $r_0>0$, and
$x_0\in\dot H^{\frac{1-\epsilon}{2}}$, we define $U_{r_0}\colon\dot{H}^{\frac{1-\epsilon}{2}}\to \R_+$ by $U_{r_0}(y):=r_0(1+\|(-A)^{\frac{1-\epsilon}{4}}y\|^2)^{1/2}$ for $y\in\dot H^{\frac{1-\epsilon}{2}}$. Then there exist positive constants $\varsigma_4$ and $C$ dependent on $r_0$ such that  for all $N\in\N_+$ and $\tau\in(0,\tau_N)$ with $\tau_N$ given in \eqref{def:tau_0},
\begin{align}\label{momentbounds-fulldiscre}
\sup_{t\in[0,T]} \E\Big[\exp\big(\frac{U_{r_0}(\tilde{X}^{N}(t))}{e^{\varsigma_4t}}\big)\Big]\le e^{CT}e^{U_{r_0}(x_0)}.
\end{align}
\end{lemma}

The proof of Lemma  \ref{lem:momentbounds-fulldiscre} is rather technical and is deferred to the Appendix.

\begin{cor}\label{cor:exponentialintegrability}
Let Assumptions \ref{assum:nonlinearity-F}--\ref{assum:nonlinearity-G} be satisfied. Then for arbitrary fixed constants $M_0>0$, $k\ge 1$, $\epsilon\in(0,1]$ and for all $\beta_2\in[0, 1)$ and $x_0\in\dot{H}^{\frac{1-\epsilon}{2}}$, there exist positive constants $L_3$ and $C_{L_3,\beta_2,T,M_0,k}$ such that for all $N\in\N_+$ and all $\tau\in(0,\tau_N)$ with $\tau_N$ given in \eqref{def:tau_0}, 
\begin{align*}
\sup_{t\in[0,T]}\Big(\E\big[\exp\big(M_0k\|(-A)^{\frac{1-\epsilon}{4}}\tilde{X}^{N}(t)\|^{\beta_2}\big)\big]\Big)^{\frac{1}{k}}\le \exp\big(C_{L_3, \beta_2,T,M_0,k}+\frac{1}{4ke^{L_3T}}\|(-A)^{\frac{1-\epsilon}{4}}x_0\|\big).
\end{align*}	
In particular, for all $q\ge 1$, $\beta_1\in[0,\frac14)$ and for all $x_0\in\dot{H}^{2\beta_1}$,  there exists a constant $C_{L_3,q,T}>0$ such that for all $N\in\N_+$ and $\tau\in(0,\tau_N)$, 
\[
\sup_{t\in[0,T]}\E\big[\big\|(-A)^{\beta_1}\tilde{X}^{N}(t)\big\|^{2q}\big]\le C_{L_3,q,T}\exp\Big(\frac{\|(-A)^{\beta_1}x_0\|}{4e^{L_3T}}\Big).
\]
\end{cor}

\begin{proof}
The proof is based on the application of Lemma \ref{lem:momentbounds-fulldiscre}. From the proof of Lemma \ref{lem:momentbounds-fulldiscre} in the Appendix, the constant $\varsigma_4:=\varsigma_4(r_0)$ given in \eqref{def:varsigma4} is an increasing affine function of $r_0$. Throughout the remainder of the paper, we choose $L_3>0$ sufficiently large so that $L_3\ge \varsigma_4(r_0)$ for all values of $r_0$ used below. In this case, \eqref{momentbounds-fulldiscre} yields that 
\begin{align*}
	\sup_{t\in[0,T]} \E\Big[\exp\big(\frac{U_{r_0}(\tilde{X}^{N}(t))}{e^{L_3t}}\big)\Big]\le \sup_{t\in[0,T]} \E\Big[\exp\big(\frac{U_{r_0}(\tilde{X}^{N}(t))}{e^{\varsigma_4(r_0)t}}\big)\Big] \le e^{CT}e^{U_{r_0}(x_0)}.
\end{align*}
Below we take $r_0:=\frac{1}{4e^{L_3T}}$ which belongs to $(0,1]$, and write $U$ for $U_{r_0}$ since $r_0$ is fixed. Applying Lemma \ref{lem:momentbounds-fulldiscre} with such $r_0$ and setting $\eta:=\frac{1}{4M_0k}e^{-2L_3T}$,
we derive  for all  $N\in\N_+$ and $\beta_2\in[0, 1)$ that
\begin{align*}
&~~~\sup_{t\in[0,T]}\E\big[\exp\big(M_0k\|(-A)^{\frac{1-\epsilon}{4}}\tilde{X}^{N}(t)\|^{\beta_2}\big)\big]\le \sup_{t\in[0,T]}\E\big[\exp\big(M_0k\eta\|(-A)^{\frac{1-\epsilon}{4}}\tilde{X}^{N}(t)\|+C_{\eta,\beta_2}\big)\big]\\
&\le e^{C_{\eta,\beta_2}}\times \sup_{t\in[0,T]}\E\Big[\exp\big(\frac{U( \tilde{X}^{N}(t))}{e^{L_3t}}\big)\Big]\le e^{C_{\eta,\beta_2}+CT}\E[e^{U(x_0)}]= e^{C_{\eta,\beta_2,T}} \exp\Big(\frac{\sqrt{1+\|(-A)^{\frac{1-\epsilon}{4}}x_0\|^2}}{4e^{L_3T}}\Big).
\end{align*}	
In particular, for all $q\ge 1$ and $\beta_1\in[0,\frac14)$, there exist positive constants $C_q, C_{L_3,q,T}$ such that
\begin{align*}
\sup_{t\in[0,T]}\E\big[\big\|(-A)^{\beta_1}\tilde{X}^{N}(t)\big\|^{2q}\big]\le C_q \sup_{t\in[0,T]}\E\big[\exp\big(\frac12\|(-A)^{\beta_1}\tilde{X}^{N}(t)\|^{\frac12}\big)\big]\le C_{L_3,q,T}\exp\Big(\frac{\|(-A)^{\beta_1}x_0\|}{4e^{L_3T}}\Big),
\end{align*}  
which completes the proof.
\end{proof}

Similar to Lemma \ref{lem:momentbounds-fulldiscre} and Corollary \ref{cor:exponentialintegrability}, the corresponding estimates for $X^N(t)$ are stated in the following lemma.

\begin{lemma}\label{lem:exponential-semidiscrete}
Let Assumptions \ref{assum:nonlinearity-F}--\ref{assum:nonlinearity-G} be satisfied. For any fixed $\epsilon\in(0,1]$, $r_0>0$, and $x_0\in\dot H^{\frac{1-\epsilon}{2}}$, there exist positive constants $\varsigma_5$ and $C$  dependent on $r_0$  such that for all $N\in\N_+$,
\begin{align}\label{momentbounds-semidiscre}
\sup_{t\in[0,T]} \E\Big[\exp\big(\frac{U_{r_0}(X^{N}(t))}{e^{\varsigma_5t}}\big)\Big]\le e^{CT}e^{U_{r_0}(x_0)}.
\end{align}
Moreover, for arbitrary fixed constants $M_0>0$, $k\ge 1$, $\epsilon\in(0,1]$ and for all $\beta_4\in[0, 1)$, there exist positive constants $L_4$ and $C_{L_4, \beta_4,T}$ such that for all $N\in\N_+$ and $x_0\in\dot{H}^{\frac{1-\epsilon}{2}}$
\begin{align}\label{eq:exponentialintegrability-exact1}
\sup_{t\in[0,T]}\Big(\E\big[\exp\big(M_0k\|(-A)^{\frac{1-\epsilon}{4}}X^{N}(t)\|^{\beta_4}\big)\big]\Big)^{\frac{1}{k}}\le \exp\big(C_{L_4, \beta_4,T}+\frac{1}{4ke^{L_4T}}\|(-A)^{\frac{1-\epsilon}{4}}x_0\|\big).
\end{align}	
It implies that for all $p\ge 2$ and $\beta_3\in[0,\frac14)$,  there exists $C_{L_4,p,T}>0$ such that for all $N\in\N_+$ and $x_0\in\dot{H}^{2\beta_3}$
\begin{align}\label{eq:exponentialintegrability-exact}
&\sup_{t\in[0,T]}\|(-A)^{\beta_3}X^N(t)\|_{L^p(\Omega;H)}\le C_{L_4,p,T}\exp\Big(\frac{\|(-A)^{\beta_3}x_0\|}{4pe^{L_4T}}\Big).
\end{align} 
\end{lemma}

The proof of Lemma \ref{lem:exponential-semidiscrete} can be obtained through a slight modification of the proofs of Lemma \ref{lem:momentbounds-fulldiscre} and Corollary \ref{cor:exponentialintegrability}, and hence is omitted here for brevity. We remark that   combining  the estimates for the stochastic integral with respect to $W^Q$ (cf. \cite[Theorem 4.36]{GJ-2014}) and the one with respect to $\tilde{N}$ (cf. \cite[Proposition 3.3]{CCM-2010}),  we can obtain directly the uniform  moment bounds of the fully discrete numerical solution  in state space $H$ under Assumptions \ref{assum:nonlinearity-F}--\ref{assum:nonlinearity-G}. More precisely, for any $K\in\mathbb N_+$, there exists a constant
$C_{K,T}>0$, independent of $N$ and $\tau$, such that 
\begin{align}\label{eq:estimate-X}
\sup_{m\in\{0,1,\ldots,M\}}\|\bar{X}^{N}_m\|_{L^{2K}(\Omega;H)}\le C_{K,T}(1+\|x_0\|).   
\end{align}

\subsection{Estimates of Malliavin derivatives}

We next investigate the estimates for the  Wiener–Malliavin and Poisson–Malliavin derivatives of the solution process $\{\tilde{X}^N(t)\}_{t\in[0,T]}$.

\begin{lemma}\label{lem:Malliavinderivative-B1}
Let Assumptions \ref{assum:nonlinearity-F}--\ref{assum:nonlinearity-G} be satisfied. Then 

(\romannumeral1) there exists a constant $C_T>0$ such that, for all
$N\in\mathbb N_+$, $\tau\in(0,1)$,
$\ell\in\{0,\ldots,M-1\}$, $\theta\in U_0$, and
$t\in(t_\ell,T]$, estimate
\begin{align}
\big(\E[\|\ID^{W,\theta}_s \tilde{X}^{N}(t)\|^{2}]\big)^{\frac{1}{2}}
	&\le \! C_{T}(\E[\|B(\bar{X}^{N}_\ell)\theta\|^{2}])^{\frac{1}{2}}\label{eq:Malliavinderivative-4}    
\end{align}
holds for a.e. $s\in(t_\ell,t_{\ell+1})\cap(0,t]$.

(\romannumeral2) for any $p\ge1$, there exists a constant $C_{p,T}>0$ such that, for all $N\in\mathbb N_+$, $\tau\in(0,1)$, $\ell\in\{0,\ldots,M-1\}$, and $t\in(t_\ell,T]$, estimate 
\begin{align}\label{ineq:DN1}
\big(\E[\| \ID^{\mathbf{N}}_{(s,v)}\tilde{X}^{N}(t)\|^{2p}]\big)^{\frac{1}{2p}}&\le C_{p,T}|\bar{c}(v)|
\end{align}
holds for
$(ds\otimes\nu(dv))$-a.e.
$(s,v)\in((t_\ell,t_{\ell+1})\cap(0,t])\times\chi$. Here  function $\bar{c}\colon \chi\to\R$ is defined by Assumption~\ref{assum:nonlinearity-G}.
\end{lemma}

\begin{proof}
 \textbf{(\romannumeral1)} By expression \eqref{eq:one-step} of solution $\{\bar{X}^N_n\}_{n\in\{\ell+1,\ldots, M\}}$, properties \eqref{eq:chainrule} and \eqref{eq:zero-W}, and Malliavin derivatives of stochastic integrals (cf. \cite[Proposition 4.16]{Kruse-2014}), we can deduce the Wiener--Malliavin differentiability of  $\{\bar{X}^N_n\}_{n\in\{\ell+1,\ldots,M\}}$. Consequently,  for all $n\in\{\ell+2,\ldots, M\}$, we have
\begin{align*}
\mathcal D^{W,\theta}_s \bar{X}^N_{n}&=E((n-\ell)\tau)P_NB(\bar{X}^{N}_\ell)\theta+\tau\sum_{j=\ell+1}^{n-1}E((n-j)\tau )P_NDF(\bar{X}^{N}_{j}).\mathcal D_s^{W,\theta}\bar{X}^N_{j}\\
&~~~\quad+\sum_{j=\ell+1}^{n-1}E((n-j)\tau )P_NDB(\bar{X}^{N}_{j}).\mathcal D_s^{W,\theta}\bar{X}^N_{j}\Delta W^Q_{j}\\
&~~~\quad +\sum_{j=\ell+1}^{n-1}E((n-j)\tau) \int_{t_{j}}^{t_{j+1}}\int_\chi P_NDG(\bar{X}^{N}_{j},z).\ID_s^{W,\theta}\bar{X}^N_{j}\,\tilde{\mathbf{N}}(dz,dt).
\end{align*}
In particular, it holds that $
\ID^{W,\theta}_s \bar{X}^{N}_{\ell+1}=E(\tau)P_NB(\bar{X}^{N}_\ell)\theta$ for $s\in (t_\ell,t_{\ell+1})$. Under Assumption \ref{assum:nonlinearity-F} (ii), Assumption \ref{assum:nonlinearity-B} (ii) with $\varsigma_1=0$, and Assumption \ref{assum:nonlinearity-G} (ii) with $\varsigma_2=0$, the It\^o isometries (cf. \cite[Proposition 4.20]{GJ-2014} and \cite[Proposition 2.16]{TE-2016}) yield
\begin{align*}
\E[\|\mathcal D^{W,\theta}_s \bar{X}^N_{n}\|^2]
&\le C\E\big[\big\|B(\bar{X}^{N}_\ell)\theta\big\|^2\big]+C_T\tau\sum_{j=\ell+1}^{n-1}\E\Big[\|DF(\bar{X}^{N}_{j}).\mathcal D^{W,\theta}_s \bar{X}^N_{j}\|^2+\|DB(\bar{X}^{N}_{j}).\mathcal D^{W,\theta}_s \bar{X}^N_{j}\|_{\IL_2^0(H)}^2\Big]\\
&~~~+C\tau\sum_{j=\ell+1}^{n-1}\E\Big[\int_{\chi}\|DG(\bar{X}^{N}_{j},z).\ID_s^{W,\theta}\bar{X}^N_{j}\|^2\,\nu(dz)\Big] \\
&\le C\E\big[\big\|B(\bar{X}^{N}_\ell)\theta\big\|^2\big]+C_T\tau\sum_{j=\ell+1}^{n-1} \E[\|\mathcal D^{W,\theta}_s \bar{X}^N_{j}\|^2].
\end{align*}
Then the discrete Gr\"onwall inequality gives that 
$
\E[\|\mathcal D^{W,\theta}_s \bar{X}^N_{n}\|^2]\le C_T\E[\|B(\bar{X}^{N}_\ell)\theta\|^{2}].
$
Similarly, for all $t\in[s,T]$, based on the expression \eqref{eq:continuoustime} of solution $\{\tilde{X}^N(t)\}_{t\in[s,T]}$, we derive
\begin{align*}
\E[\|\mathcal D^{W,\theta}_s \tilde{X}^N(t)\|^2]&\le \!C_T\Big(\E\big[\|\mathcal D^{W,\theta}_s \bar{X}^N_{\lfloor t\rfloor/\tau}\|^2\big]\!+\!(t-\lfloor t\rfloor)\E\big[\|\mathcal D^{W,\theta}_s \bar{X}^N_{\lfloor t\rfloor/\tau}\|^2\big]\!+\!\E[\|B(\bar{X}^{N}_\ell)\theta\|^{2}]\Big)\le \!C_T\E[\|B(\bar{X}^{N}_\ell)\theta\|^{2}].
\end{align*}

 \textbf{(\romannumeral2)} For $(ds\otimes\nu(dv))$-a.e.
$(s,v)\in(t_\ell,t_{\ell+1})\times\chi$, by properties \eqref{eq:chainrule} and \eqref{eq:zero-N}, together with the Malliavin derivative formula for stochastic integrals  (cf. \cite[Proposition 2.8]{AF2019}),  we deduce the Poisson--Malliavin differentiability of $\{\bar{X}^N_n\}_{n\in\{\ell+1,\ldots,M\}}$. Consequently, applying the Poisson--Malliavin derivative $\mathcal D^N_{(s,v)}$ to both sides of the recursion \eqref{eq:one-step}, we obtain, for all $n>\ell+1$, 
\begin{align*}
\ID^{\mathbf{N}}_{(s,v)}\bar{X}^{N}_n&=E((n-\ell)\tau)G_N(\bar{X}^{N}_{\ell},v)+\tau\sum_{j=\ell+1}^{n-1} E((n-j)\tau )(F_N(\bar{X}^{N}_{j}+\ID^{\mathbf{N}}_{(s,v)}\bar{X}^{N}_j)-F_N(\bar{X}^{N}_j))\\
&~~~+\sum_{j=\ell+1}^{n-1}E((n\!-\!j)\tau ) (P_NB(\bar{X}^{N}_{j}\!+\!\ID^{\mathbf{N}}_{(s,v)}\bar{X}^{N}_j)\!-\!P_NB(\bar{X}^{N}_j))\Delta W^Q_{j}\\
&~~~+\sum_{j=\ell+1}^{n-1}\int_{t_{j}}^{t_{j+1}}\!\!\int_\chi E((n\!-\!j)\tau )\big(G_N(\bar{X}^{N}_{j}+\ID^{\mathbf{N}}_{(s,v)}\bar{X}^{N}_{j},z)-G_N(\bar{X}^{N}_{j},z)\big)\,\tilde{\mathbf{N}}(dz,dt).
\end{align*}
 In particular,  we derive for $n=\ell+1$ that $\ID^{\mathbf{N}}_{(s,v)}\bar{X}^{N}_{\ell+1}=E(\tau)G_N(\bar{X}^{N}_{\ell},v)$. Under Assumptions \ref{assum:nonlinearity-F}--\ref{assum:nonlinearity-G}, using the H\"older and
Burkholder--Davis--Gundy
inequalities yields
\begin{align*}
\E[\| \ID^{\mathbf{N}}_{(s,v)}\bar{X}^{N}_n\|^{2p}]&\le C_{p}\E\big[\|G(\bar{X}^{N}_{\ell},v)\|^{2p}\big]+C_{p,T}\int_{t_{\ell+1}}^{t_{n}}\int_0^1\E\Big[\big\| DF(\bar{X}^{N}_{\lfloor t\rfloor/\tau}+\theta_1\ID^{\mathbf{N}}_{(s,v)}\bar{X}^{N}_{\lfloor t\rfloor/\tau}).\ID^{\mathbf{N}}_{(s,v)}\bar{X}^{N}_{\lfloor t\rfloor/\tau}\big\|^{2p}\\
&~~~+\big\| DB(\bar{X}^{N}_{\lfloor t\rfloor/\tau}+\theta_1\ID^{\mathbf{N}}_{(s,v)}\bar{X}^{N}_{\lfloor t\rfloor/\tau}).\ID^{\mathbf{N}}_{(s,v)}\bar{X}^{N}_{\lfloor t\rfloor/\tau}\big\|_{\IL_2^0(H)}^{2p}\Big]\,d\theta_1 dt\\
&+C_p\int_{t_{\ell+1}}^{t_{n}}\E\Big[\Big(\int_{\chi}\int_0^1\big\| DG(\bar{X}^{N}_{\lfloor t\rfloor/\tau}+\theta_1\ID^{\mathbf{N}}_{(s,v)}\bar{X}^{N}_{\lfloor t\rfloor/\tau},z).\ID^{\mathbf{N}}_{(s,v)}\bar{X}^{N}_{\lfloor t\rfloor/\tau}\big\|^{2}\,d\theta_1 \nu(dz)\Big)^p\\
&~~~+\int_{\chi}\int_0^1\big\| DG(\bar{X}^{N}_{\lfloor t\rfloor/\tau}+\theta_1\ID^{\mathbf{N}}_{(s,v)}\bar{X}^{N}_{\lfloor t\rfloor/\tau},z).\ID^{\mathbf{N}}_{(s,v)}\bar{X}^{N}_{\lfloor t\rfloor/\tau}\big\|^{2p}\,d\theta_1 \nu(dz)\Big]\,dt\\
&\le C_p|\bar{c}(v)|^{2p}+C_{p,T}\int_{t_{\ell+1}}^{t_{n}} \E\big[\|\ID^{\mathbf{N}}_{(s,v)}\bar{X}^{N}_{\lfloor t\rfloor/\tau}\|^{2p}\big]\,dt,
\end{align*}
and then we derive that $\E[\| \ID^{\mathbf{N}}_{(s,v)}\bar{X}^{N}_n\|^{2p}]\le C_{p,T}|\bar{c}(v)|^{2p}$ for all $n\in\{\ell+1,\ldots,M\}$.  Similarly, from \eqref{eq:continuoustime}, we obtain
\begin{align*}
\E\big[
\|\mathcal D^{\mathbf{N}}_{(s,v)}\tilde X^N(t)\|^{2p}
\big]
\le C_{p,T}\Big((1+t-\lfloor t\rfloor)
\E\big[
\|\mathcal D^{\mathbf{N}}_{(s,v)}
\bar X^N_{\lfloor t\rfloor/\tau}\|^{2p}
\big] +
\mathbf 1_{(\lfloor t\rfloor,t]}(s)|\bar c(v)|^{2p}
\Big)
\le C_{p,T}|\bar c(v)|^{2p}.
\end{align*}
The proof is completed.
\end{proof}

\section{Kolmogorov equations and regularity estimates}\label{sec:4}

In this section, we study the Kolmogorov equations associated to the spatial semi-discretization \eqref{eq:spectralgalerkin}. 
We establish new regularity estimates for the derivatives of their solutions, which are crucial for obtaining the weak convergence results (see Theorem \ref{the:weakconvergence-full}) for the numerical methods proposed in Section~\ref{sec:3}.

\subsection{Finite-dimensional setting}

For any $N\in\mathbb N_+$ and $\varphi\in \IB_b(H;\mathbb R)$, define $u^N(t,x):=\E[\varphi(X^{N,(0,x)}(t))]$ for any $t\in[0,T]$ and $x\in H_N$, where $X^{N,(0,x)}(t)$ denotes the solution of \eqref{eq:spectralgalerkin} starting from $X^N(0)=x$. For sufficiently regular test functions $\varphi$, the function $u^N$ is associated with the following Kolmogorov equation:
\begin{align}\label{eq:kolmogorovequation}
\frac{\partial u^N(t,x)}{\partial t}&=\mathbb{L^N}u^N(t,x):= \langle Du^N(t,x), Ax+F_N(x)\rangle +\frac12\mbox{Tr}\big((P_NB(x)Q^{\frac12})^{*}D^2u^N(t,x)(P_NB(x)Q^{\frac12})\big)\nonumber\\
&~~~+ \int_\chi\big(u^N(t,x+G_N(x,z))-u^N(t,x)-\langle Du^N(t,x), G_N(x,z)\rangle\big)\nu(dz),\quad t\in(0,T]
\end{align}
starting from $u^N(0,\cdot)=\varphi(\cdot)$, where $D$ and $D^2$ stand for the first and second order Fr\'echet derivatives with respect to $x$, respectively. If $\varphi\in \IC^2_{\mathfrak{p}}(H;\R)$, then for all $x,y,h_1,h_2\in H_N$ and $t\in[0,T]$, it holds that
\begin{align}
&Du^{N}(t,x).y=\E\big[D\varphi(X^{N}(t)).\eta^{N}_{y,x}(t)\big],\label{expression:uN-1}\\
&D^2 u^{N}(t,x).(h_1,h_2)=\E\big[D\varphi(X^{N}(t)). \xi^{N}_{h_1,h_2,x}(t)\big]+\E\big[D^2\varphi(X^{N}(t)).(\eta^{N}_{h_1,x}(t),\eta^{N}_{h_2,x}(t))\big]\label{expression:uN-2}
\end{align}
according to the chain rule, 
where $\eta^N_{y,x}(t):=DX^{N,(0,x)}(t).y$ solves the first order variational equation
 \begin{align*}
d\eta^{N }_{y,x}(t)=&\,A \eta^{N }_{y,x}(t)dt\!+\!P_NDF(X^{N }(t)).\eta^{N }_{y,x}(t)dt+P_NDB(X^{N }(t)).\eta^{N }_{y,x}(t)\,dW^Q(t)\\
&+\int_\chi P_NDG(X^{N }(t),z).\eta^{N }_{y,x}(t)\,\tilde{\mathbf{N}}(dz,dt)
\end{align*}	
with initial value $\eta^{N }_{y,x}(0)=y$, and $\xi^{N}_{h_1,h_2,x}(t):=D^2X^{N,(0,x)}(t).(h_1,h_2)$ satisfies $\xi^{N}_{h_1,h_2,x}(0)=0$ and
\begin{align*}
d\xi^{N}_{h_1,h_2,x}(t)=&\big(A\xi^{N}_{h_1,h_2,x}(t)+P_N(DF(X^{N}(t)).\xi^{N}_{h_1,h_2,x}(t)+D^2F(X^N(t)).(\eta^{N}_{h_1,x}(t),\eta^{N}_{h_2,x}(t)))\big)dt\\
&+P_N\big(DB(X^N(t)).\xi^{N}_{h_1,h_2,x}(t)\!+\!D^2B(X^N(t)).(\eta^{N}_{h_1,x}(t),\eta^{N}_{h_2,x}(t))\big)dW^Q(t)\\
&+\int_\chi P_N\big(DG(X^{N }(t),z).\xi^{N}_{h_1,h_2,x}(t)+D^2G(X^{N }(t),z).(\eta^{N}_{h_1,x}(t),\eta^{N}_{h_2,x}(t))\big)\,\tilde{\mathbf{N}}(dz,dt).
\end{align*}

For $\varphi\in C_{\mathfrak{p}}^2(H;\mathbb R)$, the differentiability of the finite-dimensional stochastic flow with respect to its initial value, together with the corresponding moment estimates, yields \eqref{expression:uN-1}--\eqref{expression:uN-2}. These representations provide the spatial regularity required below. In the following, we establish refined regularity estimates for the derivatives of $u^N$, which are the key ingredients in the weak convergence analysis.


\begin{propo}
 \label{cor:firstderivative2}
Let Assumptions \ref{assum:nonlinearity-F}--\ref{assum:nonlinearity-G} hold and $\varphi\in\IC^2_{\mathfrak{p}}(H;\R)$. Then for all $t\in(0,T]$, $\epsilon\in(0,\frac12\wedge\frac{1}{(\lambda_1\wedge 1)^{-3/4}-1}\wedge(\frac{1}{\alpha_0}-1))$, $\epsilon_1\in(0,\epsilon]$, and $\gamma_0\in[0,\frac{1-\epsilon}{4}+\frac{1-\epsilon_1}{2(1+\epsilon_1)(1+2\epsilon_1)})$, there exist positive constants $L_3, C_{L_3,\epsilon,T}, $ and $ C_{\mathfrak{p},T}$ such that for all $N\in\N_+$ and $x\in  H_N$
\begin{align}
\|(-A)^{\gamma_0}D u^N(t,x)\|&\le  C_{L_3,\epsilon,T}\exp\big(\frac{\|(-A)^{\frac{1-\epsilon}{4}}x\|}{2e^{L_3T}}\big)(1+t^{-\gamma_0}),\label{eq29}\\
\|Du^N(t,x)\|&\le C_{\mathfrak{p},T}(1+\|x\|^{\mathfrak{p}}),\label{ineq:normal-first}
\end{align}	
where $\alpha_0$ is given by Assumption~\ref{assum:nonlinearity-F}.
   \end{propo}

\begin{propo}\label{cor-secondderivative}
Let Assumptions \ref{assum:nonlinearity-F}--\ref{assum:nonlinearity-G} hold and $\varphi\in\IC^3_{\mathfrak{p}}(H;\R)$. For all $t\in(0,T]$, $\epsilon\in(0,\frac{1}{(\lambda_1\wedge 1)^{-3/4}-1}\wedge\frac{1}{4(\lambda_1\wedge 1)^{-1/2}-1}\wedge(\frac{1}{\alpha_0}-1))$, if there exists some $p>2$ such that $\gamma_1\in[\frac{1-\epsilon}{4},\frac{1-\epsilon}{4}+\frac{1-\epsilon}{2(1+2\epsilon)(1+\epsilon)}\wedge\frac{1}{p(1+2\epsilon)})$, then for all $\gamma_2\in[\frac{1-\epsilon}{4},\frac{1-\epsilon}{4}+\frac{1-\epsilon}{2(1+2\epsilon)(1+\epsilon)}\wedge\frac{p-2}{2p(1+2\epsilon)})$, there exist positive constants $L_3, C_{L_3,\epsilon,T}, C_{\mathfrak{p},T}$ such that for all $N\in\N_+$ and $x\in  H_N$
\begin{align}
\big\|(-A)^{\gamma_1}D^2 u^{N}(t,x)(-A)^{\gamma_2}\big\|_{\IL(H_N)}\!\!&\le  C_{L_3,\epsilon,T}\exp\big(\frac{2\|(-A)^{\frac{1-\epsilon}{4}}x\|}{e^{L_3T}}\big)(1+t^{-(\gamma_1+\gamma_2)}),\label{eq30}\\
\|D^2 u^{N}(t,x)\|_{\IL(H_N)}&\le C_{\mathfrak{p},T}\big(1+\|x\|^{\mathfrak{p}}\big),\label{ineq:normal-second}
\end{align}
where $\alpha_0$ is given by Assumption~\ref{assum:nonlinearity-F}.
\end{propo}

 We remark that, for the purely Gaussian multiplicative noise case, \cite{AD-2011} established regularity estimates analogous to \eqref{eq29} and \eqref{eq30} for the solution $u$ of the Kolmogorov equation associated to the original SPDE with $\gamma_0\in[0,\frac12)$ and $\gamma_1,\gamma_2\in[0,\frac12)$. Subsequently, \cite{CE-2018} improved result \eqref{eq29} by extending the admissible range of $\gamma_0$ to $[0,1)$ based on a discrete-time analysis framework and Duhamel's principle. To the best of our knowledge, \eqref{eq29} and \eqref{eq30} are the first results in the literature providing regularity estimates for the first and second order derivatives of the solutions to \eqref{eq:kolmogorovequation} in the presence of both Gaussian and Poisson multiplicative noise. Inspired by \cite{CE-2018}, to prove Propositions~\ref{cor:firstderivative2} and \ref{cor-secondderivative}, we introduce a discrete-time approximation of $u^N$, denoted by $u^{N,\tau}$, and 
obtain the regularity estimates for derivatives of $u^{N,\tau}$ (see Propositions~\ref{the:expression-first} and \ref{the:expression-second}) in the following subsection. By subsequently letting $\tau\to0$, we can obtain the corresponding estimates for $u^N$ (i.e., Propositions ~\ref{cor:firstderivative2} and \ref{cor-secondderivative}). The detailed proofs are provided in Subsection~\ref{sec:4.3}.

\subsection{Discrete-time approximation}\label{subsec:3.3.2}

For the discrete-time approximation, let $\mathcal T_\tau:=\{t_m: m=0,1,\ldots,M\}$. We define the mapping $u^{N,\tau}\colon \mathcal T_\tau\times H_N\to \R$ by $u^{N,\tau}(t_m,x):=\E[\varphi(\bar{X}^{N,(0,x)}_m)]$ for $m\in\{0,1,\ldots,M\}$ and $x\in H_N$, where $\bar{X}^{N,(0,x)}_m$ denotes the solution of \eqref{eq:one-step} at time $t_m$ starting from $\bar{X}^{N,(0,x)}_0=x$. 
For simplicity, we also write  $\bar X_m^N$ when the dependence on the initial data $x$ is not emphasized. Similar to \eqref{expression:uN-1}--\eqref{expression:uN-2}, it holds for all $x,y,h_1,h_2\in H_N$ and $m\in\{0,1,\ldots,M\}$ that 
\begin{align}
&Du^{N,\tau}(t_m,x).y=\E\big[D\varphi(\bar{X}^{N}_m).\eta^{N}_{y,x,m}\big],\label{expression:u-1}\\
&D^2 u^{N,\tau}(t_m,x).(h_1,h_2)=\E\big[D\varphi(\bar{X}^{N}_m). \xi^{N}_{h_1,h_2,x, m}\big]+\E\big[D^2\varphi(\bar{X}^{N}_m).(\eta^{N}_{h_1,x,m},\eta^{N}_{h_2,x,m})\big]\label{expression:u-2},
\end{align}
where $\eta^{N}_{y,x,m}:=D\bar{X}^{N,(0,x) }_m.y$ solves the problem
\begin{align}\label{expression:eta1}
\eta^{N }_{y,x,m}&=\!E(\tau)\Big(\eta^{N }_{y,x,m-1}\!+\!\tau P_NDF(\bar{X}^{N }_{m-1}).\eta^{N }_{y,x,m-1}\!+\!P_NDB(\bar{X}^{N }_{m-1}).\eta^{N }_{y,x,m-1}\Delta W^Q_{m-1}\nonumber\\
&~~~\quad+ \int_{t_{m-1}}^{t_{m}}\!\!\int_\chi  P_NDG(\bar{X}^{N }_{m-1},z).\eta^{N }_{y,x,m-1}\,\tilde{\mathbf{N}}(dz,ds)\Big)
\end{align}
starting from $\eta^{N }_{y,x,0}=y$ and $\xi^{N}_{h_1,h_2,x,m}:=D^2\bar{X}^{N,(0,x)}_m.(h_1,h_2)$ satisfies 
\begin{align}
\!\!\!\!&\xi^{N}_{h_1,h_2,x,m}=\tau \sum_{j=0}^{m-1}E((m-j)\tau)P_NDF(\bar{X}^{N}_{j}).\xi^{N}_{h_1,h_2,x,j}+\tau\sum_{j=0}^{m-1}E((m-j)\tau)P_ND^2F(\bar{X}^{N}_{j}).(\eta^{N}_{h_1,x,j},\eta^{N}_{h_2,x,j})\nonumber\\
&\quad\quad+\sum_{j=0}^{m-1}\!E((m-j)\tau)P_NDB(\bar{X}^{N}_{j}).\xi^{N}_{h_1,h_2,x,j}\Delta W^Q_{j}+\sum_{j=0}^{m-1}\!E((m-j)\tau)P_ND^2B(\bar{X}^{N}_{j}).(\eta^{N}_{h_1,x,j},\eta^{N}_{h_2,x,j})\Delta W^Q_{j}\nonumber\\
&\quad\quad +\sum_{j=0}^{m-1}E((m-j)\tau)\int_{t_j}^{t_{j+1}}\int_\chi P_N\big(DG(\bar{X}^{N}_{j},z).\xi^{N}_{h_1,h_2,x,j}+D^2G(\bar{X}^{N}_{j},z).(\eta^{N}_{h_1,x,j},\eta^{N}_{h_2,x,j})\big) \,\tilde{\mathbf{N}}(dz,dt)\label{expression:xi}
\end{align}
with initial value $\xi^{N}_{h_1,h_2,x,0}=0$. Then we have the following a priori estimates for $\eta^{N}_{y,x,m}$ and $\xi^{N}_{h_1,h_2,x,m}$. 

\begin{lemma}\label{lem:normal-1-2}
Let Assumptions \ref{assum:nonlinearity-F}--\ref{assum:nonlinearity-G} be satisfied. For any $\kappa_0\ge 1$ and $\tilde{\delta}\in(0,\frac12)$, there exist positive constants $C_{\kappa_0,T}$ and $C_{\tilde{\delta},T}$ such that, for all $N\in\N_+$, $\tau\in(0,1)$, $m\in\{0,1,\ldots,M\}$, and $x, y, h_1,h_2\in H_N$, 
\begin{align}
\|\eta^{N}_{y,x,m}\|_{L^{2\kappa_0}(\Omega;H)}&\le C_{\kappa_0,T}\|y\|,\label{estimate:eta-kappa}\\
\|\xi^{N}_{h_1,h_2,x,m}\|_{L^{\frac{4}{1+\tilde{\delta}}}(\Omega;H)}&\le C_{\tilde{\delta},T}\|h_1\|\|h_2\|.\label{estimate:xi-4-}
\end{align}
\end{lemma}

The proof of Lemma~\ref{lem:normal-1-2} is similar to that for \eqref{eq:estimate-X}, and thus is omitted here.
We emphasize that results in Lemma \ref{lem:normal-1-2} are not sharp enough to derive the regularity estimates for derivatives of the discrete-time approximation $u^{N,\tau}$, for instance, in the form
\begin{align}\label{estimate:u-first-full}
\|(-A)^{\gamma_0}Du^{N,\tau}(t_m,x)\|\le C(t_{m})K(x)\Big(1+\tau^{\bar{\alpha}}\times\sup_{y\in H_N\setminus\{0\}}\frac{\|y\|}{\|(-A)^{-\gamma_0}y\|}\Big),\quad m\in\{1,2,\ldots,M\}
\end{align}
with some $\gamma_0>0$, $\bar{\alpha}>0$, and positive constants $C(t_{m})$, $K(x)$ independent of $\tau$,
which is indispensable for proving \eqref{eq29}. 
Furthermore, in contrast to \cite{CE-2018}, where purely Gaussian multiplicative noise was treated, the passage to the setting with both the Gaussian and Poisson multiplicative noise in the present paper is not a mere extension: each power of the jump coefficient introduces an additional singular integral in the Malliavin-weighted Duhamel expansion. Consequently,  the estimates on the derivatives of $u^{N,\tau}$ cannot be inherited from \cite{CE-2018}, and new arguments  are required:   Rather than relying solely on the energy estimates in  state space $H$, we establish refined moment estimates for $\eta^{N}_{y,x,m}$ and $\xi^{N}_{h_1,h_2,x,m}$ in appropriate negative Sobolev norms, as shown in \eqref{eq31} and \eqref{eq32} of Lemmas~\ref{lem:estimate-eta} and \ref{lem:estimate-xi}.

\begin{lemma}\label{lem:estimate-eta}
Let Assumptions \ref{assum:nonlinearity-F}--\ref{assum:nonlinearity-G} be satisfied. Then for all $\tilde{p}\ge 2$, $\alpha_2\in[0,\frac{1}{\tilde{p}(1+2\epsilon_1)})$ with $\epsilon_1\in(0,\epsilon]$ and $\epsilon\in(0,\frac12\wedge\frac{1}{(\lambda_1\wedge 1)^{-3/4}-1}\wedge (\frac{1}{\alpha_0}-1))$, there exist constants $L_3>0$ and $C_{L_3,\tilde{p},\epsilon,T}>0$ such that for all $N\in\N_+$, $\tau\in(0, \tau_N)$ with $\tau_N$ given in \eqref{def:tau_0}, $m\in\{1,2,\ldots, M\}$, and $x,y\in H_N$ 
\begin{align}
&~~~\E\big[\big\|(-A)^{-\frac{1-\epsilon}{4}}\eta^{N}_{y,x,m}\big\|^{\tilde{p}}\big]\nonumber\\
&\le C_{L_3,\tilde{p},\epsilon,T}\big(\exp\big(\frac{\|(-A)^{\frac{1-\epsilon}{4}}x\|}{2e^{L_3T}}\big)+t_m^{-\alpha_2\tilde{p}}\big)\big(\tau^{\frac{1}{1+2\epsilon_1}}\big\|(-A)^{-\frac{1-\epsilon}{4}}y\big\|^{\tilde{p}}+ \big\|(-A)^{-\frac{1-\epsilon}{4}-\alpha_2}y\big\|^{\tilde{p}}\big),\label{eq31}
\end{align}
where $\alpha_0$ is given by Assumption~\ref{assum:nonlinearity-F}. In particular, for all $\iota\in[0,\frac{1}{2(1+\epsilon_1)(1+2\epsilon_1)})$ and $\iota_1\in[0,\frac{1-\epsilon_1}{2(1+\epsilon_1)(1+2\epsilon_1)})$, there exists a constant $C_{L_3,\epsilon,T}>0$ such that for all $N\in\N_+$, $\tau\in(0, \tau_N)$ with $\tau_N$ given in \eqref{def:tau_0}, $m\in\{1,2,\ldots, M\}$, and $x,y\in H_N$
\begin{align}
&\!\!\big\|(-A)^{\frac{1-\epsilon}{4}}\eta^{N}_{y,x,m}\big\|_{L^{2}(\Omega;H)}\le C_{L_3,\epsilon,T}\big(\exp\big(\frac{\|(-A)^{\frac{1-\epsilon}{4}}x\|}{4e^{L_3T}}\big)+t_m^{-\frac{1-\epsilon}{2}-\iota}+t_{m-1}^{-\iota+\frac{\epsilon}{2}}\times \mathbf 1_{\{m>1\}}\big)\nonumber\\
&\quad\quad\quad\quad\quad\quad\quad\quad\quad\quad\quad\quad\times \big(\tau^{\frac{1}{2(1+\epsilon_1)(1+2\epsilon_1)}}\big\|(-A)^{-\frac{1-\epsilon}{4}}y\big\|+ \big\|(-A)^{-\frac{1-\epsilon}{4}-\iota}y\big\|\big);\nonumber\\
&\big\|\eta^{N}_{y,x,m}\big\|_{L^{2(1+\epsilon_1)}(\Omega;H)}\le C_{L_3,\epsilon,T}\exp\big(\frac{\|(-A)^{\frac{1-\epsilon}{4}}x\|}{4(1+\epsilon_1)e^{L_3T}}\big)\big(1+t_m^{-\frac{1-\epsilon}{4}-\iota_1}+t_{m-1}^{-\frac{\iota_1}{2}+\frac{\epsilon}{4}}(1+t_m^{-\frac{\iota_1}{2}})\times \mathbf 1_{\{m>1\}}\big)\nonumber\\
&\quad\quad\quad\quad\quad\quad\quad\quad\quad\quad\quad\times\big(\tau^{\frac{1-\epsilon_1}{2(1+\epsilon_1)(1+2\epsilon_1)}}\big\|(-A)^{-\frac{1-\epsilon}{4}}y\big\|\!+\! \big\|(-A)^{-\frac{1-\epsilon}{4}-\iota_1}y\big\|\big).\label{estimate-eta3}
\end{align}	
\end{lemma}

\begin{lemma}\label{lem:estimate-xi}
Let Assumptions \ref{assum:nonlinearity-F}--\ref{assum:nonlinearity-G} be satisfied. Then for all $\epsilon\in(0, \frac{1}{(\lambda_1\wedge 1)^{-3/4}-1}\wedge\frac{1}{4(\lambda_1\wedge 1)^{-1/2}-1}\wedge (\frac{1}{\alpha_0}-1))$ and $\varepsilon_0\in(0,\epsilon)$, if there exists some $p>2$ such that $\varpi_1\in[0,\frac{1-\epsilon}{2(1+2\epsilon)(1+\epsilon)}\wedge\frac{1}{p(1+2\epsilon)})$, then for all $\varpi_2\in[0,\frac{1-\epsilon}{2(1+2\epsilon)(1+\epsilon)}\wedge\frac{p-2}{2p(1+2\epsilon)})$ such that $\varpi_1+\varpi_2\in[0,\frac{1}{2(1+2\epsilon)})$, 
there exist constants $L_3>0$ and  $C_{L_3,\epsilon,T}>0$ such that for all $N\in\N_+$, $\tau\in(0, \tau_N)$ with $\tau_N$ given in \eqref{def:tau_0}, and $x\in H_N$, $h_1,h_2\in H_N\setminus\{0\}$
\begin{align}\label{eq32}
\sup_{m\in\{0,1,\ldots, M\}}\big\|(-A)^{-\frac{1-\epsilon}{4}}&\xi^{N}_{h_1,h_2,x,m}\big\|_{L^{1+\varepsilon_0}(\Omega;H)}
\le C_{L_3,\epsilon,T}\exp\big(\frac{\|(-A)^{\frac{1-\epsilon}{4}}x\|}{e^{L_3T}}\big)\\
&\times\Big(\big\|(-A)^{-\frac{1-\epsilon}{4}-\varpi_1}h_1\big\|\big\|(-A)^{-\frac{1-\epsilon}{4}-\varpi_2}h_2\big\|+C_{h_1,h_2}\tau^{\frac{1-\epsilon}{2(1+\epsilon)(1+2\epsilon)}\wedge \frac{2\wedge (p-2)}{p(1+\epsilon)(1+2\epsilon)}\wedge\frac{1-\epsilon}{4}}\Big),\nonumber
\end{align}
where  $C_{h_1,h_2}:=\frac{\|h_1\|^{\frac{2}{1+\epsilon}}\|h_2\|^{\frac{2}{1+\epsilon}}}{\|(-A)^{-\frac{1-\epsilon}{4}-\varpi_1}h_1\|^{\frac{1-\epsilon}{1+\epsilon}}\|(-A)^{-\frac{1-\epsilon}{4}-\varpi_2}h_2\|^{\frac{1-\epsilon}{1+\epsilon}}}+\|h_1\|\|h_2\|$  and $\alpha_0$ is given by Assumption~\ref{assum:nonlinearity-F}. 
\end{lemma}


The proofs of Lemmas~\ref{lem:estimate-eta} and \ref{lem:estimate-xi} are based on a stochastic variation-of-constants decomposition. By exploiting the smoothing property of the semigroup, this strategy effectively compensates for the singular effect induced by the jump coefficient and yields refined estimates for the variational processes in negative Sobolev norms. For the reader's convenience, we present the detailed proof of Lemma~\ref{lem:estimate-eta} in the Appendix. The proof of Lemma \ref{lem:estimate-xi} follows by the same arguments and is therefore omitted.


Combined with the Malliavin--Duhamel mechanism, these bounds of \eqref{eq31} and \eqref{eq32} allow for enlarged admissible ranges of $\gamma_0$ in \eqref{estimate:u-first-full}, as well as wider admissible ranges of $\gamma_1,\gamma_2$ in the estimates for the second order derivatives of the solution $u^{N,\tau}$. Propositions \ref{the:expression-first} and \ref{the:expression-second} below identify the admissible ranges of the indices $\gamma_0$ and $\gamma_1,\gamma_2$.

\begin{propo}\label{the:expression-first}
Let Assumptions \ref{assum:nonlinearity-F}--\ref{assum:nonlinearity-G} hold and $\varphi\in\IC^2_{\mathfrak{p}}(H;\R)$. Then for all $\epsilon\in(0,\frac12\wedge\frac{1}{(\lambda_1\wedge 1)^{-3/4}-1}\wedge (\frac{1}{\alpha_0}-1))$, $\epsilon_1\in(0,\epsilon]$, and $\gamma_0\in[0,\frac{1-\epsilon}{4}+\frac{1-\epsilon_1}{2(1+\epsilon_1)(1+2\epsilon_1)})$, there exist  constants $L_3>0$ and $C_{L_3,\epsilon,T}>0$ such that for all $N\in\N_+$, $\tau\in(0, \tau_N)$ with $\tau_N$ given in  \eqref{def:tau_0}, $m\in\{1,2,\ldots,M\}$, and $x,y\in H_N$  
\begin{align*}
|Du^{N,\tau}(t_m,x).y|&\le C_{L_3,\epsilon,T}\exp\big(\frac{\|(-A)^{\frac{1-\epsilon}{4}}x\|}{2e^{L_3T}}\big)\times\tau^{\frac{1-\epsilon_1}{2(1+\epsilon_1)(1+2\epsilon_1)}}\|y\|\nonumber\\
&~~~+C_{L_3,\epsilon,T}\exp\big(\frac{\|(-A)^{\frac{1-\epsilon}{4}}x\|}{2e^{L_3T}}\big) \times (1+t_{m}^{-\gamma_0})\|(-A)^{-\gamma_0}y\|,
\end{align*}
where $\alpha_0$ is given by Assumption~\ref{assum:nonlinearity-F}. In particular, there exists a constant $C_{\mathfrak{p},T}>0$ such that for all $N\in\N_+$, $\tau\in(0,1)$, $m\in\{0,1,\ldots,M\}$, and $x,y\in H_N$
\begin{align}\label{ineq:normal-first-discrete}
|Du^{N,\tau}(t_m,x).y|\le C_{\mathfrak{p},T}\big(1+\|x\|^{\mathfrak{p}}\big)\|y\|.
\end{align}
\end{propo}

\begin{propo}\label{the:expression-second}
Let Assumptions \ref{assum:nonlinearity-F}--\ref{assum:nonlinearity-G} hold and $\varphi\in\IC^2_{\mathfrak{p}}(H;\R)$. For all $\epsilon\in(0,\frac{1}{(\lambda_1\wedge 1)^{-3/4}-1}\wedge\frac{1}{4(\lambda_1\wedge 1)^{-1/2}-1}\wedge (\frac{1}{\alpha_0}-1))$, if there exists some $p>2$ such that $\gamma_1\in[\frac{1-\epsilon}{4},\frac{1-\epsilon}{4}+\frac{1-\epsilon}{2(1+2\epsilon)(1+\epsilon)}\wedge\frac{1}{p(1+2\epsilon)})$, then for all $\gamma_2\in[\frac{1-\epsilon}{4},\frac{1-\epsilon}{4}+\frac{1-\epsilon}{2(1+2\epsilon)(1+\epsilon)}\wedge\frac{p-2}{2p(1+2\epsilon)})$, there exist  constants $L_3>0$ and $C_{L_3,\epsilon,T}>0$  such that for all $N\in\N_+$, $\tau\in(0, \tau_N)$ with $\tau_N$ given in \eqref{def:tau_0}, $m\in\{1,2,\ldots,M\}$, and $x\in H_N$, $h_1,h_2\in H_N\setminus\{0\}$ 
\begin{align*}
|D^2 u^{N,\tau}(t_m,x).(h_1,h_2)|&\le C_{L_3,\epsilon,T}\exp\big(\frac{2\|(-A)^{\frac{1-\epsilon}{4}}x\|}{e^{L_3T}}\big)\big(1+t_{m}^{-(\gamma_1+\gamma_2)}+t_{m-1}^{-(\gamma_1+\gamma_2)}\mathbf{1}_{\{m>1\}}+t_{m-2}^{-(\gamma_1+\gamma_2)}\mathbf 1_{\{m>2\}}\big)\\
&\quad\quad \times \big(\|(-A)^{-\gamma_1}h_1\|\|(-A)^{-\gamma_2}h_2\|+C_{h_1,h_2}\tau^{\frac{1-\epsilon}{2(1+\epsilon)(1+2\epsilon)}\wedge \frac{2\wedge (p-2)}{p(1+\epsilon)(1+2\epsilon)}\wedge\frac{1-\epsilon}{4}}\big),
\end{align*}
where $C_{h_1,h_2}$ is given in Lemma \ref{lem:estimate-xi} and $\alpha_0$ is given by Assumption~\ref{assum:nonlinearity-F}. In particular, there exists a constant $C_{\mathfrak{p},T}>0$ such that for all $N\in\N_+$, $\tau\in(0,1)$, $m\in\{0,1,\ldots,M\}$, and $x,h_1,h_2\in H_N$
    \begin{align}\label{ineq:normal-second-discrete}
        |D^2u^{N,\tau}(t_m,x).(h_1,h_2)|\le C_{\mathfrak{p},T}\big(1+\|x\|^{\mathfrak{p}}\big)\|h_1\|\|h_2\|.
    \end{align}
\end{propo}

\subsection{Proofs of Propositions in Section~\ref{sec:4}
}\label{sec:4.3}

This subsection provides the proofs of the main results of Section~\ref{sec:4}, namely,  Propositions \ref{cor:firstderivative2} and \ref{cor-secondderivative} for the spatial approximation $u^{N}$. The proofs rely on intermediate results established for the discrete-time approximation $u^{N,\tau}$ in Propositions \ref{the:expression-first} and \ref{the:expression-second}. Hence, we first present the proofs of Propositions \ref{the:expression-first} and \ref{the:expression-second} in Subsections \ref{subsection4.3.1} and \ref{subsection4.3.2}, respectively, 
and then prove Propositions \ref{cor:firstderivative2} and \ref{cor-secondderivative} in Subsections \ref{subsection4.3.3} and \ref{subsection4.3.4} by letting $\tau\to 0$.

\subsubsection{Proof of Proposition \ref{the:expression-first} }\label{subsection4.3.1}


By expression \eqref{expression:eta1} of $\{\eta^{N}_{y,x,m}\}_{m\in\{0,1,\cdots,M\}}$, we have 
\begin{align}\label{decompose-T}
Du^{N,\tau}(t_m,x).y&=\E\big[D\varphi(\bar{X}^{N}_m).\eta^{N}_{y,x,m}\big]=:\mathbb{T}_1+\mathbb{T}_2+\mathbb{T}_3,
\end{align}
where  $\mathbb{T}_1=\E\big[D\varphi(\bar{X}^{N}_m).E(m\tau)y\big]+\sum_{\ell=0}^{m-1}\mathbb{T}_{1,\ell}$, $\mathbb{T}_2=\sum_{\ell=0}^{m-1}\mathbb{T}_{2,\ell}$, and $\mathbb{T}_3=\sum_{\ell=0}^{m-1}\mathbb{T}_{3,\ell}$ with 
\begin{align*}
\mathbb{T}_{1,\ell}&:=\E\Big[\int_{t_\ell}^{t_{\ell+1}}\!\!\big\langle D\varphi(\bar{X}^{N}_m),E((m-\ell)\tau)P_NDF(\bar{X}^{N}_{\ell}).\eta^{N}_{y,x,\ell}\big\rangle \,ds\Big],\nonumber\\
	\mathbb{T}_{2,\ell}&:=\E\Big[\int_{t_\ell}^{t_{\ell+1}}\!\!\Big\langle D\varphi(\bar{X}^{N}_m),E((m-\ell)\tau)P_NDB(\bar{X}^{N}_{\ell}).\eta^{N}_{y,x,\ell}\,dW^Q(s)\Big\rangle \Big],\nonumber\\
	\mathbb{T}_{3,\ell}&:=\E\Big[\int_{t_\ell}^{t_{\ell+1}}\Big\langle D\varphi(\bar{X}^{N}_m), E((m-\ell)\tau) \int_\chi\!P_NDG(\bar{X}^{N}_{\ell},z).\eta^{N}_{y,x,\ell}\,\tilde{\mathbf{N}}(dz,ds)\Big\rangle \Big].
\end{align*}

\textbf{Estimate of term $\mathbb{T}_1$}. Under condition on $DF$ in Assumption \ref{assum:nonlinearity-F} (ii) with $\varsigma_0=\varsigma=\epsilon$, we derive from \eqref{prop:semigroup1}, the H\"older inequality, Corollary~\ref{cor:exponentialintegrability}, and Lemma~\ref{lem:estimate-eta} with $\alpha_2\in[0,\frac{1}{2(1+2\epsilon_1)})$ and $\epsilon_1\in(0,\epsilon]$ that for all $\ell\in\{1,2,\ldots,m-1\}$
\begin{align}\label{ineq:T1}
	|\mathbb{T}_{1,\ell}|&\le C_{\mathfrak{p},\epsilon}\big(1\!+\!\big(\E[\|\bar{X}^{N}_m\|^{\mathfrak{p}(1+1/\epsilon)}]\big)^{\frac{\epsilon}{1+\epsilon}}\big)  \int_{t_\ell}^{t_{\ell+1}}t_{m-\ell}^{-\frac{1+\epsilon}{4}}\times \Big|\big(1\!+\!\|(-A)^{\frac{1-\epsilon}{4}}\bar{X}^{N}_{\ell}\|\big)\big\|(-A)^{-\frac{1-\epsilon}{4}}\eta^{N}_{y,x,\ell}\big\|\Big|_{L^{1+\epsilon}(\Omega;\R)}ds\nonumber\\
	&\le C_{L_3,\epsilon,T}\exp\big(\frac{\|(-A)^{\frac{1-\epsilon}{4}}x\|}{4e^{L_3T}}\big)\int_{t_\ell}^{t_{\ell+1}}t_{m-\ell}^{-\frac{1+\epsilon}{4}}\times \big\|(-A)^{-\frac{1-\epsilon}{4}}\eta^{N}_{y,x,\ell}\big\|_{L^{2}(\Omega;H)}\,ds\\
	&\le C_{L_3,\epsilon,T}\exp\big(\frac{\|(-A)^{\frac{1-\epsilon}{4}}x\|}{2e^{L_3T}}\big)\int_{t_\ell}^{t_{\ell+1}}t_{m-\ell}^{-\frac{1+\epsilon}{4}}(1+t_{\ell}^{-\alpha_2})\,ds\times \big(\tau^{\frac{1}{2(1+2\epsilon_1)}}\big\|(-A)^{-\frac{1-\epsilon}{4}}y\big\|+ \big\|(-A)^{-\frac{1-\epsilon}{4}-\alpha_2}y\big\|\big).\nonumber
\end{align}
In particular,  it holds that $|\mathbb{T}_{1,0}|\le C_{\mathfrak{p}}(1+\E[\|\bar{X}^N_m\|^{\mathfrak{p}}])\int_0^{t_1}\|DF(P_Nx).y\|\,ds\le C_{\mathfrak{p},T}\tau(1+\|x\|^{\mathfrak{p}})\|y\|$ by \eqref{eq:estimate-X}. Then combining the smoothing property of  the semigroup and estimate \eqref{ineq:T1} gives
\begin{align*}
	|\mathbb{T}_1|&\le \E\big[|D\varphi(\bar{X}^{N}_m).E(m\tau)y|\big]+\sum_{\ell=0}^{m-1} |\mathbb{T}_{1,\ell}|\le C_{L_3,\epsilon,T}\exp\big(\frac{\|(-A)^{\frac{1-\epsilon}{4}}x\|}{2e^{L_3T}}\big)\times\Big\{t_m^{-\gamma_0}\|(-A)^{-\gamma_0}y\|+\tau \|y\|\\
	&~~~+\int_{t_1}^{t_m}(t_{m}-\lfloor s\rfloor)^{-\frac{1+\epsilon}{4}}(1+\lfloor s\rfloor^{-\alpha_2})\,ds\times \big(\tau^{\frac{1}{2(1+2\epsilon_1)}}\big\|(-A)^{-\frac{1-\epsilon}{4}}y\big\|+ \big\|(-A)^{-\frac{1-\epsilon}{4}-\alpha_2}y\big\|\big)\Big\},
\end{align*}
where $\gamma_0\le \frac{1-\epsilon}{4}+\alpha_2<\frac{1-\epsilon}{4}+\frac{1}{2(1+2\epsilon_1)}$.  
By the scaling property of the beta function (cf. \cite[Chapter 6]{AS-1964}), we obtain that  for $m>1$,
\begin{align}\label{ineq:beta}
\int_{t_1}^{t_m}(t_m-\lfloor s\rfloor)^{-\frac{1+\epsilon}{4}}(1+\lfloor s\rfloor^{-\alpha_2})\,ds&\le C\int_{t_1}^{t_m}(1+(t_m-s)^{-\frac{1+\epsilon}{4}})(1+(s-\tau)^{-\alpha_2})\,ds\nonumber\\
&\le (C_T+\mathcal{B}(\frac{3-\epsilon}{4},1-\alpha_2))(1+t_{m-1}^{\frac{3-\epsilon}{4}-\alpha_2}).
\end{align}
Hence it follows that
$
	|\mathbb{T}_1|
	\le C_{L_3,\epsilon,T}\exp\big(\frac{\|(-A)^{\frac{1-\epsilon}{4}}x\|}{2e^{L_3T}}\big)\times\big(\tau^{\frac{1}{2(1+2\epsilon_1)}}\big\|y\big\|+(1+t_m^{-\gamma_0})\times\|(-A)^{-\gamma_0}y\|\big).
$

\textbf{Estimate of term $\mathbb{T}_2$.} Applying the Wiener duality formula in Proposition~\ref{prop:dualityformula} together
with the chain rule \eqref{eq:chainrule}, we obtain 
\begin{align*}
\mathbb{T}_{2,\ell}=\sum_{i\in\N_+}\E\Big[\int_{t_\ell}^{t_{\ell+1}}\!\!\!D^2\varphi(\bar{X}^{N}_m).\Big(\ID^{W,Q^{\frac12}e_i}_s\bar{X}^{N}_m,E((m-\ell)\tau)P_NDB(\bar{X}^{N}_{\ell}).\eta^{N}_{y,x,\ell}Q^{\frac12}e_i\Big)\,ds\Big].
\end{align*}
For all $\ell\in\{1,2,\ldots,m-1\}$, under Assumption \ref{assum:nonlinearity-B}, using the H\"older inequality,   Lemma~\ref{lem:Malliavinderivative-B1}, and Corollary~\ref{cor:exponentialintegrability} yields  
\begin{align}\label{ineq:T2}
	&|\mathbb{T}_{2,\ell}|=\Bigl|\sum_{i\in\N_+}\E\Big[\int_{t_\ell}^{t_{\ell+1}}\!\!\!D^2\varphi(\bar{X}^{N}_m).\Big(\ID^{W,Q^{\frac12}e_i}_s\bar{X}^{N}_m,E((m-\ell)\tau)P_NDB(\bar{X}^{N}_{\ell}).\eta^{N}_{y,x,\ell}Q^{\frac12}e_i\Big)\,ds\Big]\Bigr|\\
	&\le C_{\mathfrak{p}}\int_{t_\ell}^{t_{\ell+1}}\E\Big[(1+\|\bar{X}^{N}_m\|^{\mathfrak{p}})\times\sum_{i\in\N_+}\big\|\ID^{W,Q^{\frac12}e_i}_s\bar{X}^{N}_m\big\|\big\|E((m-\ell)\tau)P_NDB(\bar{X}^{N}_{\ell}).\eta^{N}_{y,x,\ell}Q^{\frac12}e_i\big\|\Big]\,ds\nonumber\\
    &\le C_{\mathfrak{p}}\big(1\!+\!\big(\E[\|\bar{X}^{N}_m\|^{\mathfrak{p}(1+1/\epsilon_1)}]\big)^{\frac{\epsilon_1}{1+\epsilon_1}}\big)\!\int_{t_\ell}^{t_{\ell+1}}\!\!\Big|\sum_{i\in\N_+}\big\|\ID^{W,Q^{\frac12}e_i}_s\bar{X}^{N}_m\big\|\big\|E((m-\ell)\tau)P_NDB(\bar{X}^{N}_{\ell}).\eta^{N}_{y,x,\ell}Q^{\frac12}e_i\big\|\Big|_{L^{1+\epsilon_1}(\Omega;\R)}\!\!\!\!ds\nonumber\\
    &\le C_{\mathfrak{p},T}\big(1\!+\!\big(\E[\|\bar{X}^{N}_m\|^{\mathfrak{p}(1+1/\epsilon_1)}]\big)^{\frac{\epsilon_1}{1+\epsilon_1}}\big)\int_{t_\ell}^{t_{\ell+1}}\!\!t_{m-\ell}^{-\frac{1-\epsilon}{4}}\big\|(-A)^{-\frac{1-\epsilon}{4}}\eta^{N}_{y,x,\ell}\big\|_{L^{\frac{2(1+\epsilon_1)}{1-\epsilon_1}}(\Omega;H)}\big|\|B(\bar{X}^{N}_{\ell})\|_{\IL_2^0(H)}\big|_{L^2(\Omega;\R)}ds\nonumber\\
    &\le C_{L_3,\epsilon,T}\exp\big(\frac{\|(-A)^{\frac{1-\epsilon}{4}}x\|}{2e^{L_3T}}\big)\int_{t_\ell}^{t_{\ell+1}}\!\!t_{m-\ell}^{-\frac{1-\epsilon}{4}}(1+t_{\ell}^{-\gamma_0+\frac{1-\epsilon}{4}})\,ds \times \big(\tau^{\frac{1-\epsilon_1}{2(1+\epsilon_1)(1+2\epsilon_1)}}\big\|(-A)^{-\frac{1-\epsilon}{4}}y\big\|+ \big\|(-A)^{-\gamma_0}y\big\|\big),\nonumber
\end{align}
where $\gamma_0\in[0,\frac{1-\epsilon}{4}+\frac{1-\epsilon_1}{2(1+\epsilon_1)(1+2\epsilon_1)})$. When $\ell=0$, we derive similarly that 
\begin{small}\begin{align}\label{ineq:T-2,0}
	|\mathbb{T}_{2,0}|&\le C_{\mathfrak{p}}\big(1\!+\!\big(\E[\|\bar{X}^{N}_m\|^{\mathfrak{p}(1+1/\epsilon_1)}]\big)^{\frac{\epsilon_1}{1+\epsilon_1}}\big)\int_0^{t_1}\Big|\sum_{i\in\N_+}\big\|\ID^{W,Q^{\frac12}e_i}_s\bar{X}^{N}_m\big\|\big\|E(t_m)P_N(DB(P_Nx).y)Q^{\frac12}e_i\big\|\Big|_{L^{1+\epsilon_1}(\Omega;\R)}ds\nonumber\\
	&\le C_{\mathfrak{p}}\big(1\!+\!\big(\E[\|\bar{X}^{N}_m\|^{\mathfrak{p}(1+1/\epsilon_1)}]\big)^{\frac{\epsilon_1}{1+\epsilon_1}}\big)\int_0^{t_1} \big\|E(t_m)P_N(DB(P_Nx).y)\big\|_{\IL_2^0(H)} \Big(\sum_{i\in\N_+}\E\big[\big\|\ID^{W,Q^{\frac12}e_i}_s\bar{X}^{N}_m\big\|^2\big]\Big)^{\frac12}ds\nonumber\\
	&\le C_{L_3,\epsilon,T}\exp\big(\frac{\|(-A)^{\frac{1-\epsilon}{4}}x\|}{4e^{L_3T}}\big)\times\tau\|y\|\|B(P_Nx)\|_{\IL_2^0(H)}.
\end{align}\end{small}
Hence combining~\eqref{ineq:T2} and \eqref{ineq:T-2,0} yields
\begin{align*}
	|\mathbb{T}_2|\le\sum_{\ell=0}^{m-1}|\mathbb{T}_{2,\ell}|&\le  C_{L_3,\epsilon,T}\exp\big(\frac{\|(-A)^{\frac{1-\epsilon}{4}}x\|}{2e^{L_3T}}\big)\times\Big(\tau\|y\|+\int_{t_1}^{t_{m}}(t_m-\lfloor s\rfloor)^{-\frac{1-\epsilon}{4}}(1+\lfloor s\rfloor^{-\gamma_0+\frac{1-\epsilon}{4}})\,ds\nonumber\\
&\quad\quad\times \big(\tau^{\frac{1-\epsilon_1}{2(1+\epsilon_1)(1+2\epsilon_1)}}\big\|(-A)^{-\frac{1-\epsilon}{4}}y\big\|+ \big\|(-A)^{-\gamma_0}y\big\|\big)\Big)\nonumber\\
&\le   C_{L_3,\epsilon,T}\exp\big(\frac{\|(-A)^{\frac{1-\epsilon}{4}}x\|}{2e^{L_3T}}\big)\times\big(\tau^{\frac{1-\epsilon_1}{2(1+\epsilon_1)(1+2\epsilon_1)}}\|y\|+\|(-A)^{-\gamma_0}y\|\big),
\end{align*}
where in the last step we used a similar argument as in \eqref{ineq:beta}.

\textbf{Estimate of term $\mathbb{T}_3$.} Applying the Poisson duality formula in Proposition~\ref{prop:dualityformula}  together
with the chain rule \eqref{eq:chainrule}, we obtain
\begin{align*}
\mathbb{T}_{3,\ell}=\E\Big[\int_{t_\ell}^{t_{\ell+1}}\!\!\int_\chi \big\langle D\varphi\big(\bar{X}^{N}_m+\ID^{\mathbf{N}}_{(s,z)}\bar{X}^{N}_m\big)-D\varphi(\bar{X}^{N}_m),E((m-\ell)\tau)P_NDG(\bar{X}^{N}_{\ell},z).\eta^{N}_{y,x,\ell} \big\rangle\,\nu(dz)ds\Big].
\end{align*}
Using  the result \eqref{ineq:DN1} in Lemma~\ref{lem:Malliavinderivative-B1},  \eqref{prop:semigroup1}, \eqref{eq:estimate-X}, and Lemma~\ref{lem:estimate-eta} with $\alpha_2\in[0,\frac{1}{2(1+2\epsilon_1)})$ and $\epsilon_1\in(0,\epsilon]$ yields
\begin{align*}
		|\mathbb{T}_3|&\le C_{\mathfrak{p}}\!\sum_{\ell=0}^{m-1}\int_{t_\ell}^{t_{\ell+1}}\!\!\int_\chi  \E\Big[(1+\|\bar{X}^{N}_m\|^{\mathfrak{p}}+\|\ID^{\mathbf{N}}_{(s,z)}\bar{X}^{N}_m\|^{\mathfrak{p}})\times  \big\|E((m-\ell)\tau)P_NDG(\bar{X}^{N}_{\ell},z).\eta^{N}_{y,x,\ell} \big\|\Big] \nu(dz)ds\\
		&\le C_{\mathfrak{p}}\!\sum_{\ell=0}^{m-1}\int_{t_\ell}^{t_{\ell+1}}\!\! \E\Big[(1+\|\bar{X}^{N}_m\|^{\mathfrak{p}})  \int_\chi \big\|E((m-\ell)\tau)P_NDG(\bar{X}^{N}_{\ell},z).\eta^{N}_{y,x,\ell} \big\| \nu(dz)\Big]ds\\
		&~~~+C_{\mathfrak{p}}\!\sum_{\ell=0}^{m-1}\int_{t_\ell}^{t_{\ell+1}}\!\!\Big(\int_\chi \E\big[ \|\ID^{\mathbf{N}}_{(s,z)}\bar{X}^{N}_m\|^{2\mathfrak{p}}\big]\nu(dz)\Big)^{\frac12} \Big(\E\Big[\int_\chi  \big\|E((m-\ell)\tau)P_NDG(\bar{X}^{N}_{\ell},z).\eta^{N}_{y,x,\ell} \big\|^2\nu(dz)\Big]\Big)^{\frac12}ds\\
		&\le C_{\mathfrak{p}}\int^{t_{1}}_{0}(1+(\E[\|\bar{X}^{N}_m\|^{2\mathfrak{p}}])^{\frac12})\times \Bigl|\int_{\chi}\|E(m\tau) DG(P_Nx,z).y\|\,\nu(dz)\Bigr|_{L^2(\Omega;\R)}ds\\
		&~~~+C_{\mathfrak{p}}\int_{t_{1}}^{t_{m}}(t_{m}\!-\!\lfloor s\rfloor)^{-\frac{1-\epsilon}{4}}(1+(\E[\|\bar{X}^{N}_m\|^{2\mathfrak{p}}])^{\frac12})\times \Big|\int_{\chi}\|(-A)^{-\frac{1-\epsilon}{4}}DG(\bar{X}^{N}_{\lfloor s\rfloor/\tau},z).\eta^{N}_{y,x,\lfloor s\rfloor/\tau} \big\|\nu(dz)\Big|_{L^2(\Omega;\R)}ds\\
		&~~~+C_{\mathfrak{p}}\int^{t_{1}}_{0}\Big(\int_\chi \E\big[ \|\ID^{\mathbf{N}}_{(s,z)}\bar{X}^{N}_m\|^{2\mathfrak{p}}\big]\,\nu(dz)\Big)^{\frac12}\Big( \int_{\chi} \|E(m\tau) DG(P_Nx,z).y\|^2\,\nu(dz)\Big)^{\frac12}\,ds\\
		&~~~+C_{\mathfrak{p}}\int_{t_{1}}^{t_{m}}\!\!(t_{m}\!-\!\lfloor s\rfloor)^{-\frac{1-\epsilon}{4}}\Big(\int_\chi\!\! \E\big[ \|\ID^{\mathbf{N}}_{(s,z)}\bar{X}^{N}_m\|^{2\mathfrak{p}}\big]\,\nu(dz)\Big)^{\frac12}\Big(\E\Big[\int_{\chi}\!\!\|(-A)^{-\frac{1-\epsilon}{4}}DG(\bar{X}^{N}_{\lfloor s\rfloor/\tau},z).\eta^{N}_{y,x,\lfloor s\rfloor/\tau} \big\|^2\nu(dz)\Big]\Big)^{\frac12}ds\\
		&\le C_{\mathfrak{p},T}\Big(1+\|x\|^{\mathfrak{p}}+\big(\int_{\chi}|\bar{c}(z)|^{2\mathfrak{p}}\nu(dz)\big)^{\frac12}\Big)\times \tau\|y\|\\
		&\quad+C_{\mathfrak{p},T}\Big(1+\|x\|^{\mathfrak{p}}+\big(\int_{\chi}|\bar{c}(z)|^{2\mathfrak{p}}\nu(dz)\big)^{\frac12}\Big)\times \int_{t_{1}}^{t_{m}} (t_{m}\!-\!\lfloor s\rfloor)^{-\frac{1-\epsilon}{4}}\|(-A)^{-\frac{1-\epsilon}{4}}\eta^{N}_{y,x,\lfloor s\rfloor/\tau} \big\|_{L^2(\Omega;H)}\,ds\\
		&\le C_{L_3,\epsilon,T}\exp\big(\frac{\|(-A)^{\frac{1-\epsilon}{4}}x\|}{2e^{L_3T}}\big)\times\big(\tau^{\frac{1}{2(1+2\epsilon_1)}}\|y\|+\|(-A)^{-\gamma_0}y\|\big),
\end{align*}
where $\gamma_0\in[0,\frac{1-\epsilon}{4}+\frac{1}{2(1+2\epsilon_1)})$. 

Consequently, combining the estimates of $\mathbb T_1$, $\mathbb T_2$, and $\mathbb T_3$ with decomposition \eqref{decompose-T}, we obtain that for all $\gamma_0\in[0,\frac{1-\epsilon}{4}+\frac{1-\epsilon_1}{2(1+\epsilon_1)(1+2\epsilon_1)})$
\begin{align*}
|Du^{N,\tau}(t_m,x).y|&\le C_{L_3,\epsilon,T}\exp\big(\frac{\|(-A)^{\frac{1-\epsilon}{4}}x\|}{2e^{L_3T}}\big)\times\tau^{\frac{1-\epsilon_1}{2(1+\epsilon_1)(1+2\epsilon_1)}}\|y\|\nonumber\\
	&\quad+C_{L_3,\epsilon,T}\exp\big(\frac{\|(-A)^{\frac{1-\epsilon}{4}}x\|}{2e^{L_3T}}\big)(1+t_{m}^{-\gamma_0})  \|(-A)^{-\gamma_0}y\|.
\end{align*}	
In particular, we derive from \eqref{eq:estimate-X} and \eqref{estimate:eta-kappa} that
\begin{align*}
 |Du^{N,\tau}(t_m,x).y|&=\big|\E\big[D\varphi(\bar{X}^{N,(0,x)}_m).\eta^{N}_{y,x,m}\big]\big|\le C_{\mathfrak{p}}\E\big[\big(1+\|\bar{X}^{N,(0,x)}_m\|^{\mathfrak{p}}\big)\times \|\eta^{N}_{y,x,m}\|\big]\le C_{\mathfrak{p},T}\|y\|(1+\|x\|^{\mathfrak{p}}).
\end{align*}
The proof is completed. 

\subsubsection{Proof of Proposition \ref{the:expression-second}}\label{subsection4.3.2}

By expression \eqref{expression:xi} of $\{\xi^{N}_{h_1,h_2,x,m}\}_{m\in\{0,1,\cdots,M\}}$, we have
\begin{align}\label{decompose-S}
\E\big[D\varphi(\bar{X}^{N}_m).\xi^{N}_{h_1,h_2,x,m}\big]=: \mathbb{S}_1+\mathbb{S}_2+\mathbb{S}_3,
\end{align}
where for $m\ge 1$, $\mathbb{S}_i=\sum_{\ell=0}^{m-1}\mathbb{S}_{i,\ell}$, $i=1,2,3$, with the terms given by
\begin{align*}
	\mathbb{S}_{1,\ell}&\!:=\!\tau\times \E\Big[\Big\langle D\varphi(\bar{X}^{N}_m), E((m\!-\!\ell)\tau)P_N\big(DF(\bar{X}^{N}_{\ell}).\xi^{N}_{h_1,h_2,x,\ell}+D^2F(\bar{X}^{N}_{\ell}).(\eta^{N}_{h_1,x,\ell},\eta^{N}_{h_2,x,\ell})\big)\Big\rangle\Big],\\
	\mathbb{S}_{2,\ell}&\!:=\!\E\Big[\Big\langle D\varphi(\bar{X}^{N}_m),\int_{t_\ell}^{t_{\ell+1}}E((m-\ell)\tau)P_N\big(DB(\bar{X}^{N}_{\ell}).\xi^{N}_{h_1,h_2,x,\ell}+D^2B(\bar{X}^{N}_{\ell}).(\eta^{N}_{h_1,x,\ell},\eta^{N}_{h_2,x,\ell})\big)dW^Q(s)\Big\rangle\Big],\\
	\mathbb{S}_{3,\ell}&\!:=\!\E\Big[\Big\langle \!D\varphi(\bar{X}^{N}_m),\int_{t_\ell}^{t_{\ell+1}}\!\!\!\!\int_\chi \!\!E((m\!-\!\ell)\tau)P_N\big(DG(\bar{X}^{N}_{\ell},z).\xi^{N}_{h_1,h_2,x,\ell}\!+\!D^2G(\bar{X}^{N}_{\ell},z).(\eta^{N}_{h_1,x,\ell},\eta^{N}_{h_2,x,\ell})\big)\tilde{\mathbf{N}}(dz,ds)\Big\rangle\Big].
\end{align*}

\textbf{Estimate of term $\mathbb{S}_1$.} 
By the assumptions on $DF$ and $D^2F$ in Assumption \ref{assum:nonlinearity-F} (ii) with $\varsigma_0\in(0,1)$ and $\varsigma=\epsilon$, we derive for all $\ell\in\{1,2,\ldots,m-1\}$ that
\begin{align*}
	|\mathbb{S}_{1,\ell}|&\le C_{\mathfrak{p}}\tau\E\big[(1+\|\bar{X}^{N}_m\|)^{\mathfrak{p}}\times \big\|E((m\!-\!\ell)\tau)P_N\big(DF(\bar{X}^{N}_{\ell}).\xi^{N}_{h_1,h_2,x,\ell}\!+\!D^2F(\bar{X}^{N}_{\ell}).(\eta^{N}_{h_1,x,\ell},\eta^{N}_{h_2,x,\ell})\big)\big\|\big]\\
&\le C_{\mathfrak{p}}(1+(\E[\|\bar{X}^{N}_{m}\|^{\mathfrak{p}(1+1/\delta)}])^{\frac{\delta}{1+\delta}})\times\tau \big(\big\|E((m-\ell)\tau)P_NDF(\bar{X}^{N}_{\ell}).\xi^{N}_{h_1,h_2,x,\ell}\big\|_{L^{1+\delta}(\Omega;H)}\\
& \quad+\big\|E((m-\ell)\tau)P_ND^2F(\bar{X}^{N}_{\ell}).(\eta^{N}_{h_1,x,\ell},\eta^{N}_{h_2,x,\ell})\big\|_{L^{1+\delta}(\Omega;H)}\big)\\
&\le C_{\mathfrak{p}}(1+(\E[\|\bar{X}^{N}_{m}\|^{\mathfrak{p}(1+1/\delta)}])^{\frac{\delta}{1+\delta}})\times\int_{t_\ell}^{t_{\ell+1}}t_{m-\ell}^{-\frac{1+\varsigma_0}{4}}\times \big[\big(1+\|(-A)^{\frac{1-\epsilon}{4}}\bar{X}^{N}_{\ell}\|_{L^{\frac{(1+\delta)(1+2\delta)}{\delta}}(\Omega;H)}\big)\\
&~~~\times \|(-A)^{-\frac{1-\epsilon}{4}}\xi^{N}_{h_1,h_2,x,\ell}\|_{L^{1+2\delta}(\Omega;H)}+\|\eta^{N}_{h_1,x,\ell}\|_{L^{2(1+\delta)}(\Omega;H)}\|\eta^{N}_{h_2,x,\ell}\|_{L^{2(1+\delta)}(\Omega;H)}\big]\,ds.
\end{align*}
Now using  the result \eqref{estimate-eta3} in  Lemma~\ref{lem:estimate-eta} separately with $\epsilon_1=\delta$, $\iota_1=\varpi_1$ and $\epsilon_1=\delta$, $\iota_1=\varpi_2$, and Lemma~\ref{lem:estimate-xi} with $\varepsilon_0=2\delta\in(0,\epsilon)$ yields that, if there exists some $p>2$ such that $\varpi_1\in[0,\frac{1-\epsilon}{2(1+2\epsilon)(1+\epsilon)}\wedge\frac{1}{p(1+2\epsilon)})$, then for all $\varpi_2\in[0,\frac{1-\epsilon}{2(1+2\epsilon)(1+\epsilon)}\wedge\frac{p-2}{2p(1+2\epsilon)})$ such that $\varpi_1+\varpi_2\in[0,\frac{1}{2(1+2\epsilon)})$,
\begin{align*}
	|\mathbb{S}_{1,\ell}|&\le C_{L_3,\epsilon,T}\exp\big(\frac{2\|(-A)^{\frac{1-\epsilon}{4}}x\|}{e^{L_3T}}\big)\times\bigg\{\int_{t_\ell}^{t_{\ell+1}}t_{m-\ell}^{-\frac{1+\varsigma_0}{4}}(1+t_{\ell}^{-\frac{1-\epsilon}{2}-(\varpi_1+\varpi_2)}+t_{\ell-1}^{-(\varpi_1+\varpi_2)+\frac{\epsilon}{2}}\times \mathbf 1_{\{\ell>1\}}\big)\,ds\\
	& \times \big(\tau^{\frac{1-\delta}{2(1+\delta)(1+2\delta)}}\big\|(-A)^{-\frac{1-\epsilon}{4}}h_1\big\|\!+\! \big\|(-A)^{-\frac{1-\epsilon}{4}-\varpi_1}h_1\big\|\big) \big(\tau^{\frac{1-\delta}{2(1+\delta)(1+2\delta)}}\big\|(-A)^{-\frac{1-\epsilon}{4}}h_2\big\|\!+\! \big\|(-A)^{-\frac{1-\epsilon}{4}-\varpi_2}h_2\big\|\big)\\
&+\int_{t_\ell}^{t_{\ell+1}}t_{m-\ell}^{-\frac{1+\varsigma_0}{4}}\big(\big\|(-A)^{-\frac{1-\epsilon}{4}-\varpi_1}h_1\big\|\big\|(-A)^{-\frac{1-\epsilon}{4}-\varpi_2}h_2\big\|+ C_{h_1,h_2}\tau^{\frac{1-\epsilon}{2(1+\epsilon)(1+2\epsilon)}\wedge \frac{2\wedge (p-2)}{p(1+\epsilon)(1+2\epsilon)}\wedge\frac{1-\epsilon}{4}}\big)\,ds\bigg\}.
\end{align*}
In particular, under the condition on $D^2F$ in Assumption \ref{assum:nonlinearity-F} (ii) with $\varsigma_0=1-2\epsilon\in(0,1)$, we derive 
\begin{align*}
	|\mathbb{S}_{1,0}|&\le C_{\mathfrak{p}} \tau (1+(\E[\|\bar{X}^{N}_{m}\|^{\mathfrak{p}(1+1/\delta)}])^{\frac{\delta}{1+\delta}})\times\big\|E(t_m)P_ND^2F(P_Nx).(h_1,h_2)\big\|_{L^{1+\delta}(\Omega;H)}\\
	&\le C_{L_3,\epsilon,T}\exp\big(\frac{\|(-A)^{\frac{1-\epsilon}{4}}x\|}{e^{L_3T}}\big)\times \tau t_m^{-\frac{1-\epsilon}{2}}\|h_1\|\|h_2\|.
\end{align*}
Combining the estimates on terms $\mathbb{S}_{1,\ell}$, $\ell\in\{0,1,\ldots,m-1\}$ gives 
\begin{align*}
	|\mathbb{S}_1|&\le C_{L_3,\epsilon,T}\exp\big(\frac{2\|(-A)^{\frac{1-\epsilon}{4}}x\|}{e^{L_3T}}\big)\Big\{\big(1+t_m^{-\frac{1-\epsilon}{2}}\!+\!t_{m-1}^{-\frac{3+\varsigma_0-2\epsilon}{4}+1-(\varpi_1+\varpi_2)}\mathbf 1_{\{m>1\}}\!+\!t_{m-2}^{-\frac{1+\varsigma_0-2\epsilon}{4}+1-(\varpi_1+\varpi_2)}\mathbf 1_{\{m>2\}}\big)\\
&\quad \times\tau^{\frac{1-\delta}{2(1+\delta)(1+2\delta)}}\|h_1\|\|h_2\|+C_{h_1,h_2}\tau^{\frac{1-\epsilon}{2(1+\epsilon)(1+2\epsilon)}\wedge \frac{2\wedge (p-2)}{p(1+\epsilon)(1+2\epsilon)}\wedge\frac{1-\epsilon}{4}}+\|(-A)^{-\frac{1-\epsilon}{4}-\varpi_1}h_1\|\|(-A)^{-\frac{1-\epsilon}{4}-\varpi_2}h_2\|\\
&\quad \times\bigl(2+t_{m-1}^{-\frac{3+\varsigma_0-2\epsilon}{4}+1-(\varpi_1+\varpi_2)}\mathbf 1_{\{m>1\}}+t_{m-2}^{-\frac{1+\varsigma_0-2\epsilon}{4}+1-(\varpi_1+\varpi_2)}\mathbf 1_{\{m>2\}}\bigr)\Big\}\\
&\le C_{L_3,\epsilon,T}\exp\big(\frac{2\|(-A)^{\frac{1-\epsilon}{4}}x\|}{e^{L_3T}}\big) (1+t_m^{-\frac{1-\epsilon}{2}}+t_{m-1}^{-(\gamma_1+\gamma_2)}\mathbf 1_{\{m>1\}}+t_{m-2}^{-(\gamma_1+\gamma_2)}\mathbf 1_{\{m>2\}}) \\
&\quad\times \big(C_{h_1,h_2}\tau^{\frac{1-\epsilon}{2(1+\epsilon)(1+2\epsilon)}\wedge \frac{2\wedge (p-2)}{p(1+\epsilon)(1+2\epsilon)}\wedge\frac{1-\epsilon}{4}}+ \|(-A)^{-\gamma_1}h_1\|\|(-A)^{-\gamma_2}h_2\|\big),
\end{align*}
where $\gamma_1:=\frac{1-\epsilon}{4}+\varpi_1$ and $\gamma_2:=\frac{1-\epsilon}{4}+\varpi_2$, that is, $\gamma_1\in[\frac{1-\epsilon}{4},\frac{1-\epsilon}{4}+\frac{1-\epsilon}{2(1+\epsilon)(1+2\epsilon)}\wedge\frac{1}{p(1+2\epsilon)})$ and $\gamma_2\in[\frac{1-\epsilon}{4},\frac{1-\epsilon}{4}+\frac{1-\epsilon}{2(1+\epsilon)(1+2\epsilon)}\wedge\frac{p-2}{2p(1+2\epsilon)})$ where $\gamma_1+\gamma_2<\frac{1-\epsilon}{2}+\frac{1}{2(1+\epsilon)}$.

For $\ell\in\{0,1,\ldots,m-1\}$, let
$\bar X_j^{N,(t_\ell,u)}$, $j=\ell,\ldots,m$, denote the numerical
solution restarted at time $t_\ell$ from the deterministic initial
value $u\in H_N$ and driven by the Wiener--Poisson increments after
$t_\ell$. Since these future increments are independent of
$\mathcal F_{t_\ell}$, the numerical solution has the discrete-time
Markov property. In particular, for $s\in(t_\ell,t_{\ell+1})$,
\begin{align*}
\E\Big[\Big(\sum_{i\in\mathbb N_+}
\|\mathcal D_s^{W,Q^{1/2}e_i}\bar X_m^N\|^2
\Big)^{\frac{1+\delta}{2}}
\Big|\mathcal F_{t_\ell}\Big]=\E\Big[\Big(
\sum_{i\in\mathbb N_+}\|\mathcal D_s^{W,Q^{1/2}e_i}\bar X_m^{N,(t_\ell,u)}\|^2
\Big)^{\frac{1+\delta}{2}}
\Big]\Big|_{u=\bar X_\ell^N},
\end{align*}
where the equality is understood for a.e.~$s$.

\textbf{Estimate of term $\mathbb{S}_2$.} Using Proposition \ref{prop:Malliavinderivative-W},  the discrete-time Markov property of the numerical solution, and  Lemma~\ref{lem:Malliavinderivative-B1}, under Assumption~\ref{assum:nonlinearity-B} with $\varsigma_1=\frac{1-\epsilon}{4}\in[0,\frac14)$, we obtain, for $\ell\in\{1,2,\ldots,m-1\}$,
\begin{align*}
	|\mathbb{S}_{2,\ell}|&\le \Big|\sum_{i\in\N_+}\E\Big[\int_{t_\ell}^{t_{\ell+1}}\!\!\!\!D^2\varphi(\bar{X}^{N}_m).\Big(\ID^{W,Q^{\frac12}e_i}_s\bar{X}^{N}_m, \\
	&\quad \quad\quad \quad\quad \quad E((m-\ell)\tau)P_N\big(DB(\bar{X}^{N}_{\ell}).\xi^{N}_{h_1,h_2,x,\ell}\!+\!D^2B(\bar{X}^{N}_{\ell}).(\eta^{N}_{h_1,x,\ell},\eta^{N}_{h_2,x,\ell})\big)Q^{\frac12}e_i\Big)ds\Big]\Big|\\
   &\le C_{\mathfrak{p}}\big(1\!+\!\big(\E[\|\bar{X}^{N}_m\|^{\mathfrak{p}(1+1/\delta)}])^{\frac{\delta}{1+\delta}}\big)\int_{t_\ell}^{t_{\ell+1}}\!\!\Big(\E\Big[~\E\Big[\Big(\sum_{i\in\N_+}\big\|\ID^{W,Q^{\frac12}e_i}_s\bar{X}^{N}_m\big\|^2\Big)^{\frac{1+\delta}{2}}\Big|~\mathcal F_{t_{\ell}}\Big]\\
   &\quad\times \big\|E(t_{m-\ell})P_N\big(DB(\bar{X}^{N}_{\ell}).\xi^{N}_{h_1,h_2,x,\ell} +D^2B(\bar{X}^{N}_{\ell}).(\eta^{N}_{h_1,x,\ell},\eta^{N}_{h_2,x,\ell})\big)\big\|_{\IL_2^0(H)}^{1+\delta}\Big]\Big)^{\frac{1}{1+\delta}}\,ds\\
   &\le C_{\mathfrak{p}}\big(1\!+\!\big(\E[\|\bar{X}^{N}_m\|^{\mathfrak{p}(1+1/\delta)}])^{\frac{\delta}{1+\delta}}\big)\int_{t_\ell}^{t_{\ell+1}}\!\! t_{m-\ell}^{-\frac{1-\epsilon}{4}}\times \Big(\E\Big[~\E\Big[\Big(\sum_{i\in\N_+}\big\|\ID^{W,Q^{\frac12}e_i}_s\bar{X}^{N,(t_\ell,u)}_m\big\|^2\Big)^{\frac{1+\delta}{2}}\Big]\Big|_{u=\bar{X}^N_{\ell}}\\
   &\quad\times\big(\big\|(-A)^{-\frac{1-\epsilon}{4}}\xi^{N}_{h_1,h_2,x,\ell}\big\|\!+\!\|(-A)^{-\frac{1-\epsilon}{4}}\eta^{N}_{h_1,x,\ell}\|\|(-A)^{-\frac{1-\epsilon}{4}}\eta^{N}_{h_2,x,\ell}\big\|\big)^{1+\delta}\Big]\Big)^{\frac{1}{1+\delta}}ds\\
   &\le C_{\mathfrak{p}}\big(1\!+\!\big(\E[\|\bar{X}^{N}_m\|^{\mathfrak{p}(1+1/\delta)}])^{\frac{\delta}{1+\delta}}\big)\int_{t_\ell}^{t_{\ell+1}}\!\!\!t_{m-\ell}^{-\frac{1-\epsilon}{4}}\times\big(1+\|\bar{X}^N_{\ell}\|_{L^{\frac{(1+\delta)(1+2\delta)}{\delta}}(\Omega;H)}\big)\\
   &\quad \times \big(\big\|(-A)^{-\frac{1-\epsilon}{4}}\xi^{N}_{h_1,h_2,x,\ell}\big\|_{L^{1+2\delta}(\Omega;H)}+\big|\big\|(-A)^{-\frac{1-\epsilon}{4}}\eta^{N}_{h_1,x,\ell}\big\|\big\|(-A)^{-\frac{1-\epsilon}{4}}\eta^{N}_{h_2,x,\ell}\big\|\big|_{L^{1+2\delta}(\Omega;\R)}\big)ds.
\end{align*}
In particular, we obtain that 
$|\mathbb{S}_{2,0}|\le C_{\mathfrak{p}}\big(1\!+\!\big(\E[\|\bar{X}^{N}_m\|^{\mathfrak{p}(1+1/\delta)}])^{\frac{\delta}{1+\delta}}\big)\times\tau\|h_1\|\|h_2\| \|B(P_Nx)\|_{\IL_2^0(H)}.$ 
Now combining the estimates of terms $|\mathbb{S}_{2,\ell}|$, $\ell\in\{0,1,\ldots,m-1\}$, and using Corollary~\ref{cor:exponentialintegrability}, and Lemmas~\ref{lem:estimate-eta} and \ref{lem:estimate-xi},  similarly, we derive that for the same ranges of $\gamma_1$ and $\gamma_2$ as in the estimate of $\mathbb{S}_1$, it holds that
\begin{align*}
	&|\mathbb{S}_2|\le C_{\mathfrak{p}}\big(1\!+\!\big(\E[\|\bar{X}^{N}_m\|^{\mathfrak{p}(1+1/\delta)}])^{\frac{\delta}{1+\delta}}\big)\Big\{\tau\|h_1\|\|h_2\|+\!\int_{t_1}^{t_{m}}\!\!\!(t_{m}-\lfloor s\rfloor)^{-\frac{1-\epsilon}{4}}\times\big(1+\|\bar{X}^N_{\lfloor s\rfloor/\tau}\|_{L^{\frac{(1+\delta)(1+2\delta)}{\delta}}(\Omega;H)}\big)\\
   &\quad \times \big(\big\|(-A)^{-\frac{1-\epsilon}{4}}\xi^{N}_{h_1,h_2,x,\lfloor s\rfloor/\tau}\big\|_{L^{1+2\delta}(\Omega;H)}+\big|\big\|(-A)^{-\frac{1-\epsilon}{4}}\eta^{N}_{h_1,x,\lfloor s\rfloor/\tau}\big\|\big\|(-A)^{-\frac{1-\epsilon}{4}}\eta^{N}_{h_2,x,\lfloor s\rfloor/\tau}\big\|\big|_{L^{1+2\delta}(\Omega;\R)}\big)ds\Big\}\\
&\le C_{L_3,\epsilon,T}\exp\big(\frac{2\|(-A)^{\frac{1-\epsilon}{4}}x\|}{e^{L_3T}}\big)(1\!+\!t_{m-1}^{-\frac{1-\epsilon}{4}}\mathbf{1}_{\{m>1\}}) \big(C_{h_1,h_2}\tau^{\frac{1-\epsilon}{2(1+\epsilon)(1+2\epsilon)}\wedge \frac{2\wedge (p-2)}{p(1+\epsilon)(1+2\epsilon)}\wedge\frac{1-\epsilon}{4}}\!+\! \|(-A)^{-\gamma_1}h_1\|\|(-A)^{-\gamma_2}h_2\|\big).
\end{align*}

\textbf{Estimate of term $\mathbb{S}_3$.} Applying the Poisson duality formula and the difference chain rule, and then using Assumption~\ref{assum:nonlinearity-G} together with Lemmas \ref{lem:Malliavinderivative-B1}, \ref{lem:estimate-eta}, and \ref{lem:estimate-xi}, we obtain
\begin{align*}
	&|\mathbb{S}_3|\le \sum_{\ell=0}^{m-1}\int_{t_\ell}^{t_{\ell+1}}\int_\chi \Big|\E\Big[\big\langle D\varphi(\bar{X}^{N}_m+\ID^{\mathbf{N}}_{(s,z)}\bar{X}^{N}_m)-D\varphi(\bar{X}^{N}_m), \\
	&~~~\quad\quad\quad  E((m-\ell)\tau) P_N\big(DG(\bar{X}^{N}_{\ell},z).\xi^{N}_{h_1,h_2,x,\ell}\!+\!D^2G(\bar{X}^{N}_{\ell},z).(\eta^{N}_{h_1,x,\ell},\eta^{N}_{h_2,x,\ell})\big)\big\rangle\Big]\Big|\nu(dz)ds\\
&\le C_{\mathfrak{p}}\big(1\!+\!(\E[\|\bar{X}^{N}_m\|^{\mathfrak{p}(1+1/\delta)}])^{\frac{\delta}{1+\delta}}\big)\Big\{\tau\|h_1\|\|h_2\|+\int_{t_1}^{t_m}\!\!(t_m-\lfloor s\rfloor)^{-\frac{1-\epsilon}{4}} \\
&\quad \times \Big(\big\|(-A)^{-\frac{1-\epsilon}{4}}\xi^{N}_{h_1,h_2,x,\lfloor s\rfloor/\tau}\big\|_{L^{1+\delta}(\Omega;H)}+\big|\big\|(-A)^{-\frac{1-\epsilon}{4}}\eta^{N}_{h_1,x,\lfloor s\rfloor/\tau}\big\|\big\|(-A)^{-\frac{1-\epsilon}{4}}\eta^{N}_{h_2,x,\lfloor s\rfloor/\tau}\big\|\big|_{L^{1+\delta}(\Omega;\R)}\Big)ds\Big\}\\
&\le C_{L_3,\epsilon,T}\exp\big(\frac{2\|(-A)^{\frac{1-\epsilon}{4}}x\|}{e^{L_3T}}\big)(1+t_{m-1}^{-\frac{1-\epsilon}{4}}\mathbf{1}_{\{m>1\}})\big(C_{h_1,h_2}\tau^{\frac{1-\epsilon}{2(1+\epsilon)(1+2\epsilon)}\wedge \frac{2\wedge (p-2)}{p(1+\epsilon)(1+2\epsilon)}\wedge\frac{1-\epsilon}{4}}\!+\!\|(-A)^{-\gamma_1}h_1\|\|(-A)^{-\gamma_2}h_2\|\big).
\end{align*}
Therefore, combining the above results yields that for $\gamma_1\in[\frac{1-\epsilon}{4},\frac{1-\epsilon}{4}+\frac{1-\epsilon}{2(1+\epsilon)(1+2\epsilon)}\wedge\frac{1}{p(1+2\epsilon)})$ and $\gamma_2\in[\frac{1-\epsilon}{4},\frac{1-\epsilon}{4}+\frac{1-\epsilon}{2(1+\epsilon)(1+2\epsilon)}\wedge\frac{p-2}{2p(1+2\epsilon)})$ with $\gamma_1+\gamma_2<\frac{1-\epsilon}{2}+\frac{1}{2(1+\epsilon)}$,
\begin{align*}
\Big|\E\big[D\varphi(\bar{X}^{N}_m).\xi^{N}_{h_1,h_2,x,m}\big]\Big|&\le C_{L_3,\epsilon,T}\exp\big(\frac{2\|(-A)^{\frac{1-\epsilon}{4}}x\|}{e^{L_3T}}\big)\bigl(1+t_m^{-\frac{1-\epsilon}{2}}+t_{m-1}^{-(\gamma_1+\gamma_2)}\mathbf 1_{\{m>1\}}+t_{m-2}^{-(\gamma_1+\gamma_2)}\mathbf 1_{\{m>2\}}\bigr)\\
&\quad\times\big(C_{h_1,h_2}\tau^{\frac{1-\epsilon}{2(1+\epsilon)(1+2\epsilon)}\wedge \frac{2\wedge (p-2)}{p(1+\epsilon)(1+2\epsilon)}\wedge\frac{1-\epsilon}{4}}+\|(-A)^{-\gamma_1}h_1\|\|(-A)^{-\gamma_2}h_2\|\big).
\end{align*}
Using \eqref{estimate-eta3} of Lemma \ref{lem:estimate-eta} with $\epsilon_1=\delta$ gives that for any $\varphi\in\IC^2_{\mathfrak{p}}(H;\R)$ and $\gamma_1, \gamma_2\in[\frac{1-\epsilon}{4},\frac{1-\epsilon}{4}+\frac{1-\epsilon}{2(1+2\epsilon)(1+\epsilon)})$,
\begin{align*}
\!\!&\big|\E\big[D^2\varphi(\bar{X}^{N}_m).(\eta^{N}_{h_1,x,m},\eta^{N}_{h_2,x,m})\big]\big|\le C_{\mathfrak{p}}\big(1\!+\!(\E[\|\bar{X}^{N}_m\|^{\mathfrak{p}(1+1/\delta)}])^{\frac{\delta}{1+\delta}})\|\eta^{N}_{h_1,x,m}\|_{L^{2(1+\delta)}(\Omega;H)}\|\eta^{N}_{h_2,x,m}\|_{L^{2(1+\delta)}(\Omega;H)}\\
&\le C_{L_3,\epsilon,T}\exp\big(\frac{2\|(-A)^{\frac{1-\epsilon}{4}}x\|}{e^{L_3T}}\big)\big(1+t_{m}^{-(\gamma_1+\gamma_2)}+t_{m-1}^{-(\gamma_1+\gamma_2)}\mathbf{1}_{\{m>1\}}\big)\\
&\quad \times\big( \tau^{\frac{1-\epsilon}{2(1+\epsilon)(1+2\epsilon)}}\|h_1\|\|h_2\| +\|(-A)^{-\gamma_1}h_1\|\|(-A)^{-\gamma_2}h_2\|\big)
\end{align*}
Hence, it follows that
\begin{align*}
&~~~|D^2 u^{N,\tau}(t_m,x).(h_1,h_2)|=  \big|\E\big[D\varphi(\bar{X}^{N}_m). \xi^{N}_{h_1,h_2,x, m}\big]+\E\big[D^2\varphi(\bar{X}^{N}_m).(\eta^{N}_{h_1,x,m},\eta^{N}_{h_2,x,m})\big]\big|\\
&\le C_{L_3,\epsilon,T}\exp\big(\frac{2\|(-A)^{\frac{1-\epsilon}{4}}x\|}{e^{L_3T}}\big)\big(1+t_{m}^{-(\gamma_1+\gamma_2)}+t_{m-1}^{-(\gamma_1+\gamma_2)}\mathbf{1}_{\{m>1\}}+t_{m-2}^{-(\gamma_1+\gamma_2)}\mathbf 1_{\{m>2\}}\big)\\
&\quad \times\big(\|(-A)^{-\gamma_1}h_1\|\|(-A)^{-\gamma_2}h_2\|+C_{h_1,h_2}\tau^{\frac{1-\epsilon}{2(1+\epsilon)(1+2\epsilon)}\wedge \frac{2\wedge (p-2)}{p(1+\epsilon)(1+2\epsilon)}\wedge\frac{1-\epsilon}{4}}\big).
\end{align*}
In particular, based on $\varphi\in\IC^2_{\mathfrak{p}}(H;\R)$, for all $m\in\{0,1,\ldots,M\}$ and any $x,h_1,h_2\in H_N$, 
\begin{align*}
 |D^2u^{N,\tau}(t_m,x).(h_1,h_2)|&=  \big|\E\big[D\varphi(\bar{X}^{N}_m). \xi^{N}_{h_1,h_2,x, m}\big]+\E\big[D^2\varphi(\bar{X}^{N}_m).(\eta^{N}_{h_1,x,m},\eta^{N}_{h_2,x,m})\big]\big|\\
  &\le C_{\mathfrak{p}}\E\big[\big(1+\|\bar{X}^{N}_m\|^{\mathfrak{p}}\big)\times \|\xi^{N}_{h_1,h_2,x, m}\|+\big(1+\|\bar{X}^{N}_m\|^{\mathfrak{p}}\big)\times \|\eta^{N}_{h_1,x,m}\|\|\eta^{N}_{h_2,x,m}\|\big],
\end{align*}
which, combined with \eqref{eq:estimate-X} and \eqref{estimate:xi-4-}, yields \eqref{ineq:normal-second-discrete}. The proof is completed.

\subsubsection{Proof of Proposition \ref{cor:firstderivative2}}\label{subsection4.3.3}

Fix $N\in\mathbb N_+$, $x\in H_N$, and $t\in(0,T]$.
For any integer $m>\max\big\{\frac{t}{\tau_N},\,t\big\}$, set $\tau:=t/m$, so that $\tau\in(0,\tau_N)$ and $t_m=t$. We claim that 
the following estimates hold for all $\delta_0\in(0,1]$ and $p\ge 2$,
\begin{align}
	\|X^N(t)-\bar{X}^N_m\|_{L^p(\Omega;H)}&\le C_{L_3,\epsilon,T}\exp\big(\frac{\|(-A)^{\frac{1-\delta_0}{4}}x\|}{4pe^{L_3T}} \big)\times\tau^{\frac{1-\delta_0}{p}}\label{error:time},\\
	\big\|\eta^N_{y,x}(t)-\eta^N_{y,x,m}\big\|_{L^2(\Omega;H)}&\le C_{L_3,\epsilon,T}\exp\big(\frac{\|(-A)^{\frac{1-\delta_0}{4}}x\|}{8e^{L_3T}} \big)\|y\|\times\tau^{\frac{1-\delta_0}{4}}.\label{convergence:time-eta}
\end{align}
The proof is standard and omitted here. In what follows, the limit $\tau\to0$ is understood along the
sequence $\tau=t/m$ as $m\to\infty$. Recall that $u^{N,\tau}(t_m,x)=\E[\varphi(\bar{X}^{N,(0,x)}_m)]$. By Lemma \ref{lem:exponential-semidiscrete} and results \eqref{error:time} and \eqref{convergence:time-eta}, 
\begin{align*}
	\!\!\!\!\!\!\!\!\!&\lim_{\tau\to 0}\big\|(-A)^{\gamma_0}\big(D u^{N,\tau}(t_m,x)\!-\!D u^N(t_m,x)\big)\big\|_{\IL(H_N;\R)}\!\!\!=\!\lim_{\tau\to 0}\sup_{y\in H_N\setminus\{0\}}\!\!\!\!\frac{\!\!\big|D\E[\varphi(\bar{X}^{N,(0,x)}_m)].y\!-\!D\E[\varphi(X^{N,(0,x)}(t_m))].y\big|}{\|(-A)^{-\gamma_0}y\|}  \\
   &\le \lim_{\tau\to 0}\sup_{y\in H_N\setminus\{0\}}\frac{\big|\E\big[D\varphi(\bar{X}^N_m).\big(\eta^N_{y,x,m}-\eta^N_{y,x}(t_m)\big)\big]\big|}{\|(-A)^{-\gamma_0}y\|}\\
   &\quad+\lim_{\tau\to 0}\sup_{y\in H_N\setminus\{0\}}\frac{\big|\int_0^1\E\big[D^2\varphi(X^N(t_m)+\bar{\theta}(\bar{X}^N_m-X^N(t_m))). \big(\bar{X}^N_m-X^N(t_m), \eta^N_{y,x}(t_m)\big)\big]\,d\bar{\theta}\big|}{\|(-A)^{-\gamma_0}y\|}\\
	&\le C_{\mathfrak{p}}\lim_{\tau\to 0}\sup_{y\in H_N\setminus\{0\}}\big(1+\big(\E[\|\bar{X}^N_m\|^{2\mathfrak{p}}]\big)^{\frac{1}{2}}\big)\times \frac{\big\|\eta_{y,x}^N(t_m)-\eta^N_{y,x,m}\big\|_{L^{2}(\Omega;H)}}{\|(-A)^{-\gamma_0}y\|}\\
	&+C_{\mathfrak{p}}\lim_{\tau\to 0}\sup_{y\in H_N\setminus\{0\}}\!\!\big(1\!+\!\big(\E[\|\bar{X}^N_m\|^{4\mathfrak{p}}]\big)^{\frac14}\!+\!\big(\E[\|X^N(t_m)\|^{4\mathfrak{p}}]\big)^{\frac14}\big)\frac{\big\|\eta^N_{y,x}(t_m)\big\|_{L^4(\Omega;H)}}{\|(-A)^{-\gamma_0}y\|}\big\|\bar{X}^N_m\!-\!X^N(t_m)\big\|_{L^2(\Omega;H)}=0.
	\end{align*}
For $\gamma_0\ge0$, we have $\sup_{y\in H_N\setminus\{0\}}
\frac{\|y\|}
{\|(-A)^{-\gamma_0}y\|}
=
\lambda_N^{\gamma_0}<\infty$. Since the power of $\tau$ in the remainder term of  Proposition~\ref{the:expression-first} is strictly positive, 
\[
\tau^{\frac{1-\epsilon_1}
{2(1+\epsilon_1)(1+2\epsilon_1)}}
\sup_{y\in H_N\setminus\{0\}}
\frac{\|y\|}
{\|(-A)^{-\gamma_0}y\|}
\longrightarrow0,
\qquad \tau\to0.
\]
Hence, letting $\tau\to 0$,
\begin{align*}
	&\big\|(-A)^{\gamma_0}D u^{N}(t,x)\big\|_{\IL(H_N;\R)}=\lim_{\tau\to 0}\big\|(-A)^{\gamma_0}D u^{N,\tau}(t,x)\big\|_{\IL(H_N;\R)}	\le C_{L_3,\epsilon,T}\exp\big(\frac{\|(-A)^{\frac{1-\epsilon}{4}}x\|}{2e^{L_3T}}\big)(1+t^{-\gamma_0}).
\end{align*}
In particular, we obtain by \eqref{ineq:normal-first-discrete} that
$
	\big\|D u^{N}(t,x)\big\|_{\IL(H_N;\R)}=\lim_{\tau\to 0}\big\|D u^{N,\tau}(t,x)\big\|_{\IL(H_N;\R)}\le C_{\mathfrak{p},T} (1+\|x\|^{\mathfrak{p}}).
$
The proof is completed.

\subsubsection{Proof of Proposition \ref{cor-secondderivative}}\label{subsection4.3.4}

We emphasize that Proposition~\ref{the:expression-second} only requires $\varphi\in C_p^2(H;\mathbb R)$, since its proof involves only $D\varphi$ and $D^2\varphi$. The additional assumption $\varphi\in C_p^3(H;\mathbb R)$ in Proposition~\ref{cor-secondderivative} is used only in the passage $\tau\to0$, where the variation of $D^2\varphi$ has to be controlled.

Fix $N\in\mathbb N_+$, $x\in H_N$, and $t\in(0,T]$. For any integer $m>\max\big\{\frac{t}{\tau_N},\,t\big\}$, set $\tau:=t/m$, so that $\tau\in(0,\tau_N)$ and $t_m=t$. We claim that, for all $\varrho_1\in(\frac12,1)$ and $\delta\in(0,\frac12)$,
\begin{align}
        &	 \lim_{\tau\to 0}\|\xi^{N}_{h_1,h_2,x}(t)-\xi^{N}_{h_1,h_2,x, m}\|_{L^{1+\delta}(\Omega;H)}=0,\label{eq:convergence-xi}\\
		&\|(-A)^{\frac{\varrho_1}{2}}\eta^{N}_{h_2,x,m}\|_{L^2(\Omega;H)}\le C_{T}(1+t_m^{-\frac{\varrho_1}{2}})\|h_2\|\label{eq:regularity-eta}.
\end{align}
The proof is standard and is omitted. Arguing as in the proof of Proposition~\ref{cor:firstderivative2}, and using \eqref{error:time}--\eqref{eq:regularity-eta} together with the polynomial-growth bounds on $D\varphi$, $D^2\varphi$, and $D^3\varphi$, we obtain, for each fixed $N$,
\[
\lim_{\tau\to0}
\|(-A)^{\gamma_1}
\bigl(D^2u^N(t,x)-D^2u^{N,\tau}(t,x)\bigr)
(-A)^{\gamma_2}\|_{\mathcal L(H_N)}
=0.
\]
Moreover, by the equivalence of norms on the finite-dimensional
space $H_N$,
\[
K_N:=
\sup_{h_1,h_2\in H_N\setminus\{0\}}
\frac{C_{h_1,h_2}}
{\|(-A)^{-\gamma_1}h_1\|
 \|(-A)^{-\gamma_2}h_2\|}
<\infty.
\]
Denote by $\rho>0$ the exponent of $\tau$ appearing in the remainder term of \eqref{eq32}. Then
\[
K_N\tau^\rho\longrightarrow0,
\qquad \tau\to0.
\]
Consequently, letting $\tau\to0$ in Proposition~\ref{the:expression-second} yields
\[
\begin{aligned}
\|(-A)^{\gamma_1}D^2u^N(t,x)(-A)^{\gamma_2}\|_{\mathcal L(H_N)}
&=
\lim_{\tau\to0}
\|(-A)^{\gamma_1}D^2u^{N,\tau}(t,x)
(-A)^{\gamma_2}\|_{\mathcal L(H_N)}
\\
&\le
C_{L_3,\epsilon,T}
\exp\left(
\frac{2\|(-A)^{\frac{1-\epsilon}{4}}x\|}
{e^{L_3T}}
\right)
\left(1+t^{-(\gamma_1+\gamma_2)}\right).
\end{aligned}
\]
Since the constants on the right-hand side are independent of $N$, the estimate is uniform with respect to the spatial discretization parameter. In particular, we derive by \eqref{ineq:normal-second-discrete} that
\[
  \|D^2 u^{N}(t,x)\|_{\IL(H_N)}=\lim_{\tau\to 0}\|D^2 u^{N,\tau}(t,x)\|_{\IL(H_N)}\le C_{\mathfrak{p},T} (1+\|x\|^{\mathfrak{p}}). 
\]
The proof is completed.

\section{Proofs of weak convergence results}\label{sec.5}

This section presents the proof of the main result of this paper, namely, Theorem~\ref{the:weakconvergence-full}. The proof is split into two steps, based on the error between the solutions of \eqref{SPDE:MultiplicativeNoise} and the semi-discrete system \eqref{eq:spectralgalerkin}, and the error between the solutions of  \eqref{eq:spectralgalerkin} and the fully discrete scheme \eqref{eq:one-step-I}. As explained before, the main idea is to utilize the properties of the Kolmogorov equation associated with the finite-dimensional approximation \eqref{eq:spectralgalerkin}, which has been established in Section~\ref{sec:4}. Since the test functions  $\varphi$ are assumed to belong to $\IC^3_{\mathfrak{p}}(H;\R)$, the exponential integrability of the solution processes established in Corollary~\ref{cor:exponentialintegrability},  Lemmas~\ref{lem:momentbounds-fulldiscre} and~\ref{lem:exponential-semidiscrete} is also required in the proof.

Throughout this section, we fix $\epsilon>0$ sufficiently small so that all the restrictions in Propositions~\ref{cor:firstderivative2} and~\ref{cor-secondderivative} are satisfied, and choose $0<\epsilon_1<
\frac{\epsilon}{2(2-\epsilon)}$. Then $\frac{3(1-\epsilon)}4
<
\frac{1-\epsilon}{4}
+
\frac{1-\epsilon_1}
{2(1+\epsilon_1)(1+2\epsilon_1)}$, so that $\gamma_0=\frac{3(1-\epsilon)}4$ is admissible in
Proposition~\ref{cor:firstderivative2}. For convenience, since only finitely many values of $r_0$ are used
in the sequel, we fix a sufficiently large constant $L>0$ such that
\begin{align}\label{L}
L\ge
\max_{r_0\in\mathcal R}
\bigl\{\varsigma_4(r_0),\varsigma_5(r_0)\bigr\},
\end{align}
where $\mathcal R$ denotes the set of all values of $r_0$ appearing
below.

We first present the weak error estimate between solutions of \eqref{SPDE:MultiplicativeNoise} and semi-discrete system \eqref{eq:spectralgalerkin}.

\subsection{Spatial semi-discretization}

We fix an integer $N_0\ge N$ and obtain
\begin{align*}
	\big|\E[\varphi(X(T))]&-\E[\varphi(X^N(T))]\big|\le \big|\E[\varphi(X(T))]-\E[\varphi(X^{N_0}(T))]\big|+\big|\E[\varphi(X^{N_0}(T))]-\E[\varphi(X^N(T))]\big|.
\end{align*}	
By the fact that $X^{N_0}(T)$ converges to $X(T)$ in mean square sense as $N_0 \to\infty$ (see e.g. \cite[Theorem 5.1]{JA-2019}) and the uniform moment estimate of $X^{N_0}$ in Lemma~\ref{lem:exponential-semidiscrete}, we have
\begin{align}\label{limit-spatial}
	\big|\E[\varphi(X(T))]-\E[\varphi(X^N(T))]\big|\le \limsup_{N_0\to\infty}\big|\E[\varphi(X^{N_0}(T))]-\E[\varphi(X^N(T))]\big|.
\end{align}
To handle the term on the right-hand side of the above inequality,  the Kolmogorov equation \eqref{eq:kolmogorovequation} is used here, which implies that
\begin{align}\label{decompositionerror-spatial}
	\E[\varphi(X^{N_0}(T))]\!-\!\E[\varphi(X^N(T))]&=\!u^{N_0}(T,P_{N_0}x_0)\!-\!u^{N_0}(T,P_Nx_0)\nonumber\\
	&~~~+\E\big[u^{N_0}(T,P_Nx_0)\!-\!u^{N_0}(0,X^N(T))\big]
\end{align}
by employing solution $u^N(t,x)=\E[\varphi(X^{N,(0,x)}(t))]$ for all $N\in\N_+$, $t\in[0,T]$, and $x\in H_N$. Then using Taylor's formula for the first term of \eqref{decompositionerror-spatial} yields
\begin{align}\label{eq:error-onestep-Taylor}
	\big|&u^{N_0}(T,P_{N_0}x_0)-u^{N_0}(T,P_Nx_0)\big|=\Big|\int_0^1Du^{N_0}(T,P_Nx_0+\theta(P_{N_0}-P_N)x_0).(P_{N_0}-P_N)x_0\,d\theta\Big|.
\end{align}
If we take $\gamma_0=\frac{3(1-\epsilon)}{4}<\frac{1-\epsilon}{4}+\frac{1-\epsilon_1}{2(1+\epsilon_1)(1+2\epsilon_1)}$ from Proposition~\ref{cor:firstderivative2} when $\epsilon_1\in(0,\frac{\epsilon}{2(2-\epsilon)})$, then we obtain for $x_0\in\dot{H}^{\frac{1-\epsilon}{2}}$ that  
\begin{align*}
	\big|u^{N_0}(T,P_{N_0}x_0)\!-\!u^{N_0}(T,P_Nx_0)\big|\le C_{L_3,\epsilon,T}\exp\big(\frac{\|(-A)^{\frac{1-\epsilon}{4}}x_0\|}{2e^{L_3T}}\big)\big\|(-A)^{-\frac{3(1-\epsilon)}{4}}(P_{N_0}\!-\!P_N)x_0\big\|\le C_{L_3,\epsilon,T} \lambda_N^{-1+\epsilon}.
\end{align*}
For the second term on the right-hand side of \eqref{decompositionerror-spatial}, applying  It\^o's formula and Kolmogorov equation \eqref{eq:kolmogorovequation} gives 
\begin{align}\label{eq:error-secondterm-spatial}
	& \E\big[u^{N_0}(0,X^N(T))\!-\!u^{N_0}(T,P_Nx_0)\big] = \int_0^T\!\!\E\Big[D u^{N_0}(T-t, X^N(t)).(F_N(X^N(t))-F_{N_0}(X^N(t)))\Big]dt\nonumber\\
	&~~~+\frac12\int_0^T\E\Big[\sum_{i=1}^\infty D^2 u^{N_0}(T-t, X^N(t)).\big((P_N-P_{N_0})B(X^N(t))Q^{\frac12}e_i, P_NB(X^N(t))Q^{\frac12}e_i\big)\nonumber\\
	&\quad  +\sum_{i=1}^\infty D^2 u^{N_0}(T-t, X^N(t)).\big(P_{N_0}B(X^N(t))Q^{\frac12}e_i, (P_N-P_{N_0})B(X^N(t))Q^{\frac12}e_i\big)\Big]dt\nonumber\\
	&~~~+\int_0^T\int_{\chi}\!\!\E\Big[\Big(u^{N_0}(T-t, X^N(t)+G_N(X^N(t),z))-u^{N_0}(T-t, X^N(t)+G_{N_0}(X^N(t),z))\nonumber\\
	&\quad +Du^{N_0}(T-t,X^N(t)).G_{N_0}(X^N(t),z)\!-\!Du^{N_0}(T-t,X^N(t)).G_{N}(X^N(t),z)\Big)\Big]\nu(dz)dt \nonumber\\
	&=:\mathbb{J}_1+\mathbb{J}_2+\mathbb{J}_3.
\end{align}

\textbf{Estimate of term $\mathbb{J}_1$.} Indeed, under Assumptions  \ref{assum:nonlinearity-F}--\ref{assum:nonlinearity-G}, it holds that,  for all $\epsilon\in (0,1]$ and initial value $x_0\in\dot{H}^{\frac{1-\epsilon}{2}}$, there exists a  constant  $C_{\epsilon,T}>0$ such that for all $N\in\N_+$ 
	\begin{align}\label{estimate-eta5}
		\|(-A)^{\frac14}X^N(t)\|_{L^2(\Omega;H)}\le C_{\epsilon,T}(1+t^{-\frac{\epsilon}{4}}), \quad  t\in(0,T].
	\end{align}
The proof of \eqref{estimate-eta5} is postponed to Subsection~\ref{subsec:5.1.1}. Then under Assumption~\ref{assum:nonlinearity-F}(i), using Proposition~\ref{cor:firstderivative2} with $\gamma_0=\frac{3(1-\epsilon)}{4}$, Lemma~\ref{lem:exponential-semidiscrete} with $r_0=1$, we obtain a constant $C_{\epsilon,T}>0$ such that
\begin{align*}
	&|\mathbb{J}_1|\le \int_0^T\E\big[\big|D u^{N_0}(T-t, X^N(t)).((P_N-P_{N_0})F(X^N(t)))\big|\big]\,dt\\
	&\le C_{L,\epsilon,T}\int_0^T(1+(T-t)^{-\frac{3(1-\epsilon)}{4}})\times\E\Big[\exp\big(\frac{\|(-A)^{\frac{1-\epsilon}{4}}X^N(t)\|}{2e^{LT}}\big)\times\big\|(-A)^{-\frac{3(1-\epsilon)}{4}}(P_N-P_{N_0})F(X^N(t))\big\|\Big]\,dt\\
	&\le C_{L,\epsilon,T}\lambda_N^{-1+\epsilon}\sup_{t\in[0,T]}\!\Big(\E\big[\exp\big(\frac{\|(-A)^{\frac{1-\epsilon}{4}}X^N(t)\|}{e^{LT}}\big)\big]\Big)^{\frac12}\int_0^T(1+(T-t)^{-\frac{3(1-\epsilon)}{4}})\|(-A)^{\frac{1-\epsilon}{4}}F(X^N(t))\|_{L^2(\Omega;H)}\,dt\\
	&\le C_{L,\epsilon,T}\lambda_N^{-1+\epsilon}\sup_{t\in[0,T]}\!\Big(\E\big[\exp\big(\frac{\|(-A)^{\frac{1-\epsilon}{4}}X^N(t)\|}{e^{LT}}\big)\big]\Big)^{\frac12}\int_0^T(1+(T-t)^{-\frac{3(1-\epsilon)}{4}})(1+\| (-A)^{\frac14}X^N(t)\|_{L^2(\Omega;H)})\,dt \\
	&\le C_{\epsilon,T}\lambda_N^{-1+\epsilon}\exp\big(\frac12\|(-A)^{\frac{1-\epsilon}{4}}x_0\|\big)\times \int_0^T(1+(T-t)^{-\frac{3(1-\epsilon)}{4}})(1+t^{-\frac{\epsilon}{4}})\,dt\le C_{\epsilon,T}\lambda_N^{-1+\epsilon},
\end{align*}

\textbf{Estimate of term $\mathbb{J}_2$.}
Under Assumption \ref{assum:nonlinearity-B} (ii) with $\varsigma_1=0$, using Proposition~\ref{cor-secondderivative} with $p=\frac{2(1+\epsilon)}{1-\epsilon}$, $\gamma_1=\frac{3(1-3\epsilon)}{4}<\frac{1-\epsilon}{4}+\frac{1-\epsilon}{2(1+2\epsilon)(1+\epsilon)}$ and $\gamma_2=\frac{1-\epsilon}{4}$, the H\"older inequality, Lemma~\ref{lem:exponential-semidiscrete} with $r_0=8$, and Corollary \ref{cor:exponentialintegrability} gives that for $x_0\in\dot{H}^{\frac{1-\epsilon}{2}}$
\begin{align*}
	&|\mathbb{J}_2|\le\frac12 \int_0^T\E\Big[\Big|\sum_{i=1}^\infty D^2 u^{N_0}(T-t, X^N(t)).\big((P_N-P_{N_0})B(X^N(t))Q^{\frac12}e_i, P_NB(X^N(t))Q^{\frac12}e_i\big)\\
	&\quad\quad \quad +\sum_{i=1}^\infty D^2 u^{N_0}(T-t, X^N(t)).\big(P_{N_0}B(X^N(t))Q^{\frac12}e_i, (P_N-P_{N_0})B(X^N(t))Q^{\frac12}e_i\big)\Big|\Big]\,dt\\
	&\le  \int_0^T\!\E\Big[\big\|(-A)^{\frac{3(1-3\epsilon)}{4}}D^2 u^{N_0}(T-t, X^N(t))(-A)^{\frac{1-\epsilon}{4}}\big\|_{\IL(H_{N_0})} \big\|(-A)^{-\frac{3(1-3\epsilon)}{4}}(P_N\!-\!P_{N_0})(-A)^{-\frac{1-\epsilon}{4}}\big\|_{\IL(H_{N_0})}\\
	&\quad\times\sum_{i=1}^\infty \big\|(-A)^{\frac{1-\epsilon}{4}}B(X^N(t))Q^{\frac12}e_i\big\|\big\|(-A)^{-\frac{1-\epsilon}{4}}B(X^N(t))Q^{\frac12}e_i\big\|\Big]\,dt\\
	&\le C_{L,\epsilon,T}\lambda_N^{-1+3\epsilon}\int_0^T(1+(T-t)^{-1+\frac{5}{2}\epsilon})\times\big(\E\big[\exp\big(\frac{8\|(-A)^{\frac{1-\epsilon}{4}}X^N(t)\|}{e^{LT}}\big)\big]\big)^{\frac14}\\
	&~~~\times \big|\|(-A)^{\frac{1-\epsilon}{4}}B(X^N(t))\|_{\IL_2^0(H)}\big|_{L^2(\Omega;\R)}\big|\|B(X^N(t))\|_{\IL_2^0(H)}\big|_{L^4(\Omega;\R)}dt\\
	&\le C_{ \epsilon,T}\lambda_N^{-1+3\epsilon}\exp\big(3\|(-A)^{\frac{1-\epsilon}{4}}x_0\|\big)\times\int_0^T(1+(T-t)^{-1+\frac{5}{2}\epsilon}) \,dt\le C_{ \epsilon,T}\lambda_N^{-1+3\epsilon}.
\end{align*}

\textbf{Estimate of term $\mathbb{J}_3$.}
Under Assumption~\ref{assum:nonlinearity-G}, by Proposition~\ref{cor:firstderivative2}, Lemma~\ref{lem:exponential-semidiscrete} with $r_0=1$, we obtain
\begin{align*}
&	|\mathbb{J}_3|\le \int_0^T\int_{\chi}\Big|\E\Big[\Big(u^{N_0}(T-t, X^N(t)+G_N(X^N(t),z))-u^{N_0}(T-t,  X^N(t)+G_{N_0}(X^N(t),z))\\
	&~~~+Du^{N_0}(T\!-\!t,X^N(t)).G_{N_0}(X^N(t),z)-Du^{N_0}(T\!-\!t,X^N(t)).G_N(X^N(t),z)\Big)\Big]\Big|\nu(dz)dt\\
	&\le \int_0^T\!\!\int_{\chi}\Big|\E\Big[\int_0^1 Du^{N_0}\big(T\!-\!t, X^N(t)+G_N(X^N(t),z)+\theta(P_{N_0}-P_N)G(X^N(t),z)\big)\\
	&\quad .(P_{N_0}-P_N)G(X^N(t),z)\,d\theta\Big]\Big|\nu(dz)dt+\int_0^T\!\!\int_{\chi}\Big|\E\Big[Du^{N_0}(T\!-\!t,X^N(t)).(P_{N_0}-P_N)G(X^N(t),z)\Big]\Big|\,\nu(dz)dt\\
	&\le C_{L,\epsilon,T}\int_0^T\!\!(1\!+\!(T\!-\!t)^{-\frac{3(1-\epsilon)}{4}})\times \int_{\chi}\E\Big[\exp\big(\frac{\|(-A)^{\frac{1-\epsilon}{4}}X^N(t)\|+2\|(-A)^{\frac{1-\epsilon}{4}}G(X^N(t),z)\|}{2e^{LT}}\big)\\
	&\quad \times \big\|(-A)^{-\frac{3(1-\epsilon)}{4}}(P_{N_0}-P_N)(-A)^{-\frac{1}{4}}\big\|_{\IL(H_{N_0})}\big\|(-A)^{\frac{1}{4}}G(X^N(t),z)\big\|\Big]\,\nu(dz)dt\\
	&\le C_{ \epsilon,T}\exp\big(\frac12\|(-A)^{\frac{1-\epsilon}{4}}x_0\|\big)  \lambda_N^{-1+3\epsilon/4} \int_{\chi}|\bar{c}(z)|\times e^{\lambda_1^{-\frac{\epsilon}{4}}\bar{c}(z)}\nu(dz)\le C_{ \epsilon,T}\lambda_N^{-1+\frac{3\epsilon}{4}}.
\end{align*}
Combining \eqref{limit-spatial}--\eqref{eq:error-secondterm-spatial} with the above estimates yields
\[
\begin{aligned}
\left|
\E[\varphi(X(T))]
-\E[\varphi(X^N(T))]
\right|
&\le
C_{\epsilon,T}
\left(
\lambda_N^{-1+\epsilon}
+\lambda_N^{-1+3\epsilon}
+\lambda_N^{-1+\frac34\epsilon}
\right)\le
C_{\epsilon,T}\lambda_N^{-1+3\epsilon}.
\end{aligned}
\]

\subsubsection{Proof of the inequality \eqref{estimate-eta5}}\label{subsec:5.1.1}


	For all $t\in(0,T]$, under Assumptions \ref{assum:nonlinearity-F}--\ref{assum:nonlinearity-G}, it follows from \eqref{prop:semigroup1} and \eqref{eq:exponentialintegrability-exact}  with $p=2$ and $\beta_3=0$ that for $x_0\in\dot{H}^{\frac{1-\epsilon}{2}}$
    \begin{align*}
		&\|(-A)^{\frac14}X^N(t)\|_{L^2(\Omega;H)}\le C\|(-A)^{\frac14}E(t)P_Nx_0\|+C\int_0^t\|(-A)^{\frac14}E(t-s)F_N(X^N(s))\|_{L^2(\Omega;H)}\,ds\\
		&\quad+C\Big(\int_0^t\E\big[\big\|(-A)^{\frac14}E(t-s)P_NB(X^N(s))\big\|_{\IL_2^0(H)}^2\big]\,ds\Big)^{\frac12}\\
		&\quad+C\Big(\int_0^t\int_\chi\E\big[\big\|(-A)^{\frac14}E(t-s) G_N(X^N(s),z)\big\|^2\big]\,\nu(dz)ds\Big)^{\frac12}\\
		&\le C_{\epsilon}t^{-\frac{\epsilon}{4}}\|(-A)^{\frac{1-\epsilon}{4}}x_0\|+C\int_0^t(t-s)^{-\frac14}(1+\|X^N(s)\|_{L^2(\Omega;H)})\,ds\\
		&\quad+C\Big(\int_0^t(t-s)^{-\frac{1}{2}}\times\E\big[\big\|B(X^N(s))\big\|_{\IL_2^0(H)}^2\big]\,ds\Big)^{\frac12} +C\Big(\int_0^t\int_\chi |\bar{c}(z)|^2 \,\nu(dz)ds\Big)^{\frac12}\\
		&\le C_{\epsilon}t^{-\frac{\epsilon}{4}}+C_T\big(1+\sup_{s\in[0,T]}\|X^N(s)\|_{L^2(\Omega;H)}\big)\le C_{\epsilon,T}(1+t^{-\frac{\epsilon}{4}}).
	\end{align*}
The proof is completed.


We next establish the weak error estimate between solutions of   \eqref{eq:spectralgalerkin} and the fully discrete scheme \eqref{eq:one-step-I}.

\subsection{Spatio-temporal full-discretization}

The following lemmas provide, respectively, moment estimates and local temporal increment estimates for the continuous-time interpolation $\tilde X^N$, as well as a regularity estimate for the third-order derivative of the Kolmogorov solution $u^N$. They will be used repeatedly in the estimate of temporal weak error. The proofs of Lemmas~\ref{lem:localerror} and \ref{lem:thirdderivative} are postponed to Subsections~\ref{subsec5.2.2} and ~\ref{subsec5.2.3} respectively.

 \begin{lemma}\label{lem:localerror}
 	Let Assumptions \ref{assum:nonlinearity-F}--\ref{assum:nonlinearity-G} hold. Then for all $\epsilon\in(0,\frac12)$ and $x_0\in \dot{H}^{\frac{1-\epsilon/2}{2}}$, there exists a constant $C_{\epsilon,T}>0$ such that for all $N\in\N_+$ and $\tau\in(0, \tau_N)$  with $\tau_N$ given in \eqref{def:tau_0},
     \begin{align*}
 		\|(-A)^{\frac{3-\epsilon}{4}}\tilde{X}^N(t)\|_{L^2(\Omega;H)}&\le C_{\epsilon,T}(1+t^{-\frac{1}{2}} ),\quad t\in(0,T],\\ 
 		 	\int_{\lfloor t\rfloor}^t \big\|(-A)^{\frac{3}{4}}\tilde{X}^N(s)\big\|_{L^2(\Omega;H)}\,ds&\le C_{\epsilon,T}\tau^{1-\frac{\epsilon}{4}}(1+\lfloor t \rfloor^{-\frac12}),\quad t\in[t_1,T].
 	\end{align*}
 	Moreover, for all $\delta\in(0,\frac12)$, there exists  a constant $C_{\epsilon,T}>0$ such that for all $N\in\N_+$ and $m\in\{0,\ldots,M-1\}$,
 	\begin{align*}
 		\|(-A)^{-\frac{1-\epsilon}{4}}(\tilde{X}^N(t)\!-\!\bar{X}^N_m)\|_{L^{2(1+\delta)}(\Omega;H)}&\le C_{\epsilon,T}\tau^{\frac{1}{2(1+\delta)}\wedge \frac{1-\epsilon}{2}},\quad t\in(t_m,t_{m+1}].
 	\end{align*}
 \end{lemma}

 \begin{lemma}\label{lem:thirdderivative}
 	Let Assumptions~\ref{assum:nonlinearity-F}--\ref{assum:nonlinearity-G} hold and $\varphi\in\IC^3_{\mathfrak{p}}(H;\R)$. Then for the sufficiently small constant $\epsilon > 0$ fixed at the beginning of Section 5 and for $\varrho_1\in(\frac12,1)$ given by \eqref{setting-varrho1}, there exist  constants $L>0$ and $C_{L,\epsilon,T}>0$ such that for all $N\in\N_+$ and $x,y,h,k\in H_N$, 
 	\begin{align*}
 		|D^3u^N(t,x).(y,h,k)|\le C_{L,\epsilon,T}\exp\big(\frac{\|(-A)^{\frac{1-\epsilon}{4}}x\|}{4e^{LT}}\big)(1+t^{-\frac{\varrho_1}{2}}) \|y\|\|h\|\|k\|,\quad t\in(0,T].
 	\end{align*}
 \end{lemma}


\subsubsection{Estimates of error terms}

Applying It\^o's formula and the Kolmogorov equation~\eqref{eq:kolmogorovequation} gives 
\begin{align}\label{decompositionerror-temporal}
	&~~~\E[\varphi(X^N(T))]-\E[\varphi(\bar{X}^N_M)]=\E\big[u^N(T,P_Nx_0)-u^N(0,\bar{X}^N_M)\big]\nonumber\\
	&=\sum_{m=0}^{M-1}\E\big[u^N(T\!-\!t_{m},\bar{X}^N_{m})-u^N(T\!-\!t_{m+1},\bar{X}^N_{m+1})\big]=:-\sum_{m=0}^{M-1}(\mathbb{H}_{m,1}+\mathbb{H}_{m,2}+\mathbb{H}_{m,3})
\end{align}
with the terms $\mathbb{H}_{m,1}$, $\mathbb{H}_{m,2}$, and $\mathbb{H}_{m,3}$ defined by
\begin{align*}
\mathbb{H}_{m,1}&=	\int_{t_m}^{t_{m+1}}\E\big[ Du^N(T-t,\tilde{X}^N(t)).\big(E(t-t_m)F_N(\bar{X}^N_{m})-F_N(\tilde{X}^N(t))\big)\big]dt,\\
\mathbb{H}_{m,2}&=\frac12\int_{t_m}^{t_{m+1}}\Big(\E\big[\sum_{i=1}^\infty D^2u^N(T-t,\tilde{X}^N(t)).\big(E(t-t_m)P_NB(\bar{X}^N_{m})Q^{\frac12}e_i,E(t-t_m)P_NB(\bar{X}^N_{m})Q^{\frac12}e_i\big)\big]\nonumber\\
	&~-\E\big[\sum_{i=1}^\infty D^2u^N(T-t,\tilde{X}^N(t)).\big(P_NB(\tilde{X}^N(t))Q^{\frac12}e_i,P_NB(\tilde{X}^N(t))Q^{\frac12}e_i\big)\big]\Big)dt,\\
\mathbb{H}_{m,3}&=\int_{t_m}^{t_{m+1}}\int_\chi\Big\{\E\Big[u^N\big(T-t,\tilde{X}^N(t)+E(t-t_m)G_N(\bar{X}^N_m,z)\big)\nonumber\\
	&~ -\!Du^N(T-t,\tilde{X}^N(t)).E(t-t_m)G_N(\bar{X}^N_m,z)\Big]-\E\Big[u^N\big(T-t, \tilde{X}^N(t)+G_N(\tilde{X}^N(t),z)\big)\nonumber\\
	&~-Du^N(T-t,\tilde{X}^N(t)).G_N(\tilde{X}^N(t),z)\Big]\Big\}\nu(dz)dt.	
\end{align*}

\textbf{Estimate of term $\mathbb{H}_{m,1}$.}
We  further decompose $\mathbb{H}_{m,1}$ as follows:
\begin{align*}
\mathbb{H}_{m,1}
	&=-\int_{t_m}^{t_{m+1}}\E\big[ Du^N(T-t,\tilde{X}^N(t)).\big(\mathrm I-E(t-t_m)\big)F_N(\bar{X}^N_m)\big]dt\\
	&~~~-\int_{t_m}^{t_{m+1}}\E\big[ \big(Du^N(T-t,\tilde{X}^N(t))-Du^N(T-t,\bar{X}^N_m)\big).\big(F_N(\tilde{X}^N(t))-F_N(\bar{X}^N_m)\big)\big]dt\\
	&\quad-\int_{t_m}^{t_{m+1}}\E\big[ Du^N(T-t,\bar{X}^N_m).\big(F_N(\tilde{X}^N(t))-F_N(\bar{X}^N_m)\big)\big]dt=:\mathbb H_{m,1,1}+\mathbb H_{m,1,2}+\mathbb H_{m,1,3}.
\end{align*}

\textit{Estimate of  subterm $\mathbb{H}_{m,1,1}$}. 
Under Assumption \ref{assum:nonlinearity-F} (i) with $r_1=\frac{1-\epsilon}{4}$ and $\varrho=\frac{\epsilon}{8}$, using Proposition~\ref{cor:firstderivative2} with $\gamma_0=\frac{3(1-\epsilon)}{4}$ and $\gamma_0=\frac12$, respectively, Lemma~\ref{lem:momentbounds-fulldiscre} with $r_0=1$,  and  Lemma~\ref{lem:localerror} yields 
\begin{align*}
	&\sum_{m=0}^{M-1}|\mathbb{H}_{m,1,1}|\le \sum_{m=0}^{M-1}\int_{t_m}^{t_{m+1}}\!\!\!\E\big[ \big|\big\langle (-A)^{\frac{3(1-\epsilon)}{4}}Du^N(T\!-\!t,\tilde{X}^N(t)), (-A)^{-\frac{3(1-\epsilon)}{4}}\big(\mathrm I\!-\!E(t-t_m)\big)F_N(\bar{X}^N_m)\big\rangle\big|\big]\,dt\\
	&\le C_{L,\epsilon,T}\int_{t_1}^{T}\!\!\big(1\!+\!(T\!-\!t)^{-\frac{3(1-\epsilon)}{4}}\big)\times \E\Big[\exp\big(\frac{\|(-A)^{\frac{1-\epsilon}{4}}\tilde{X}^N(t)\|}{2e^{LT}}\big)\times\tau^{1-\epsilon}\big\|(-A)^{\frac{1-\epsilon}{4}}F(\bar{X}^N_{\lfloor t\rfloor/\tau})\big\|\Big]\,dt\\
	&\quad +C_{L,\epsilon,T}\int_0^{t_1}\!\!\big(1\!+\!(T\!-\!t)^{-\frac12}\big)\times\tau^{\frac12}\|F(P_Nx_0)\|\times \E\big[\exp\big(\frac{\|(-A)^{\frac{1-\epsilon}{4}}\tilde{X}^N(t)\|}{2e^{LT}}\big)\big]\,dt\\
	&\le C_{\epsilon,T}\exp\big(\frac12\|(-A)^{\frac{1-\epsilon}{4}}x_0\|\big)\Big(\tau^{1-\epsilon}\int_{t_1}^{T}\!\!\big(1\!+\!(T\!-\!t)^{-\frac{3(1-\epsilon)}{4}}\big)\bigl(1+\|(-A)^{\frac{1-\epsilon/2}{4}}\bar{X}^N_{\lfloor t\rfloor/\tau}\|_{L^2(\Omega;H)}\bigr)\,dt+\tau\Big)\\
	&\le C_{\epsilon,T}\tau^{1-\epsilon}\int_{t_1}^{T}\!\!\big(1\!+\!(T\!-\!t)^{-\frac{3(1-\epsilon)}{4}}\big)\,dt \le C_{\epsilon,T}\tau^{1-\epsilon}.
\end{align*}

\textit{Estimate of subterm $\mathbb{H}_{m,1,2}$}. 
Applying  the mean value theorem, we have
\begin{align*}
	\mathbb H_{m,1,2}&=-\int_{t_m}^{t_{m+1}}\int_0^1\E\big[ D^2u^N(T\!-\!t,\bar{X}^N_m+\theta(\tilde{X}^N(t)\!-\!\bar{X}^N_m)).\big(\tilde{X}^N(t)\!-\!\bar{X}^N_m, F_N(\tilde{X}^N(t))\!-\!F_N(\bar{X}^N_m)\big)\big]\,d\theta dt.
\end{align*}
Using \eqref{eq30} with $\gamma_1=\frac{1-\epsilon}{4}$ and $\gamma_2=\frac{1+\epsilon}{4}$, whose admissibility follows from the above choice of $\epsilon$, together with \eqref{ineq:normal-second}, and applying Lemma~\ref{lem:momentbounds-fulldiscre} with $r_0=\frac{3(1+\delta)}{\delta}=\frac3\epsilon$, $\delta:=\frac{\epsilon}{1-\epsilon}\in(0,1/2),$ whose  constant $L$ is fixed in \eqref{L}, and Lemma~\ref{lem:localerror}, we obtain
\begin{align*}
	&\!\!\sum_{m=0}^{M-1}|\mathbb{H}_{m,1,2}|\!\le\! \sum_{m=0}^{M-1}\!\int_{t_m}^{t_{m+1}}\!\!\!\int_0^1\E\big[ \big|D^2u^N(T\!-\!t,\bar{X}^N_m\!+\!\theta(\tilde{X}^N(t)\!-\!\bar{X}^N_m)).\big(\tilde{X}^N(t)\!-\!\bar{X}^N_m, F_N(\tilde{X}^N(t))\!-\!F_N(\bar{X}^N_m)\big)\big|\big]d\theta dt\\
	&\le C_{L,\epsilon,T}\int_{t_1}^{T}\!\big(1\!+\!(T\!-\!t)^{-\frac{1}{2}}\big)\E\Big[\exp\big(\frac{2(\|(-A)^{\frac{1-\epsilon}{4}}\tilde{X}^N(t)\|+\|(-A)^{\frac{1-\epsilon}{4}}\bar{X}^N_{\lfloor t\rfloor/\tau}\|)}{e^{LT}}\big)\big\|(-A)^{-\frac{1-\epsilon}{4}}(\tilde{X}^N(t)\!-\!\bar{X}^N_{\lfloor t\rfloor/\tau})\big\|^2\\
	& \times\big(1\!+\!\|(-A)^{\frac{1-\epsilon}{4}}\tilde{X}^N(t)\|\!+\!\|(-A)^{\frac{1-\epsilon}{4}}\bar{X}^N_{\lfloor t\rfloor/\tau}\|\big)\Big]\,dt +C_{\mathfrak{p},T}\int_0^{t_1}\E\Big[(1+\|\tilde{X}^N(t)\|^{\mathfrak{p}}+\|x_0\|^{\mathfrak{p}}) \big\|\tilde{X}^N(t)-P_Nx_0\big\|^2\Big]\,dt\\
	&\le C_{\epsilon,T}\Big(\tau+\int_{t_1}^{T}\big(1\!+\!(T\!-\!t)^{-\frac{1}{2}}\big)\times\big\|(-A)^{-\frac{1-\epsilon}{4}}(\tilde{X}^N(t)\!-\!\bar{X}^N_{\lfloor t\rfloor/\tau})\big\|^2_{L^{2(1+\delta)}(\Omega;H)}\,dt\Big)\le C_{\epsilon,T}\tau^{1-\epsilon}.
\end{align*}

\textit{Estimate of subterm $\mathbb{H}_{m,1,3}$}. 
Using It\^o's formula, we decompose term $\mathbb{H}_{m,1,3}$ as the following terms
\begin{align*}
	&\mathbb{H}_{m,1,3}=-\int_{t_m}^{t_{m+1}}\E\Big[ \sum_{i=1}^N\langle Du^N(T-t,\bar{X}^N_m),e_i\rangle\int_{t_m}^tDF_i(\tilde{X}^N(s)).\big(A\tilde{X}^N(s)+E(s-t_m)F_N(\bar{X}^N_m)\big)\,ds\Big]\,dt\\
	&-\int_{t_m}^{t_{m+1}}\E\Big[ \sum_{i=1}^N\langle Du^N(T-t,\bar{X}^N_m),e_i\rangle\int_{t_m}^tDF_i(\tilde{X}^N(s)).E(s-t_m)P_NB(\bar{X}^N_m)\,dW^Q(s)\Big]\,dt\\
	&-\int_{t_m}^{t_{m+1}}\!\!\E\Big[ \sum_{i=1}^N\langle Du^N(T\!-\!t,\bar{X}^N_m),e_i\rangle\int_{t_m}^t\!\!\int_{\chi}\big(F_i(\tilde{X}^N(s)\!+\!E(s-t_m)G_N(\bar{X}^N_m,z))\!-\!F_i(\tilde{X}^N(s))\big)\tilde{\mathbf{N}}(dz,ds)\Big]dt\\
	&-\frac12\int_{t_m}^{t_{m+1}}\!\E\Big[ \sum_{i=1}^N\langle Du^N(T\!-\!t,\bar{X}^N_m),e_i\rangle\int_{t_m}^t\sum_{j\in\N_+}D^2F_i(\tilde{X}^N(s)).\Big(E(s\!-\!t_m)P_NB(\bar{X}^N_m)Q^{\frac12}e_j,\\
	&~~~E(s\!-\!t_m)P_NB(\bar{X}^N_m)Q^{\frac12}e_j\Big)ds\Big]dt-\int_{t_m}^{t_{m+1}}\E\Big[ \sum_{i=1}^N\langle Du^N(T\!-\!t,\bar{X}^N_m),e_i\rangle\int_{t_m}^t\int_{\chi}\Big\{\\
	&~~~F_i\big(\tilde{X}^N(s)+E(s-t_m)G_N(\bar{X}^N_m,z)\big)-F_i(\tilde{X}^N(s))-DF_i(\tilde{X}^N(s)).E(s-t_m)G_N(\bar{X}^N_m,z)\Big\}\nu(dz)ds\Big]\,dt\\
	&=:\mathbb{H}_{m,1,3,1}+\mathbb{H}_{m,1,3,2}+\mathbb{H}_{m,1,3,3}+\mathbb{H}_{m,1,3,4}+\mathbb{H}_{m,1,3,5},
\end{align*}
where $F_i:=\langle F,e_i\rangle$ and we have $\mathbb{H}_{m,1,3,2}\equiv\mathbb{H}_{m,1,3,3}\equiv0$ due to the martingale property of stochastic Wiener and Poisson integrals.  Since $\tilde X^N(s)\in H_N$, the spectral projection gives $\|(-A)^{\frac{3+\epsilon}{4}}\tilde X^N(s)\|
\le
\lambda_N^{\epsilon/2}
\|(-A)^{\frac{3-\epsilon}{4}}\tilde X^N(s)\|.$ Hence, by Lemma~\ref{lem:localerror}, we derive
\[
\|(-A)^{\frac{3+\epsilon}{4}}\tilde X^N(s)\|_{L^2(\Omega;H)}
\le
C_{\epsilon,T}
\lambda_N^{\epsilon/2}(1+s^{-1/2}).
\]
Then using Proposition~\ref{cor:firstderivative2} with $\gamma_0=\frac{1+\epsilon}{4}$ and $\gamma_0=0$ respectively, together with  Lemma~\ref{lem:momentbounds-fulldiscre} with $r_0=3$,  and the CFL condition $\tau\lambda_N\le \frac{\sqrt{\ln 2}}{4}$, we obtain
\begin{align*}
	&\!\!\sum_{m=0}^{M-1}|\mathbb{H}_{m,1,3,1}|\!\le\! \sum_{m=0}^{M-1}\int_{t_m}^{t_{m+1}}\!\!\E\Big[\Big| \sum_{i=1}^N\langle Du^N(T\!-\!t,\bar{X}^N_m),e_i\rangle\int_{t_m}^t\!\!\!DF_i(\tilde{X}^N(s)).\big(A\tilde{X}^N(s)\!+\!E(s\!-\!t_m)F_N(\bar{X}^N_m)\big)\,ds\Big|\Big]\,dt\\
	&= \int_{0}^{T}\E\Big[\Big| \int_{\lfloor t\rfloor}^t\big\langle Du^N(T-t,\bar{X}^N_{\lfloor t\rfloor/\tau}), P_NDF(\tilde{X}^N(s)).\big(A\tilde{X}^N(s)+E(s-\lfloor t\rfloor)F_N(\bar{X}^N_{\lfloor t\rfloor/\tau})\big)\big\rangle\,ds\Big|\Big]\,dt\\
	&\le C_{L,\epsilon, T}\int_{0}^{T}(1+(T\!-\!t)^{-\frac{1+\epsilon}{4}})\times \int_{\lfloor t\rfloor}^t\E\Big[\exp\big(\frac{\|(-A)^{\frac{1-\epsilon}{4}}\bar{X}^N_{\lfloor t\rfloor/\tau}\|}{2e^{LT}}\big)\\
	&\quad \times\Big\{\big(1+\big\|(-A)^{\frac{1-\epsilon}{4}}\tilde{X}^N(s)\big\|\big)\times\big\|(-A)^{\frac{3+\epsilon}{4}}\tilde{X}^N(s)\big\|+\|F_N(\bar{X}^N_{\lfloor t\rfloor/\tau})\|\Big\}\Big]\,dsdt\\
	&\le C_{\epsilon,T}\exp\big(\|(-A)^{\frac{1-\epsilon}{4}}x_0\|\big)\Big(\tau+\lambda_N^{\frac{\epsilon}{2}}\int_{t_1}^{T}(1+(T\!-\!t)^{-\frac{1+\epsilon}{4}}) \int_{\lfloor t\rfloor}^t(1+s^{-\frac{1}{2}})\,dsdt\Big)\le C_{\epsilon,T}\tau^{1-\epsilon/2}.
\end{align*}
Under Assumptions $\ref{assum:nonlinearity-F}$ and \ref{assum:nonlinearity-B}, using Proposition~\ref{cor:firstderivative2} with $\gamma_0=\frac{1+\varsigma_0}{4}$, $\varsigma_0\in(0,1)$ and Lemma~\ref{lem:momentbounds-fulldiscre} with $r_0=1$ gives
\begin{align*}
	&\!\!\sum_{m=0}^{M-1}|\mathbb{H}_{m,1,3,4}|\!\le\! \sum_{m=0}^{M-1}\int_{t_m}^{t_{m+1}}\!\E\Big[\Big| \sum_{i=1}^N\langle Du^N(T\!-\!t,\bar{X}^N_m),e_i\rangle\int_{t_m}^t\!\sum_{j\in\N_+}D^2F_i(\tilde{X}^N(s)).\Big(E(s\!-\!t_m)P_NB(\bar{X}^N_m)Q^{\frac12}e_j,\\
	&~~~\quad\quad\quad \quad\quad \quad \quad E(s\!-\!t_m)P_NB(\bar{X}^N_m)Q^{\frac12}e_j\Big)ds\Big|\Big]\,dt\\
	&\le \int_0^T\!\!\int_{\lfloor t\rfloor}^t\!\!\E\Big[\Big|\sum_{j\in\N_+}\big\langle Du^N(T\!-\!t,\bar{X}^N_{\lfloor t\rfloor/\tau}), \\
    &~~~~\quad\quad\quad P_ND^2F(\tilde{X}^N(s)).\big(E(s\!-\lfloor t\rfloor)P_NB(\bar{X}^N_{\lfloor t\rfloor/\tau})Q^{\frac12}e_j, E(s\!-\!\lfloor t\rfloor)P_NB(\bar{X}^N_{\lfloor t\rfloor/\tau})Q^{\frac12}e_j\big)\big\rangle\Big|\Big]\,dsdt\\
	&\le C_{L,\epsilon,T}\int_0^T\!\!(1+(T-t)^{-\frac{1+\varsigma_0}{4}})\times \int_{\lfloor t\rfloor}^t\E\Big[\exp\big(\frac{\|(-A)^{\frac{1-\epsilon}{4}}\bar{X}^N_{\lfloor t\rfloor/\tau}\|}{2e^{LT}}\big)\times \sum_{j\in\N_+}\big\|E(s\!-\lfloor t\rfloor)P_NB(\bar{X}^N_{\lfloor t\rfloor/\tau})Q^{\frac12}e_j\big\|^2\Big]\,dsdt\\
	&\le C_{\epsilon,T}\exp\big(\|(-A)^{\frac{1-\epsilon}{4}}x_0\|\big)\tau\le C_{\epsilon,T}\tau.
\end{align*}
By the global bound on $DF$ in Assumption~\ref{assum:nonlinearity-F} (ii),
\[
\|F(x+h)-F(x)-DF(x).h\|
\le
\|F(x+h)-F(x)\|+\|DF(x).h\|
\le
2L_F\|h\|.
\]
Using Proposition~\ref{cor:firstderivative2} with $\gamma_0=0$ and  Lemma~\ref{lem:momentbounds-fulldiscre} with $r_0=\frac12$  yields 
\begin{align*}
	&\sum_{m=0}^{M-1}|\mathbb{H}_{m,1,3,5}|\le \sum_{m=0}^{M-1}\int_{t_m}^{t_{m+1}}\E\Big[\Big| \sum_{i=1}^N\langle Du^N(T\!-\!t,\bar{X}^N_m),e_i\rangle\int_{t_m}^t\int_{\chi}\Big\{F_i\big(\tilde{X}^N(s)+E(s-t_m)G_N(\bar{X}^N_m,z)\big)\\
	&~~~-F_i(\tilde{X}^N(s))-DF_i(\tilde{X}^N(s)).(E(s-t_m)G_N(\bar{X}^N_m,z))\Big\}\nu(dz)ds\Big|\Big]\,dt\\
	&\le \int_{0}^{T} \int_{\lfloor t\rfloor}^t\int_{\chi}\E\Big[\Big|\Big\langle Du^N(T\!-\!t,\bar{X}^N_{\lfloor t\rfloor/\tau}), P_N\Big\{F\big(\tilde{X}^N(s)+E(s-\lfloor t\rfloor)G_N(\bar{X}^N_{\lfloor t\rfloor/\tau},z)\big)-F(\tilde{X}^N(s))\\
	&~~~-DF(\tilde{X}^N(s)).(E(s-\lfloor t\rfloor)G_N(\bar{X}^N_{\lfloor t\rfloor/\tau},z))\Big\}\Big\rangle\Big|\Big] \,\nu(dz)ds\,dt\\
	&\le C_{L,\epsilon,T}\int_{0}^{T} \int_{\lfloor t\rfloor}^t\int_{\chi}\E\Big[\exp\big(\frac{\|(-A)^{\frac{1-\epsilon}{4}}\bar{X}^N_{\lfloor t\rfloor/\tau}\|}{2e^{LT}}\big)\times\|G(\bar{X}^N_{\lfloor t\rfloor/\tau},z)\|\Big]\,\nu(dz)dsdt \le C_{\epsilon,T}\tau.
\end{align*}
Hence combining the above results gives $\sum_{m=0}^{M-1}|\mathbb{H}_{m,1}|\le C_{\epsilon,T}\tau^{1-\epsilon}$.

\textbf{Estimate of term $\mathbb{H}_{m,2}$.} Using the symmetry of $D^2u^N$, we decompose $\mathbb{H}_{m,2}$ as follows:
\begin{align*}
&	\mathbb{H}_{m,2}
	=\Big(\frac12\int_{t_m}^{t_{m+1}}\E\big[\sum_{i=1}^\infty D^2u^N(T-t,\tilde{X}^N(t)).\big(E(t-t_m)P_NB(\bar{X}^N_{m})Q^{\frac12}e_i,E(t-t_m)P_NB(\bar{X}^N_{m})Q^{\frac12}e_i\big)\big]dt\\
	&~~~-\frac12\int_{t_m}^{t_{m+1}}\E\big[\sum_{i=1}^\infty D^2u^N(T-t,\tilde{X}^N(t)).\big(E(t-t_m)P_NB(\tilde{X}^N(t))Q^{\frac12}e_i,E(t-t_m)P_NB(\tilde{X}^N(t))Q^{\frac12}e_i\big)\big]dt\Big)\\
&~~~-\frac12\int_{t_m}^{t_{m+1}}\E\big[\sum_{i=1}^\infty D^2u^N(T-t,\tilde{X}^N(t)).\big((\mathrm I-E(t-t_{m}))P_NB(\tilde{X}^N(t))Q^{\frac12}e_i,\\
	&\quad\quad\quad\quad\quad (\mathrm I-E(t-t_{m}))P_NB(\tilde{X}^N(t))Q^{\frac12}e_i\big)\big]dt\\
	&~~~-\int_{t_m}^{t_{m+1}}\E\big[\sum_{i=1}^\infty D^2u^N(T\!-\!t,\tilde{X}^N(t)).\big(E(t\!-\!t_m)P_NB(\tilde{X}^N(t))Q^{\frac12}e_i,(\mathrm I\!-\!E(t\!-\!t_{m}))P_NB(\tilde{X}^N(t))Q^{\frac12}e_i\big)\big]dt\\
	&=:\mathbb{H}_{m,2,1}+\mathbb{H}_{m,2,2}+\mathbb{H}_{m,2,3}.
\end{align*}

\textit{Estimate of subterm $\mathbb{H}_{m,2,2}$.}
Using Proposition~\ref{cor-secondderivative} with $\gamma_1=\gamma_2=\frac{1-\epsilon}{4}$ and Lemma~\ref{lem:momentbounds-fulldiscre} with $r_0=3$ yields 
\begin{align*}
&\sum_{m=0}^{M-1}|\mathbb{H}_{m,2,2}|\!=\!\sum_{m=0}^{M-1}\frac12\Big|\int_{t_m}^{t_{m+1}}\E\Big[\sum_{i=1}^\infty D^2u^N(T-t,\tilde{X}^N(t)).\big((\mathrm I-E(t-t_{m}))P_NB(\tilde{X}^N(t))Q^{\frac12}e_i,\\
	&\quad\quad\quad\quad\quad (\mathrm I-E(t-t_{m}))P_NB(\tilde{X}^N(t))Q^{\frac12}e_i\big)\Big]dt\Big|\\
    &\le C_{L, \epsilon,T}\int_{0}^{T}(1+(T-t)^{-\frac{1-\epsilon}{2}})\times  \E\Big[\exp\big(\frac{2\|(-A)^{\frac{1-\epsilon}{4}}\tilde{X}^N(t)\|}{e^{LT}}\big)\\
	&\quad \times\|(-A)^{-\frac{1-\epsilon}{2}}(\mathrm I-E(t-\lfloor t\rfloor))\|^2_{\IL(H)}\times \sum_{i=1}^\infty\|(-A)^{\frac{1-\epsilon}{4}}B(\tilde{X}^N(t))Q^{\frac12}e_i\|^2\Big]\,dt\\
	&\le C_{L, \epsilon,T}\tau^{1-\epsilon}\int_{0}^{T}(1+(T-t)^{-\frac{1-\epsilon}{2}})\times \E\Big[\exp\big(\frac{2\|(-A)^{\frac{1-\epsilon}{4}}\tilde{X}^N(t)\|}{e^{LT}}\big)\times\big(1+\|(-A)^{\frac{1-\epsilon}{4}}\tilde{X}^N(t)\|^{2\gamma}\big)\Big]\,dt\\
	&\le C_{\epsilon,T}\tau^{1-\epsilon}\times\exp\big(3\|(-A)^{\frac{1-\epsilon}{4}}x_0\|\big)\le C_{\epsilon,T}\tau^{1-\epsilon}.
\end{align*}

\textit{Estimate of subterm $\mathbb{H}_{m,2,3}$.}
From Proposition~\ref{cor-secondderivative} with $p=\frac{2(1+\epsilon)}{1-\epsilon}$, $\gamma_1=\frac{3(1-3\epsilon)}{4}<\frac{1-\epsilon}{4}+\frac{1-\epsilon}{2(1+2\epsilon)(1+\epsilon)}$ and $\gamma_2=\frac{1-\epsilon}{4}$, it follows that
\begin{align*}
&\sum_{m=0}^{M-1}\!|\mathbb{H}_{m,2,3}|=\sum_{m=0}^{M-1}\Big|\int_{t_m}^{t_{m+1}}\!\!\!\E\Big[\sum_{i=1}^\infty D^2u^N(T-t,\tilde{X}^N(t)).\big(E(t-t_m)P_NB(\tilde{X}^N(t))Q^{\frac12}e_i,\\
	&\quad\quad\quad\quad\quad\quad\quad\quad\quad (\mathrm I-E(t-t_m))P_NB(\tilde{X}^N(t))Q^{\frac12}e_i\big)\Big]dt\Big|\\
	&\le C_{L,\epsilon,T}\int_{0}^{T}\!\!\!\big(1+(T-t)^{-1+\frac{5\epsilon}{2}}\big)\times\E\Big[\exp\big(\frac{2\|(-A)^{\frac{1-\epsilon}{4}}\tilde{X}^N(t)\|}{e^{LT}}\big)\times \sum_{i=1}^\infty\big\|(-A)^{-\frac{1-\epsilon}{4}}P_NB(\tilde{X}^N(t))Q^{\frac12}e_i\big\|\\
	&\quad \times \bigl\|(-A)^{-\frac{3(1-3\epsilon)}{4}}(\mathrm I-E(t-\lfloor t\rfloor))(-A)^{-\frac{1-\epsilon}{4}}\bigr\|_{\IL(H)}\|(-A)^{\frac{1-\epsilon}{4}}B(\tilde{X}^N(t))Q^{\frac12}e_i\big\|\Big]dt\\
	&\le C_{L,\epsilon,T}\tau^{1-3\epsilon}\int_0^T(1+(T-t)^{-1+\frac{5\epsilon}{2}})\times\big(\E\big[\exp\big(\frac{8\|(-A)^{\frac{1-\epsilon}{4}}\tilde{X}^N(t)\|}{e^{LT}}\big)\big]\big)^{\frac14}\\
	&\quad\times \big|\|(-A)^{\frac{1-\epsilon}{4}}B(\tilde{X}^N(t))\|_{\IL_2^0(H)}\big|_{L^2(\Omega;\R)}\big|\|B(\tilde{X}^N(t))\|_{\IL_2^0(H)}\big|_{L^4(\Omega;\R)}dt\\
	&\le C_{\epsilon,T}\tau^{1-3\epsilon}\exp\big(2\|(-A)^{\frac{1-\epsilon}{4}}x_0\|\big)\int_0^T(1+(T-t)^{-1+\frac{5\epsilon}{2}})\,dt\le C_{\epsilon,T}\tau^{1-3\epsilon}.
\end{align*}

\textit{Estimate of subterm $\mathbb{H}_{m,2,1}$.}
By the symmetry of $D^2u^N$ and the triangle inequality $\sum_{m=0}^{M-1}|\mathbb H_{m,2,1}|
\le \mathcal R_{211}+\mathcal R_{212}$, where
\begin{align*}
\mathcal R_{211}
:={}&\int_{0}^{T} \Big|\E\Big[\sum_{i=1}^\infty D^2u^N(T-t,\tilde{X}^N(t)).\big(E(t-\lfloor t\rfloor)P_N(B(\tilde{X}^N(t))-B(\bar{X}^N_{\lfloor t\rfloor/\tau}))Q^{\frac12}e_i,\\
	&\quad  E(t-\lfloor t\rfloor)P_NB(\bar{X}^N_{\lfloor t\rfloor/\tau})Q^{\frac12}e_i\big)\Big]\Big|\,dt,
\end{align*}
and
\begin{align*}
\mathcal R_{212}
:={}&\frac12\int_{0}^{T}\Big|\E\Big[\sum_{i=1}^\infty D^2u^N(T-t,\tilde{X}^N(t)).\big(E(t-\lfloor t\rfloor)P_N(B(\tilde{X}^N(t))-B(\bar{X}^N_{\lfloor t\rfloor/\tau}))Q^{\frac12}e_i,\\
	&\quad  E(t-\lfloor t\rfloor)P_N(B(\tilde{X}^N(t))-B(\bar{X}^N_{\lfloor t\rfloor/\tau}))Q^{\frac12}e_i\big)\Big]\Big|\,dt .
\end{align*}
 We split $\mathcal R_{212}
=
\mathcal R_{212}^{(0)}
+
\mathcal R_{212}^{(1)},$ corresponding to the intervals $[0,t_1]$ and $[t_1,T]$, respectively. On the first interval $[0,t_1]$, using
$\bar X_0^N=P_Nx_0$, H\"older's inequality, and the uniform
moment and exponential-moment bounds, we obtain
\begin{align*}
\mathcal R_{212}^{(0)}
&\le
C_{\epsilon,T}
\int_0^{t_1}
\left(1+(T-t)^{-\frac{1-\epsilon}{2}}\right)
\times
\left\|
(-A)^{-\frac{1-\epsilon}{4}}
(\tilde X^N(t)-P_Nx_0)
\right\|_{L^{2(1+\delta)}(\Omega;H)}^2dt
\le C_{\epsilon,T}\tau
\le C_{\epsilon,T}\tau^{1-\epsilon}.
\end{align*}
On $[t_1,T]$, using Proposition~\ref{cor-secondderivative} with $\gamma_1=\gamma_2=\frac{1-\epsilon}{4}$, Assumption~\ref{assum:nonlinearity-B} with $\varsigma_1=\frac{1-\epsilon}{4}$, and Lemma~\ref{lem:localerror} with $\delta=\frac{\epsilon}{1-\epsilon}$, yields 
\begin{align*}
&\mathcal R_{212}^{(1)}\le \frac12\int_{t_1}^{T}\Big|\E\Big[\sum_{i=1}^\infty D^2u^N(T-t,\tilde{X}^N(t)).\Big(E(t-\lfloor t\rfloor)P_N(B(\tilde{X}^N(t))-B(\bar{X}^N_{\lfloor t\rfloor/\tau}))Q^{\frac12}e_i,\\
	&\quad  E(t-\lfloor t\rfloor)P_N(B(\tilde{X}^N(t))-B(\bar{X}^N_{\lfloor t\rfloor/\tau}))Q^{\frac12}e_i\Big)\Big]\Big|\,dt\\
	&\le C_{L,\epsilon,T}\int_{t_1}^T(1\!+\!(T\!-\!t)^{-\frac{1-\epsilon}{2}})\times \E\Big[\exp\big(\frac{2\|(-A)^{\frac{1-\epsilon}{4}}\tilde{X}^N(t)\|}{e^{LT}}\big)\sum_{i=1}^\infty\big\|(-A)^{-\frac{1-\epsilon}{4}}(B(\tilde{X}^N(t))\!-\!B(\bar{X}^N_{\lfloor t\rfloor/\tau}))Q^{\frac12}e_i\big\|^2\Big]\,dt\\
	&\le  C_{L,\epsilon,T}\int_{t_1}^T(1\!+\!(T\!-\!t)^{-\frac{1-\epsilon}{2}})\times \E\Big[\exp\big(\frac{2\|(-A)^{\frac{1-\epsilon}{4}}\tilde{X}^N(t)\|}{e^{LT}}\big)\times\big\|(-A)^{-\frac{1-\epsilon}{4}}(\tilde{X}^N(t)\!-\!\bar{X}^N_{\lfloor t\rfloor/\tau})\big\|^2\Big]\,dt\\
	&\le C_{\epsilon,T}\tau^{1-\epsilon}\times\exp\big(2\|(-A)^{\frac{1-\epsilon}{4}}x_0\|\big)\le C_{\epsilon,T}\tau^{1-\epsilon}.
\end{align*}
Consequently, $\mathcal R_{212}
\le C_{\epsilon,T}\tau^{1-\epsilon}.$ 

For rigor, the following calculation is first performed with the $i$-sum truncated. Assumption~\ref{assum:nonlinearity-B} provides the required Hilbert--Schmidt summability, while the moment estimates above give integrability; hence Fubini's theorem and dominated convergence justify passing to the infinite sum. By It\^o's formula, we have
\begin{align*}
	&\mathcal R_{211}\le \int_{0}^{T}\Big|\E\Big[\sum_{i=1}^\infty D^2u^N(T-t,\tilde{X}^N(t)).\Big(E(t-\lfloor t\rfloor)P_N(B(\tilde{X}^N(t))-B(\bar{X}^N_{\lfloor t\rfloor/\tau}))Q^{\frac12}e_i,\\
	&\quad \quad\quad E(t-\lfloor t\rfloor)P_NB(\bar{X}^N_{\lfloor t\rfloor/\tau})Q^{\frac12}e_i\Big)\Big]\Big|\,dt\\
	&= \int_0^T\Big|\E\Big[\sum_{i=1}^\infty \sum_{k=1}^N\big\langle P_N(B(\tilde{X}^N(t))-B(\bar{X}^N_{\lfloor t\rfloor/\tau}))Q^{\frac12}e_i, e_k\big\rangle  \\
	&~~~\quad \quad\quad\times D^2u^N(T-t,\tilde{X}^N(t)).\big(E(t-\lfloor t\rfloor)e_k, E(t-\lfloor t\rfloor)P_NB(\bar{X}^N_{\lfloor t\rfloor/\tau})Q^{\frac12}e_i\big)\Big]\Big|\,dt\\
	&=\int_0^T\Big|\E\Big[\sum_{i=1}^\infty \sum_{k=1}^ND^2u^N(T-t,\tilde{X}^N(t)).\big(E(t-\lfloor t\rfloor)e_k,  E(t-\lfloor t\rfloor)P_NB(\bar{X}^N_{\lfloor t\rfloor/\tau})Q^{\frac12}e_i\big)\\
	&~~~\quad \times \bigg\{\int_{\lfloor t\rfloor}^tD\big(\big\langle P_NB(\cdot)Q^{\frac12}e_i, e_k\big\rangle\big)(\tilde{X}^N(s)).\big(A\tilde{X}^N(s)+E(s-\lfloor t\rfloor)F_N(\bar{X}^N_{\lfloor t\rfloor/\tau})\big)ds\\
	&~~~ +\int_{\lfloor t\rfloor}^tD\big(\big\langle P_NB(\cdot)Q^{\frac12}e_i, e_k\big\rangle\big)(\tilde{X}^N(s)).E(s-\lfloor t\rfloor)P_NB(\bar{X}^N_{\lfloor t\rfloor/\tau})dW^Q(s)\\
	&~~~+\int_{\lfloor t\rfloor}^t\int_{\chi}\Big(\big\langle P_NB\big(\tilde{X}^N(s)+E(s-\lfloor t\rfloor)G_N(\bar{X}^N_{\lfloor t\rfloor/\tau},z)\big)Q^{\frac12}e_i, e_k\big\rangle-\big\langle P_NB(\tilde{X}^N(s))Q^{\frac12}e_i, e_k\big\rangle\Big)\mathbf{\tilde{N}}(dz,ds)\\
	&~~~+\frac12\int_{\lfloor t\rfloor}^t\sum_{\ell=1}^\infty D^2\big(\big\langle P_NB(\cdot)Q^{\frac12}e_i, e_k\big\rangle\big)(\tilde{X}^N(s)).\big(E(s-\lfloor t\rfloor)P_NB(\bar{X}^N_{\lfloor t\rfloor/\tau})Q^{\frac12}e_{\ell}, E(s-\lfloor t\rfloor)P_NB(\bar{X}^N_{\lfloor t\rfloor/\tau})Q^{\frac12}e_{\ell}\big)\,ds\\
	&~~~+\int_{\lfloor t\rfloor}^t\int_{\chi}\Big(\big\langle P_NB\big(\tilde{X}^N(s)+E(s-\lfloor t\rfloor)G_N(\bar{X}^N_{\lfloor t\rfloor/\tau},z)\big)Q^{\frac12}e_i, e_k\big\rangle-\big\langle P_NB(\tilde{X}^N(s))Q^{\frac12}e_i, e_k\big\rangle\\
	&~~~\quad\quad \quad -D\big(\big\langle P_NB(\cdot)Q^{\frac12}e_i, e_k\big\rangle\big)(\tilde{X}^N(s)).E(s-\lfloor t\rfloor)G_N(\bar{X}^N_{\lfloor t\rfloor/\tau},z)\Big)\,\nu(dz)ds\bigg\}\Big]\Big|\,dt.
    \end{align*}
   Hence, it follows that
   \[
\mathcal R_{211}
\le
\mathcal R_{211,1}
+\mathcal R_{211,2}
+\mathcal R_{211,3}
+\mathcal R_{211,4}
+\mathcal R_{211,5},
\]
where
\begin{align*}
	\mathcal R_{211,1}&:=\int_0^T\Big|\E\Big[\sum_{i=1}^\infty \sum_{k=1}^N\int_{\lfloor t\rfloor}^tD^2u^N(T-t,\tilde{X}^N(t)).\big(E(t-\lfloor t\rfloor)e_k, E(t-\lfloor t\rfloor)P_NB(\bar{X}^N_{\lfloor t\rfloor/\tau})Q^{\frac12}e_i\big)\\
	&~\quad \times \big\langle P_N[DB(\tilde{X}^N(s)).\big(A\tilde{X}^N(s)+E(s-\lfloor t\rfloor)F_N(\bar{X}^N_{\lfloor t\rfloor/\tau})\big)]Q^{\frac12}e_i, e_k\big\rangle \,ds\Big]\Big|\,dt,\\
	\mathcal R_{211,2}&:=\frac12\int_0^T\Big|\E\Big[\sum_{i=1}^\infty \sum_{k=1}^N\int_{\lfloor t\rfloor}^tD^2u^N(T-t,\tilde{X}^N(t)).\big(E(t-\lfloor t\rfloor)e_k, E(t-\lfloor t\rfloor)P_NB(\bar{X}^N_{\lfloor t\rfloor/\tau})Q^{\frac12}e_i\big)\\
	&~\quad \times \sum_{\ell=1}^\infty \big\langle P_N[D^2B(\tilde{X}^N(s)).\big(E(s-\lfloor t\rfloor)P_NB(\bar{X}^N_{\lfloor t\rfloor/\tau})Q^{\frac12}e_{\ell}, E(s-\lfloor t\rfloor)P_NB(\bar{X}^N_{\lfloor t\rfloor/\tau})Q^{\frac12}e_{\ell}\big)]Q^{\frac12}e_i, e_k\big\rangle \,ds\Big]\Big|\,dt,\\
	\mathcal R_{211,3}&:=\int_0^T\Big|\E\Big[\sum_{i=1}^\infty \sum_{k=1}^N\int_{\lfloor t\rfloor}^t\int_{\chi}D^2u^N(T-t,\tilde{X}^N(t)).\big(E(t-\lfloor t\rfloor)e_k, E(t-\lfloor t\rfloor)P_NB(\bar{X}^N_{\lfloor t\rfloor/\tau})Q^{\frac12}e_i\big)\\
	&~\quad\times \Big(\big\langle P_NB\big(\tilde{X}^N(s)+E(s-\lfloor t\rfloor)G_N(\bar{X}^N_{\lfloor t\rfloor/\tau},z)\big)Q^{\frac12}e_i, e_k\big\rangle-\big\langle P_NB(\tilde{X}^N(s))Q^{\frac12}e_i, e_k\big\rangle\\
	&\quad\quad \quad -\big\langle P_N[DB(\tilde{X}^N(s)).E(s-\lfloor t\rfloor)G_N(\bar{X}^N_{\lfloor t\rfloor/\tau},z)]Q^{\frac12}e_i, e_k\big\rangle \Big)\,\nu(dz)ds \Big]\Big|\,dt,\\
	\mathcal R_{211,4}&:=\int_0^T\Big|\E\Big[\sum_{i=1}^\infty \sum_{k=1}^ND^2u^N(T-t,\tilde{X}^N(t)).\big(E(t-\lfloor t\rfloor)e_k,  E(t-\lfloor t\rfloor)P_NB(\bar{X}^N_{\lfloor t\rfloor/\tau})Q^{\frac12}e_i\big)\\
	&~\quad\times\sum_{n\in\N_+}\int_{\lfloor t\rfloor}^t\big\langle P_N[DB(\tilde{X}^N(s)).(E(s-\lfloor t\rfloor)P_NB(\bar{X}^N_{\lfloor t\rfloor/\tau})Q^{\frac12}e_n)]Q^{\frac12}e_i, e_k\big\rangle \,d\beta_n(s)\Big]\Big|\,dt,\\
	\mathcal R_{211,5}&:=\int_0^T\Big|\E\Big[\sum_{i=1}^\infty \sum_{k=1}^ND^2u^N(T-t,\tilde{X}^N(t)).\big(E(t-\lfloor t\rfloor)e_k,  E(t-\lfloor t\rfloor) P_NB(\bar{X}^N_{\lfloor t\rfloor/\tau})Q^{\frac12}e_i\big)\\
	&~\quad\times\int_{\lfloor t\rfloor}^t\int_{\chi}\int_0^1\big\langle P_N[DB\big(\tilde{X}^N(s)+\theta E(s-\lfloor t\rfloor)G_N(\bar{X}^N_{\lfloor t\rfloor/\tau},z)\big). E(s-\lfloor t\rfloor)G_N(\bar{X}^N_{\lfloor t\rfloor/\tau},z)]Q^{\frac12}e_i, e_k\big\rangle\,d\theta \mathbf{\tilde{N}}(dz,ds)\Big]\Big|\,dt.
\end{align*}
For the subterm $\mathcal R_{211,1}$, using  \eqref{ineq:normal-second} and  Proposition~\ref{cor-secondderivative} with $\gamma_1=\frac{1}{4}$, $\gamma_2=\frac{1-\epsilon}{4}$ and $p=4$ gives
\begin{align*}
	&\mathcal R_{211,1}\le \int_0^T \int_{\lfloor t\rfloor}^t\sum_{i=1}^\infty\E\Big[\Big|D^2 u^N(T-t,\tilde{X}^N(t)).\Big(E(t-\lfloor t\rfloor) P_NB(\bar{X}^N_{\lfloor t\rfloor/\tau})Q^{\frac12}e_i,\\
	&~~~\quad E(t-\lfloor t\rfloor) P_NDB(\tilde{X}^N(s)).\big(A\tilde{X}^N(s)+E(s-\lfloor t\rfloor)F_N(\bar{X}^N_{\lfloor t\rfloor/\tau})\big)Q^{\frac12}e_i\Big)\Big|\Big]\,dsdt\\
	&\le C_{\mathfrak{p},T}\int_0^{t_1}\int_0^t\E\Big[\big(1+\|\tilde{X}^N(t)\|^{\mathfrak{p}}\big)\sum_{i=1}^\infty\|B(P_Nx_0)Q^{\frac12}e_i\|\|DB(\tilde{X}^N(s)).\big(A\tilde{X}^N(s)+E(s)F_N(P_Nx_0)\big)Q^{\frac12}e_i\|\Big]\,dsdt\\
	&~~~+C_{L,\epsilon,T}\int_{t_1}^T \int_{\lfloor t\rfloor}^t(1+(T-t)^{-\frac{1-\epsilon/2}{2}})\times \E\Big[\exp\big(\frac{2\|(-A)^{\frac{1-\epsilon}{4}}\tilde{X}^N(t)\|}{e^{LT}}\big) \big\|(-A)^{-\frac{1-\epsilon}{4}}B(\bar{X}^N_{\lfloor t\rfloor/\tau})\big\|_{\IL_2^0(H)}\\
	&\quad\quad\times \big\|(-A)^{-\frac{1}{4}} DB(\tilde{X}^N(s)).\big(A\tilde{X}^N(s)+E(s-\lfloor t\rfloor)F_N(\bar{X}^N_{\lfloor t\rfloor/\tau})\big)\big\|_{\IL_2^0(H)}\Big]\,dsdt.
	\end{align*}
Then under Assumption~\ref{assum:nonlinearity-B} with $r=0$ and $\varsigma_1=\frac{1}{4}$, applying  the CFL condition $\tau\lambda_N\le \frac{\sqrt{\ln 2}}{4}$, \eqref{eq:estimate-X}, and Lemma~\ref{lem:localerror}	 we derive
\begin{align*}
\mathcal R_{211,1}&\le C_{\mathfrak{p},T}\int_0^{t_1}\int_0^t\big(1+\|\tilde{X}^N(t)\|_{L^{2\mathfrak{p}}(\Omega;H)}^{\mathfrak{p}}\big)\big\|A\tilde{X}^N(s)+E(s)F_N(P_Nx_0)\big\|_{L^{2}(\Omega;H)}\,dsdt\\
	&~~~ +C_{L,\epsilon,T}\int_{t_1}^T \int_{\lfloor t\rfloor}^t(1+(T-t)^{-\frac{1-\epsilon/2}{2}})\times \E\Big[\exp\big(\frac{2\|(-A)^{\frac{1-\epsilon}{4}}\tilde{X}^N(t)\|}{e^{LT}}\big)\\
	&\quad\quad \times(1+\|\bar{X}^N_{\lfloor t\rfloor/\tau}\|^{\gamma})\big(\big\|(-A)^{\frac{3}{4}}\tilde{X}^N(s)\big\|+1+\|\bar{X}^N_{\lfloor t\rfloor/\tau}\| \big)\Big]\,dsdt\\
	&\le C_{\mathfrak{p},T}(1+\|x_0\|)\int_0^{t_1}\int_0^t\big(\lambda_N\|\tilde{X}^N(s)\|_{L^{2}(\Omega;H)}+\|F(P_Nx_0)\|_{L^{2}(\Omega;H)}\big)\,dsdt\\
	&~~~ +C_{\epsilon,T}\exp(3\|(-A)^{\frac{1-\epsilon}{4}}x_0\|)\Bigl(\tau+\int_{t_1}^T (1+(T-t)^{-\frac{1-\epsilon/2}{2}})\int_{\lfloor t\rfloor}^t \big\|(-A)^{\frac{3}{4}}\tilde{X}^N(s)\big\|_{L^2(\Omega;H)}\,dsdt\Bigr)\\
	&\le C_{\epsilon,T}\Big(\tau+\tau^{1-\epsilon/4}\int_{t_1}^T (1+(T-t)^{-\frac{1-\epsilon/2}{2}})\big(1+(t-t_1)^{-\frac12}\big)\,dt\Big)\le C_{\epsilon,T}\tau^{1-\epsilon/4}.
\end{align*}
The estimates of subterms $\mathcal R_{211,2}$ and $\mathcal R_{211,3}$ are similar to that for $\mathcal R_{211,1}$ and are omitted here. For the subterm $\mathcal R_{211,4}$, for $s\in(\lfloor t\rfloor,t)$, $\bar X^N_{\lfloor t\rfloor/\tau}$ is $\mathcal F_{\lfloor t\rfloor}$-measurable, and hence, by \eqref{eq:zero-W}, $\mathcal D_s^{W,Q^{1/2}e_n}\bar X^N_{\lfloor t\rfloor/\tau}=0,~ ds\otimes d\mathbb P\text{-a.e.}$ Therefore, when applying the Wiener--Malliavin duality formula, the Malliavin derivative acts only on $D^2u^N(T-t,\tilde X^N(t))$ in the frozen factors. The analogous locality property $\mathcal D^N_{(s,z)}\bar X^N_{\lfloor t\rfloor/\tau}=0,~ s\in(\lfloor t\rfloor,t),$ will be used below for $\mathcal R_{211,5}$. Then by Proposition~\ref{prop:dualityformula}, we obtain
\begin{align*}
	&\mathcal R_{211,4}\le  \int_0^T\Big|\E\Big[\sum_{i=1}^\infty \sum_{k=1}^ND^2u^N(T-t,\tilde{X}^N(t)).\big(E(t-\lfloor t\rfloor)e_k,  E(t-\lfloor t\rfloor)P_NB(\bar{X}^N_{\lfloor t\rfloor/\tau})Q^{\frac12}e_i\big)\\
	&~~~\quad\times\sum_{n\in\N_+}\int_{\lfloor t\rfloor}^t\big\langle P_NDB(\tilde{X}^N(s)).(E(s-\lfloor t\rfloor)P_NB(\bar{X}^N_{\lfloor t\rfloor/\tau})Q^{\frac12}e_n)Q^{\frac12}e_i, e_k\big\rangle \,d\beta_n(s)\Big]\Big|\,dt\\
	&\le  \int_0^T\Big|\E\Big[\sum_{i=1}^\infty \sum_{k=1}^N\sum_{n\in\N_+}\int_{\lfloor t\rfloor}^tD^3u^N(T-t,\tilde{X}^N(t)).\big(\mathcal D^{W,Q^{\frac12}e_n}_s\tilde{X}^N(t), E(t-\lfloor t\rfloor)e_k, \\
	&~~~\quad E(t-\lfloor t\rfloor) P_NB(\bar{X}^N_{\lfloor t\rfloor/\tau})Q^{\frac12}e_i\big)\times\big\langle P_NDB(\tilde{X}^N(s)).(E(s-\lfloor t\rfloor)P_NB(\bar{X}^N_{\lfloor t\rfloor/\tau})Q^{\frac12}e_n)Q^{\frac12}e_i, e_k\big\rangle \,ds\Big]\Big|\,dt\\
	&\le \int_0^T\sum_{i=1}^\infty \sum_{n\in\N_+}\int_{\lfloor t\rfloor}^t\Big|\E\Big[D^3u^N(T-t,\tilde{X}^N(t)).\Big(\mathcal D^{W,Q^{\frac12}e_n}_s\tilde{X}^N(t),  E(t-\lfloor t\rfloor) P_NB(\bar{X}^N_{\lfloor t\rfloor/\tau})Q^{\frac12}e_i, \\
	&~~~\quad E(t-\lfloor t\rfloor)P_NDB(\tilde{X}^N(s)).(E(s-\lfloor t\rfloor)P_NB(\bar{X}^N_{\lfloor t\rfloor/\tau})Q^{\frac12}e_n)Q^{\frac12}e_i\Big) \Big]\Big|\,dsdt.
	\end{align*}
 Using Lemma \ref{lem:thirdderivative},  the moment estimate in \eqref{eq:estimate-X}, and Lemma \ref{lem:Malliavinderivative-B1}  gives
\begin{align*}
	&\mathcal R_{211,4}\le C_{L,\epsilon,T}\int_0^T(1+(T-t)^{-\frac{\varrho_1}{2}}) \sum_{n\in\N_+}\int_{\lfloor t\rfloor}^t\E\Big[\exp\big(\frac{\|(-A)^{\frac{1-\epsilon}{4}}\tilde{X}^N(t)\|}{4e^{LT}}\big)\times\big\|\mathcal D^{W,Q^{\frac12}e_n}_s\tilde{X}^N(t)\big\| \\
	&~~~\times\sum_{i=1}^\infty \big\|E(t-\lfloor t\rfloor) P_NB(\bar{X}^N_{\lfloor t\rfloor/\tau})Q^{\frac12}e_i\big\|\times \big\|DB(\tilde{X}^N(s)).(E(s-\lfloor t\rfloor)P_NB(\bar{X}^N_{\lfloor t\rfloor/\tau})Q^{\frac12}e_n)Q^{\frac12}e_i\big\|\Big]\,dsdt\\
	&\le C_{L,\epsilon,T}\int_0^T (1+(T-t)^{-\frac{\varrho_1}{2}})\int_{\lfloor t\rfloor}^t\E\Big[\exp\big(\frac{\|(-A)^{\frac{1-\epsilon}{4}}\tilde{X}^N(t)\|}{4e^{LT}}\big)\big(1+\|\bar{X}^N_{\lfloor t\rfloor/\tau}\|\big)\\
	&~~~\times \sum_{n\in\N_+}\big\|\mathcal D^{W,Q^{\frac12}e_n}_s\tilde{X}^N(t)\big\| \big\|B(\bar{X}^N_{\lfloor t\rfloor/\tau})Q^{\frac12}e_n\big\|\Big]\,dsdt\\
	&\le C_{L,\epsilon,T}\int_0^T (1+(T-t)^{-\frac{\varrho_1}{2}})\int_{\lfloor t\rfloor}^t\big|\exp\big(\frac{\|(-A)^{\frac{1-\epsilon}{4}}\tilde{X}^N(t)\|}{4e^{LT}}\big)\big(1+\|\bar{X}^N_{\lfloor t\rfloor/\tau}\|\big)^2\big|_{L^2(\Omega;\R)} \\
	&~~~\times\big(\sum_{n\in\N_+}\E\big[\big\|\mathcal D^{W,Q^{\frac12}e_n}_s\tilde{X}^N(t)\big\|^2\big]\big)^{\frac12}\,dsdt\le C_{\epsilon,T}\tau.
\end{align*} 
Similarly, using Proposition~\ref{prop:dualityformula}, we derive
\begin{align*}
	&\mathcal R_{211,5}\le \int_0^T\Big|\E\Big[\sum_{i=1}^\infty \sum_{k=1}^ND^2u^N(T-t,\tilde{X}^N(t)).\big(E(t-\lfloor t\rfloor)e_k,  E(t-\lfloor t\rfloor) P_NB(\bar{X}^N_{\lfloor t\rfloor/\tau})Q^{\frac12}e_i \big)\\
	&~\quad\times\int_{\lfloor t\rfloor}^t\int_{\chi}\int_0^1\big\langle P_NDB\big(\tilde{X}^N(s)+\theta E(s-\lfloor t\rfloor)G_N(\bar{X}^N_{\lfloor t\rfloor/\tau},z)\big). E(s-\lfloor t\rfloor)G_N(\bar{X}^N_{\lfloor t\rfloor/\tau},z)Q^{\frac12}e_i, e_k\big\rangle\,d\theta \mathbf{\tilde{N}}(dz,ds)\Big]\Big|\,dt\\
	&\le  \int_0^T\Big|\sum_{i=1}^\infty \sum_{k=1}^N\E\Big[\int_{\lfloor t\rfloor}^t\int_{\chi}\int_0^1\Big(D^2u^N(T-t,\tilde{X}^N(t)+\ID^{\mathbf{N}}_{(s,z)}\tilde{X}^N(t)).\big(E(t-\lfloor t\rfloor)e_k, E(t-\lfloor t\rfloor) P_NB(\bar{X}^N_{\lfloor t\rfloor/\tau})Q^{\frac12}e_i \big)\\
	&~~~-D^2 u^N(T-t,\tilde{X}^N(t)).\big(E(t-\lfloor t\rfloor)e_k,  E(t-\lfloor t\rfloor) P_NB(\bar{X}^N_{\lfloor t\rfloor/\tau})Q^{\frac12}e_i \big)\Big)\\
	&~~~~\times \big\langle P_NDB\big(\tilde{X}^N(s)+\theta E(s-\lfloor t\rfloor)G_N(\bar{X}^N_{\lfloor t\rfloor/\tau},z)\big). E(s-\lfloor t\rfloor)G_N(\bar{X}^N_{\lfloor t\rfloor/\tau},z)Q^{\frac12}e_i, e_k\big\rangle\,d\theta\nu(dz)ds\Big]\Big|\,dt\\
	&\le  \int_0^T\sum_{i=1}^\infty \Big|\E\Big[\int_{\lfloor t\rfloor}^t\int_{\chi}\int_0^1\Big(D^2 u^N(T-t,\tilde{X}^N(t)+\ID^{\mathbf{N}}_{(s,z)}\tilde{X}^N(t)).\big( E(t-\lfloor t\rfloor) P_NB(\bar{X}^N_{\lfloor t\rfloor/\tau})Q^{\frac12}e_i,\\
	&~~~E(t-\lfloor t\rfloor) P_NDB\big(\tilde{X}^N(s)+\theta E(s-\lfloor t\rfloor)G_N(\bar{X}^N_{\lfloor t\rfloor/\tau},z)\big). E(s-\lfloor t\rfloor)G_N(\bar{X}^N_{\lfloor t\rfloor/\tau},z)Q^{\frac12}e_i\big)\\
    &~~~-D^2 u^N(T-t,\tilde{X}^N(t)).\big(E(t-\lfloor t\rfloor) P_NB(\bar{X}^N_{\lfloor t\rfloor/\tau})Q^{\frac12}e_i,\\
	&~~~E(t-\lfloor t\rfloor) P_NDB\big(\tilde{X}^N(s)+\theta E(s-\lfloor t\rfloor)G_N(\bar{X}^N_{\lfloor t\rfloor/\tau},z)\big). E(s-\lfloor t\rfloor)G_N(\bar{X}^N_{\lfloor t\rfloor/\tau},z)Q^{\frac12}e_i\big)\Big)\,d\theta\nu(dz)ds\Big]\Big|\,dt.
	\end{align*}
Applying the result \eqref{ineq:normal-second} in Proposition~\ref{cor-secondderivative} and  Lemma~\ref{lem:Malliavinderivative-B1}, we obtain
\begin{align*}
	&\mathcal R_{211,5}\le C_{\mathfrak{p},T}\int_0^T\int_{\lfloor t\rfloor}^t\int_{\chi}\int_0^1\E\Big[\big(1+\|\tilde{X}^N(t)\|^{\mathfrak{p}}+\|\ID^{\mathbf{N}}_{(s,z)}\tilde{X}^N(t)\|^{\mathfrak{p}}\big) \times \sum_{i=1}^\infty\big\|B(\bar{X}^N_{\lfloor t\rfloor/\tau})Q^{\frac12}e_i\big\|\\
    &~~~\times\big\|DB\big(\tilde{X}^N(s)+\theta E(s-\lfloor t\rfloor)G_N(\bar{X}^N_{\lfloor t\rfloor/\tau},z)\big). E(s-\lfloor t\rfloor)G_N(\bar{X}^N_{\lfloor t\rfloor/\tau},z)Q^{\frac12}e_i\big\|\Big]\,d\theta\nu(dz)dsdt\\
    &\le C_{\mathfrak{p},T}\int_0^T\int_{\lfloor t\rfloor}^t\int_{\chi}\E\Big[\big(1+\|\tilde{X}^N(t)\|^{\mathfrak{p}}+\|\ID^{\mathbf{N}}_{(s,z)}\tilde{X}^N(t)\|^{\mathfrak{p}}\big) \times(1+\|\bar{X}^N_{\lfloor t\rfloor/\tau}\|)\|G(\bar{X}^N_{\lfloor t\rfloor/\tau},z)\|\Big]\,\nu(dz)dsdt \\
    &\le C_{\mathfrak{p},T}\int_0^T\int_{\lfloor t\rfloor}^t\big(1+\big(\E\big[\|\tilde{X}^N(t)\|^{2\mathfrak{p}}\big]\big)^{\frac12}\big)(1+(\E[\|\bar{X}^N_{\lfloor t\rfloor/\tau}\|^2])^{\frac12})\,dsdt\le C_{\epsilon,T}\tau.
\end{align*}
Hence combining the above results gives $\sum_{m=0}^{M-1}|\mathbb{H}_{m,2}|\le C_{\epsilon,T}\tau^{1-3\epsilon}$.

\textbf{Estimate of term $\mathbb{H}_{m,3}$.} We decompose this term into three subterms by using the mean value theorem: 
\begin{align*}
&\sum_{m=0}^{M-1}\mathbb{H}_{m,3}:=\int_{0}^{T}\!\!\int_\chi\Big\{\E\Big[u^N\big(T-t,\tilde{X}^N(t)+E(t-\lfloor t\rfloor)G_N(\bar{X}^N_{\lfloor t\rfloor/\tau},z)\big)\\
&~~~-u^N\big(T-t, \tilde{X}^N(t)+G_N(\tilde{X}^N(t),z)\big)\Big]+\E\Big[Du^N(T-t,\tilde{X}^N(t)).G_N(\tilde{X}^N(t),z)\\
	&~~~-Du^N(T-t,\tilde{X}^N(t)).(E(t-\lfloor t\rfloor)G_N(\bar{X}^N_{\lfloor t\rfloor/\tau},z))\Big]\Big\}\nu(dz)dt\\
	&=-\int_{0}^{T} \int_\chi \E\Big[Du^N\big(T-t,\tilde{X}^N(t)+\Gamma_{1}(t,z)\big).\big(G_N(\tilde{X}^N(t),z)-E(t-\lfloor t\rfloor)G_N(\bar{X}^N_{\lfloor t\rfloor/\tau},z)\big)\Big]\nu(dz)dt\\
	&~~~+\int_{0}^{T} \int_\chi  \E\Big[Du^N(T-t,\tilde{X}^N(t)).\big(G_N(\tilde{X}^N(t),z)-E(t-\lfloor t\rfloor)G_N(\bar{X}^N_{\lfloor t\rfloor/\tau},z)\big)\Big]\nu(dz)dt\\
	&=\int_{0}^{T} \int_\chi \E\Big[\Big(Du^N\big(T-t,\tilde{X}^N(t)+E(t-\lfloor t\rfloor)G_N(\bar{X}^N_{\lfloor t\rfloor/\tau},z)\big)-Du^N\big(T-t,\tilde{X}^N(t)+\Gamma_{1}(t,z)\big)\Big).\\
	&~~\quad \big(G_N(\tilde{X}^N(t),z)-E(t-\lfloor t\rfloor)G_N(\bar{X}^N_{\lfloor t\rfloor/\tau},z)\big)\Big]\nu(dz)dt\\
	&~-\int_{0}^{T} \int_\chi\!\!\E\Big[Du^N(T\!-\!t,\tilde{X}^N(t)\!+\!E(t-\lfloor t\rfloor)G_N(\bar{X}^N_{\lfloor t\rfloor/\tau},z)).\big(G_N(\tilde{X}^N(t),z)\!-\!E(t\!-\!\lfloor t\rfloor)G_N(\bar{X}^N_{\lfloor t\rfloor/\tau},z)\big)\Big]\nu(dz)dt\\
	&~+\int_{0}^{T} \int_\chi \E\Big[Du^N(T-t,\tilde{X}^N(t)).\big(G_N(\tilde{X}^N(t),z)-E(t-\lfloor t\rfloor)G_N(\bar{X}^N_{\lfloor t\rfloor/\tau},z)\big)\Big]\nu(dz)dt,
\end{align*}
where $\Gamma_{1}(t,z)\!:=\!\!E(t-\lfloor t\rfloor)G_N(\bar{X}^N_{\lfloor t\rfloor/\tau},z)+\theta_0\big(G_N(\tilde{X}^N(t),z)\!-\!E(t\!-\!\lfloor t\rfloor)G_N(\bar{X}^N_{\lfloor t\rfloor/\tau},z)\big)$ for some measurable $\theta_0=\theta_0(\omega,t,z)\in[0,1]$. 

Define $\mathbb H_3:=\sum_{m=0}^{M-1}\mathbb H_{m,3}$ and decompose $\mathbb H_3
=
\mathbb H_{3,1}
+\mathbb H_{3,2}
+\mathbb H_{3,3}$.

\textit{Estimate of subterm $\mathbb{H}_{3,1}$.} Using the mean value theorem again and applying Proposition \ref{cor-secondderivative} with $\gamma_1=\gamma_2=\frac{1-\epsilon}{4}$ yield 
\begin{align*}
	|\mathbb{H}_{3,1}|&\le \int_{0}^{T} \int_\chi \E\Big[\Big|\Big(Du^N\big(T-t,\tilde{X}^N(t)+E(t-\lfloor t\rfloor)G_N(\bar{X}^N_{\lfloor t\rfloor/\tau},z)\big)-Du^N\big(T-t,\tilde{X}^N(t)+\Gamma_{1}(t,z)\big)\Big).\\
	&~~~\quad \big(G_N(\tilde{X}^N(t),z)-E(t-\lfloor t\rfloor)G_N(\bar{X}^N_{\lfloor t\rfloor/\tau},z)\big)\Big|\Big]\nu(dz)dt\\
	&\le \int_{0}^{T} \int_\chi\E\Big[\Big|D^2u^N(T-t,\tilde{X}^N(t)+\Gamma_2(t,z)\big).\big(G_N(\tilde{X}^N(t),z)\!-\!E(t\!-\!\lfloor t\rfloor)G_N(\bar{X}^N_{\lfloor t\rfloor/\tau},z),\\
	&\quad ~~~ G_N(\tilde{X}^N(t),z)-E(t-\lfloor t\rfloor)G_N(\bar{X}^N_{\lfloor t\rfloor/\tau},z)\big)\Big|\Big]\nu(dz)dt\\
	&\le C_{L,\epsilon,T} \int_{0}^{T}(1+(T-t)^{-\frac{1-\epsilon}{2}}) \int_\chi\E\Big[\exp\big(\frac{2\|(-A)^{\frac{1-\epsilon}{4}}(\tilde{X}^N(t)+\Gamma_2(t,z))\|}{e^{LT}}\big)\\
	&~~~\times\big\|(-A)^{-\frac{1-\epsilon}{4}}\big(G_N(\tilde{X}^N(t),z)\!-\!E(t\!-\!\lfloor t\rfloor)G_N(\bar{X}^N_{\lfloor t\rfloor/\tau},z)\big)\big\|^2\Big]\nu(dz)dt\\
	&\le C_{L,\epsilon,T} \int_{0}^{T}(1+(T-t)^{-\frac{1-\epsilon}{2}}) \int_\chi\E\Big[\exp\big(\frac{2\|(-A)^{\frac{1-\epsilon}{4}}(\tilde{X}^N(t)+\Gamma_2(t,z))\|}{e^{LT}}\big)\\
	&~~~\times\Big(\big\|(-A)^{-\frac{1-\epsilon}{4}}\big(\mathrm I-E(t\!-\!\lfloor t\rfloor)\big)G_N(\bar{X}^N_{\lfloor t\rfloor/\tau},z)\big\|^2+\big\|(-A)^{-\frac{1-\epsilon}{4}}\big(G(\tilde{X}^N(t),z)\!-\!G(\bar{X}^N_{\lfloor t\rfloor/\tau},z)\big)\big\|^2\Big)\Big]\nu(dz)dt,
\end{align*}
where $\Gamma_{2}(t,z)\!:=\!\!E(t-\lfloor t\rfloor)G_N(\bar{X}^N_{\lfloor t\rfloor/\tau},z)+\theta_1\theta_0\big(G_N(\tilde{X}^N(t),z)\!-\!E(t\!-\!\lfloor t\rfloor)G_N(\bar{X}^N_{\lfloor t\rfloor/\tau},z)\big)$ for some measurable $\theta_1=\theta_1(\omega,t,z)\in[0,1]$. Similar to the estimate of subterm $\mathcal R_{212}$, using the property~\eqref{prop:semigroup2} of the semigroup, Lemma \ref{lem:localerror} with $\delta=\frac{\epsilon}{1-\epsilon}$, and Lemma \ref{lem:momentbounds-fulldiscre} gives
\begin{align*}
|\mathbb{H}_{3,1}|&\le C_{L,\epsilon,T} \int_{0}^{T}(1+(T-t)^{-\frac{1-\epsilon}{2}}) \int_\chi e^{2\lambda_1^{-\frac{\epsilon}{4}}\bar{c}(z)}\times \E\Big[\exp\big(\frac{2\|(-A)^{\frac{1-\epsilon}{4}}\tilde{X}^N(t)\|}{e^{LT}}\big)\|(-A)^{\frac14}G(\bar{X}^N_{\lfloor t\rfloor/\tau},z) \|^2\Big]\\
&\quad\times \|(-A)^{-\frac{1-\epsilon}{4}-\frac14} (\mathrm I-E(t\!-\!\lfloor t\rfloor))\|_{\IL(H)}^2\,\nu(dz)dt\\
&~+C_{L,\epsilon,T} \int_{0}^{T}(1+(T-t)^{-\frac{1-\epsilon}{2}}) \E\Big[\exp\big(\frac{2\|(-A)^{\frac{1-\epsilon}{4}}\tilde{X}^N(t)\|}{e^{LT}}\big)\times\int_\chi \int_0^1(e^{2\lambda_1^{-\frac{\epsilon}{4}}\bar{c}(z)}-1+1)\\
&\quad\times  \Big\|(-A)^{-\frac{1-\epsilon}{4}}DG(\bar{X}^N_{\lfloor t\rfloor/\tau}+\theta(\tilde{X}^N(t)-\bar{X}^N_{\lfloor t\rfloor/\tau}),z).(\tilde{X}^N(t)-\bar{X}^N_{\lfloor t\rfloor/\tau})\Big\|^2d\theta\nu(dz)\Big] dt\\
&\le C_{L,\epsilon,T}\tau^{1-\epsilon} \int_{0}^{T}(1+(T-t)^{-\frac{1-\epsilon}{2}})\E\Big[\exp\big(\frac{2\|(-A)^{\frac{1-\epsilon}{4}}\tilde{X}^N(t)\|}{e^{LT}}\big)\Big]\,dt\\
&\quad +C_{L,\epsilon,T}\int_{0}^{T}(1+(T-t)^{-\frac{1-\epsilon}{2}}) \E\Big[\exp\big(\frac{2\|(-A)^{\frac{1-\epsilon}{4}}\tilde{X}^N(t)\|}{e^{LT}}\big)\big\|(-A)^{-\frac{1-\epsilon}{4}}(\tilde{X}^N(t)-\bar{X}^N_{\lfloor t\rfloor/\tau})\big\|^2\Big]\,dt\\
&\le C_{\epsilon,T}\tau^{1-\epsilon}\exp(2\|(-A)^{\frac{1-\epsilon}{4}}x_0\|)\le C_{\epsilon,T}\tau^{1-\epsilon},
\end{align*}
where we used the H\"older inequality, the fact that $(e^x-1)^2\le x^2e^{2x}$ for $x\in\R_+$, and conditions in Assumption~\ref{assum:nonlinearity-G} (i)--(ii) in the second inequality.

The estimate of $\mathbb H_{3,3}$ follows from the same argument as that for $\mathbb H_{3,2}$, and is slightly simpler since no shifted spatial argument appears in $Du^N$. Next we focus on the estimate of subterm $\mathbb{H}_{3,2}$. 

\textit{Estimate of subterm $\mathbb{H}_{3,2}$.} We first decompose this subterm as follows,
\begin{align*}
\!\!|\mathbb{H}_{3,2}|&\le \!\!	\int_{0}^{T} \!\!\int_\chi\!\!\Big|\E\Big[Du^N(T\!-\!t,\tilde{X}^N(t)\!+\!E(t-\lfloor t\rfloor)G_N(\bar{X}^N_{\lfloor t\rfloor/\tau},z)).\big(G_N(\tilde{X}^N(t),z)\!-\!E(t\!-\!\lfloor t\rfloor)G_N(\bar{X}^N_{\lfloor t\rfloor/\tau},z)\big)\Big]\Big|\nu(dz)dt\\
	&\le \int_{0}^{T} \!\!\int_\chi\!\! \E\Big[\Big|Du^N(T\!-\!t,\tilde{X}^N(t)\!+\!E(t-\lfloor t\rfloor)G_N(\bar{X}^N_{\lfloor t\rfloor/\tau},z)).\big((\mathrm{I}-E(t\!-\!\lfloor t\rfloor)G_N(\bar{X}^N_{\lfloor t\rfloor/\tau},z)\big)\Big|\Big]\nu(dz)dt\\
	&~~+\int_{0}^{T} \!\!\int_\chi\!\! \Big|\E\Big[Du^N(T\!-\!t,\tilde{X}^N(t)\!+\!E(t-\lfloor t\rfloor)G_N(\bar{X}^N_{\lfloor t\rfloor/\tau},z)).\big(G_N(\tilde{X}^N(t),z)\!-\!G_N(\bar{X}^N_{\lfloor t\rfloor/\tau},z)\big)\Big]\Big|\nu(dz)dt\\
	&=: \mathbb{H}_{3,2,1}+\mathbb{H}_{3,2,2}.
\end{align*}
For the subterm $\mathbb{H}_{3,2,1}$, using Proposition~\ref{cor:firstderivative2} with $\gamma_0=\frac{3(1-\epsilon)}{4}$ when $\epsilon_1\in(0,\frac{\epsilon}{2(2-\epsilon)})$ gives
\begin{align*}
	\mathbb{H}_{3,2,1}&\le C_{L,\epsilon,T}\int_0^T\int_{\chi}(1+(T-t)^{-\frac{3(1-\epsilon)}{4}})\times \E\Big[\exp\big(\frac{\|(-A)^{\frac{1-\epsilon}{4}}(\tilde{X}^N(t)\!+\!E(t-\lfloor t\rfloor)G_N(\bar{X}^N_{\lfloor t\rfloor/\tau},z))\|}{2e^{LT}}\big)\\
	&~~~\times \big\|(-A)^{-\frac{3(1-\epsilon)}{4}}(\mathrm{I}-E(t\!-\!\lfloor t\rfloor)(-A)^{-\frac{1}{4}}\big\|_{\IL(H)}\big\|(-A)^{\frac{1}{4}}G(\bar{X}^N_{\lfloor t\rfloor/\tau},z)\big\|\Big]\nu(dz)dt\\
 &\le C_{L,\epsilon,T}\tau^{1-\epsilon}\int_{\chi}e^{\lambda_1^{-\frac{\epsilon}{4}}\bar{c}(z)}\times \bar{c}(z)\,\nu(dz)\sup_{t\in[0,T]}\E\Big[\exp\big(\frac{\|(-A)^{\frac{1-\epsilon}{4}}\tilde{X}^N(t)\|}{2e^{LT}}\big)\Big]\\
 &\le C_{\epsilon,T}\tau^{1-\epsilon}\exp(\|(-A)^{\frac{1-\epsilon}{4}}x_0\|)\le C_{\epsilon,T}\tau^{1-\epsilon}.
\end{align*}
For the subterm $\mathbb{H}_{3,2,2}$, define $
Y_t^z
:=
\tilde X^N(t)
+
E(t-\lfloor t\rfloor)G_N(\bar X^N_{\lfloor t\rfloor/\tau},z).$ For the Wiener martingale term, using the Wiener--Malliavin
duality formula and the fact that $D_s^{W,Q^{1/2}e_n}\bar X^N_{\lfloor t\rfloor/\tau}=0$, $s\in(\lfloor t\rfloor,t)$, we obtain
\begin{align*}
&\Big|\E\Big[D u^N(T-t,Y_t^z).e_k
\int_{\lfloor t\rfloor}^t
D\langle P_NG(\cdot,z),e_k\rangle(\tilde X^N(s))
E(s-\lfloor t\rfloor)P_NB(\bar X^N_{\lfloor t\rfloor/\tau})\,dW^Q(s)\Big]\Big|\\
&\le \sum_{n\in\mathbb N_+}
\int_{\lfloor t\rfloor}^t\Big|\E\Big[D^2u^N(T-t,Y_t^z)
\big(D_s^{W,Q^{1/2}e_n}\tilde X^N(t),e_k
\big)
\\
&\quad\quad\quad\times
\big\langle
P_NDG(\tilde X^N(s),z).E(s-\lfloor t\rfloor)P_NB(\bar X^N_{\lfloor t\rfloor/\tau})Q^{1/2}e_n,e_k\big\rangle\Big]\Big|\,ds .
\end{align*}
For the Poisson martingale term, using the Poisson--Malliavin duality formula and $D^N_{(s,v)}\bar X^N_{\lfloor t\rfloor/\tau}=0,$ $s\in(\lfloor t\rfloor,t),$ we obtain
\[
\begin{aligned}
&\E\Big[
D u^N(T-t,Y_t^z).e_k
\int_{\lfloor t\rfloor}^t\int_\chi
\Big\langle
P_N\big(
G(\tilde X^N(s)+g_s(v),z)
-G(\tilde X^N(s),z)
\big),
e_k
\Big\rangle\tilde N(dv,ds)\Big]
\\&=\E\int_{\lfloor t\rfloor}^t\int_\chi
\Big[D u^N\big(
T-t,Y_t^z+D^N_{(s,v)}\tilde X^N(t)
\big)-D u^N(T-t,Y_t^z)
\Big].e_k
\\
&\qquad\qquad\times
\Big\langle
P_N\big(
G(\tilde X^N(s)+g_s(v),z)
-G(\tilde X^N(s),z)
\big),
e_k
\Big\rangle
\,\nu(dv)\,ds,
\end{aligned}
\]
where $g_s(v):=
E(s-\lfloor t\rfloor)G_N(\bar X^N_{\lfloor t\rfloor/\tau},v).$
The drift, Gaussian correction, and Poisson compensator terms
are estimated analogously to the corresponding terms in
$\mathcal R_{211}$. Combining these estimates with the two
duality estimates above gives
\[
\mathbb H_{3,2,2}
\le C_{\epsilon,T}\tau^{1-3\epsilon}.
\]
Hence combining the above results gives $\big|\sum_{m=0}^{M-1}\mathbb{H}_{m,3}\big|\le C_{\epsilon,T}\tau^{1-3\epsilon}$. 
  
  Combining the estimates of
$\mathbb H_{m,1}$, $\mathbb H_{m,2}$, and $\mathbb H_{m,3}$,
we obtain
\[
\left|
\E[\varphi(X^N(T))]
-
\E[\varphi(\bar X_M^N)]
\right|
\le
C_{\epsilon,T}\tau^{1-3\epsilon}.
\]
Together with the spatial weak error estimate obtained in
Subsection~5.1 and the triangle inequality,
\[
\begin{aligned}
\left|
\E[\varphi(X(T))]
-
\E[\varphi(\bar X_M^N)]
\right|
&\le
\left|
\E[\varphi(X(T))]
-
\E[\varphi(X^N(T))]
\right|+
\left|
\E[\varphi(X^N(T))]
-
\E[\varphi(\bar X_M^N)]
\right|
\\
&\le
C_{\epsilon,T}
\left(
\lambda_N^{-1+3\epsilon}
+
\tau^{1-3\epsilon}
\right).
\end{aligned}
\]
This proves Theorem~\ref{the:weakconvergence-full}.

\subsubsection{Proof of Lemma \ref{lem:localerror} }\label{subsec5.2.2}

	From expression~\eqref{eq:continuoustime} of process $\tilde{X}^N$, and under Assumptions \ref{assum:nonlinearity-F}--\ref{assum:nonlinearity-G}, applying \eqref{prop:semigroup1} and \eqref{eq:estimate-X} yields that for all $t\in(0,T]$ and $x_0\in \dot{H}^{\frac{1-\epsilon/2}{2}}$
\begin{align}\label{eq:3-epsilon/4}
	&\big\|(-A)^{\frac{3-\epsilon}{4}}\tilde{X}^N(t)\big\|_{L^2(\Omega;H)}\le \|(-A)^{\frac{3-\epsilon}{4}}E(t)P_Nx_0\|+\int_0^t\big\|(-A)^{\frac{3-\epsilon}{4}}E(t-\lfloor s\rfloor)\big\|_{\IL(H)}\big\|F_N(\bar{X}^N_{\lfloor s\rfloor/\tau})\big\|_{L^2(\Omega;H)}ds \nonumber\\
		&\quad+\Big(\int_0^t\E\Big[\big\|(-A)^{\frac{1-\epsilon/4}{2}}E(t-\lfloor s\rfloor)\big\|^2_{\IL(H)}\big\|(-A)^{\frac{1-\epsilon/2}{4}}B(\bar{X}^N_{\lfloor s\rfloor/
		\tau})\big\|_{\IL_2^0(H)}^2\Big]ds\Big)^{\frac12}\nonumber\\
		&\quad+\Big(\int_0^t\int_\chi\E\Big[\big\|(-A)^{\frac{1-\epsilon/2}{2}}E(t-\lfloor s\rfloor)\big\|_{\IL(H)}^2\big\|(-A)^{\frac{1}{4}}G_N(\bar{X}^N_{\lfloor s\rfloor/
		\tau},z)\big\|^2\Big]\nu(dz)ds\Big)^{\frac12}\nonumber\\
		&\le Ct^{-\frac{1}{2}} \|(-A)^{\frac{1-\epsilon}{4}}x_0\|+C\Big[\int_0^t(t-\lfloor s\rfloor)^{-\frac{3-\epsilon}{4}}\,ds+\Big(\int_0^t(t-\lfloor s\rfloor)^{-1+\epsilon/4}\,ds\Big)^{\frac12}+\Big(\int_0^t(t-\lfloor s\rfloor)^{-1+\epsilon/2}\,ds\Big)^{\frac12}\Big]\nonumber\\
        &\quad\times\Big[1\!+\!\sup_{s\in[0,t]}\|\bar{X}^N_{\lfloor s\rfloor/
		\tau}\|_{L^2(\Omega;H)}+\Big(\sup_{s\in[0,t]}\|(-A)^{\frac{1-\epsilon/2}{4}}\bar{X}^N_{\lfloor s\rfloor/
		\tau}\|_{L^2(\Omega;H)}^2+\int_\chi|\bar{c}(z)|^2\nu(dz)\Big)^{\frac12}\Big]\nonumber\\
		&\le C_{\epsilon,T}(1+t^{-\frac{1}{2}} ).
	\end{align}
In addition,  using the It\^o formula gives that for all $s\in(0,T]$
\begin{align*}
&d\|	(-A)^{\frac{1}{4}}\tilde{X}^N(s)\|^2=\Big\{-2\|(-A)^{\frac{3}{4}}\tilde{X}^N(s)\|^2+2\big\langle (-A)^{\frac{1}{2}}\tilde{X}^N(s), E(s-\lfloor s \rfloor) F_N(\bar{X}^N_{\lfloor s \rfloor/\tau}) \big\rangle\\
&\quad\quad +\bigl\|(-A)^{\frac{1}{4}} E(s - \lfloor s \rfloor) P_N B(\bar{X}^N_{\lfloor s \rfloor/\tau})\bigr\|_{\IL^0_2(H)}^2+\int_{\chi}\bigl\|(-A)^{\frac{1}{4}} E(s - \lfloor s \rfloor) G_N(\bar{X}^N_{\lfloor s \rfloor/\tau}, z)\bigr\|^2 \nu(dz)\Bigr\}\,ds\\
&\quad\quad +2\bigl\langle (-A)^{\frac{1}{2}}\tilde{X}^N(s), E(s-\lfloor s \rfloor) P_NB(\bar{X}^N_{\lfloor s \rfloor/\tau})\,dW^Q(s)\bigr\rangle \\
&\quad\quad + \int_{\chi}\Bigl( 2\bigl\langle (-A)^{\frac{1}{2}}\tilde{X}^N(s), E(s-\lfloor s \rfloor) G_N(\bar{X}^N_{\lfloor s \rfloor/\tau}, z)\bigr \rangle +\Vert{}(-A)^{\frac{1}{4}} E(s - \lfloor s \rfloor) G_N(\bar{X}^N_{\lfloor s \rfloor/\tau}, z) \Vert{}^2 \Bigr) \,\tilde{\mathbf{N}}(dz, ds). 
\end{align*}
Integrating  both sides from $\lfloor t \rfloor$ to $t$ for $t\in[t_1,T]$ and taking expectations show that
\begin{align*}
	2\int_{\lfloor t\rfloor}^t &\E\bigl[\big\|(-A)^{\frac{3}{4}}\tilde{X}^N(s)\big\|^2\bigr]\,ds=\E\bigl[\big\|(-A)^{\frac{1}{4}}\tilde{X}^N(\lfloor t \rfloor)\big\|^2-\big\|(-A)^{\frac{1}{4}}\tilde{X}^N(t)\big\|^2\bigr]\\
	&+2\int_{\lfloor t \rfloor}^t\E\Big[\big\langle (-A)^{\frac{1-\epsilon}{4}}\tilde{X}^N(s), (-A)^{\frac{1+\epsilon}{4}}E(s-\lfloor s \rfloor) F_N(\bar{X}^N_{\lfloor s \rfloor/\tau}) \big\rangle\Big]\,ds\\
	&+\int_{\lfloor t \rfloor}^t\E\Big[\bigl\|(-A)^{\frac{1}{4}} E(s - \lfloor s \rfloor) P_N B(\bar{X}^N_{\lfloor s \rfloor/\tau})\bigr\|_{\IL^0_2(H)}^2\Big]\,ds+\int_{\lfloor t \rfloor}^t\int_{\chi}\E\Bigl[\bigl\|(-A)^{\frac{1}{4}} E(s - \lfloor s \rfloor) G_N(\bar{X}^N_{\lfloor s \rfloor/\tau}, z)\bigr\|^2\Bigr] \nu(dz)ds.
\end{align*}
Hence, under Assumptions \ref{assum:nonlinearity-F}--\ref{assum:nonlinearity-G}, by using 
the estimate \eqref{eq:3-epsilon/4}, we obtain 
\begin{align*}
\int_{\lfloor t\rfloor}^t \E\bigl[&\big\|(-A)^{\frac{3}{4}}\tilde{X}^N(s)\big\|^2\bigr]\,ds\le \E\Big[\big\langle (-A)^{\frac{1}{4}}(\tilde{X}^N(\lfloor t \rfloor)-\tilde{X}^N(t)),  (-A)^{\frac{1}{4}}\tilde{X}^N(\lfloor t \rfloor)\big\rangle \Big]\\
	&\quad + C\int_{\lfloor t \rfloor}^t(s-\lfloor s \rfloor)^{-\frac{\epsilon}{2}}\|(-A)^{\frac{1-\epsilon}{4}}\tilde{X}^N(s)\|_{L^2(\Omega;H)}\|(-A)^{\frac{1-\epsilon}{4}}F(\bar{X}^N_{\lfloor s \rfloor/\tau})\|_{L^2(\Omega;H)}\,ds\\
	&\quad +C_{\epsilon}\int_{\lfloor t \rfloor}^t(s-\lfloor s \rfloor)^{-\frac{\epsilon}{2}}\big(1+\|(-A)^{\frac{1-\epsilon}{4}}\bar{X}^N_{\lfloor s \rfloor/\tau}\|^2_{L^2(\Omega;H)}\big)\,ds +C\int_{\lfloor t \rfloor}^t\int_{\chi}|\bar{c}(z)|^2\nu(dz)ds\\
	&\le \E\Big[\big\langle (-A)^{\frac{1}{4}}(\mathrm I-E(t-\lfloor t\rfloor))\tilde{X}^N(\lfloor t \rfloor)-\int_{\lfloor t \rfloor}^t(-A)^{\frac{1}{4}}E(t-\lfloor t\rfloor)F_N(\tilde{X}^N(\lfloor t \rfloor))ds, (-A)^{\frac{1}{4}}\tilde{X}^N(\lfloor t \rfloor)\big\rangle\Big]\\
	&\quad +C_{\epsilon,T}\tau^{1-\epsilon/2}+C\int_{\lfloor t \rfloor}^t(s-\lfloor s \rfloor)^{-\frac{\epsilon}{2}}\|(-A)^{\frac{1-\epsilon}{4}}\tilde{X}^N(s)\|_{L^2(\Omega;H)}\|(-A)^{\frac{1-\epsilon}{4}}F(\bar{X}^N_{\lfloor s \rfloor/\tau})\|_{L^2(\Omega;H)}\,ds\\
	&\le C_{\epsilon,T}\tau^{1-\epsilon/2}+\|(-A)^{-1+\epsilon/2}(\mathrm I-E(t-\lfloor t\rfloor))\|_{\IL(H)}\|(-A)^{\frac{3-\epsilon}{4}}\tilde{X}^N(\lfloor t \rfloor)\|_{L^2(\Omega;H)}^2\\
	&\quad +C_\epsilon\int_{\lfloor t \rfloor}^t(s - \lfloor t \rfloor)^{-\frac{\epsilon}{2}}\|(-A)^{\frac{1-\epsilon}{4}}\tilde{X}^N(s)\|_{L^2(\Omega;H)}\big(1+\|(-A)^{\frac{1-\epsilon/2}{4}}\bar{X}^N_{\lfloor s \rfloor/\tau}\|_{L^2(\Omega;H)}\big)\,ds\\
	&\le C_{\epsilon,T}\tau^{1-\epsilon/2}(1+\lfloor t \rfloor^{-1}),
\end{align*}
where we used the condition in Assumption \ref{assum:nonlinearity-F} (i) with $r_1=\frac{1-\epsilon}{4}$ and $\varrho=\frac{\epsilon}{8}$ in the third inequality. Then
\[
\begin{aligned}
\int_{\lfloor t \rfloor}^t
\|(-A)^{\frac34}\tilde X^N(s)\|_{L^2(\Omega;H)}
\,ds &\le (t-\lfloor t \rfloor)^{1/2}
\Big( \int_{\lfloor t \rfloor}^t \E\|(-A)^{\frac34}\tilde X^N(s)\|^2ds \Big)^{1/2}\\
&\le
C_{\epsilon,T}
\tau^{1/2}
\tau^{\frac12-\frac{\epsilon}{4}}
(1+\lfloor t \rfloor^{-1})^{1/2}
\le C_{\epsilon,T} \tau^{1-\epsilon/4}
\bigl(1+\lfloor t \rfloor^{-1/2}\bigr).
\end{aligned}
\]
Finally, it holds that for any $\delta\in(0,\frac12)$ and all $t\in(t_m,t_{m+1}]$, 
\begin{align*}
	&\big\|(-A)^{-\frac{1-\epsilon}{4}}\big(\tilde{X}^N(t)-\bar{X}^N_m\big)\big\|_{L^{2(1+\delta)}(\Omega;H)}\le\big\|(E(t-t_m)-\mathrm I)(-A)^{-\frac{1-\epsilon}{2}}(-A)^{\frac{1-\epsilon}{4}}\bar{X}^N_m\big\|_{L^{2(1+\delta)}(\Omega;H)}\\
	&~~~+\int_{t_m}^t\big\|E(t-t_m)F_N(\bar{X}^N_m)\big\|_{L^{2(1+\delta)}(\Omega;H)}\,ds+C_T\Big(\int_{t_m}^t\!\E\Big[\big\|E(t-t_m)P_NB(\bar{X}^N_m)\big\|_{\IL_2^0(H)}^{2(1+\delta)}\Big]\,ds\Big)^{\frac{1}{2(1+\delta)}}\\
&~~~+C_T\Big(\int_{t_m}^t\int_{\chi}\!\E\Big[\big\|G(\bar{X}^N_{m},z)\big\|^{2(1+\delta)}\Big]\nu(dz)ds\Big)^{\frac{1}{2(1+\delta)}}\\
&\le C_T\tau^{\frac{1-\epsilon}{2}}\|(-A)^{\frac{1-\epsilon}{4}}\bar{X}^N_m\|_{L^{2(1+\delta)}(\Omega;H)}+C\tau^{\frac{1}{2(1+\delta)}}\big(1+\|\bar{X}^N_m\|_{L^{2(1+\delta)}(\Omega;H)}\big)\le C_{\epsilon,T}\tau^{\frac{1}{2(1+\delta)}\wedge \frac{1-\epsilon}{2}},
\end{align*}
where we used the inequality \eqref{prop:semigroup2} and Corollary \ref{cor:exponentialintegrability} in the last two steps. The proof is completed.

\subsubsection{Proof of Lemma  \ref{lem:thirdderivative}}\label{subsec5.2.3}

For each fixed $N\in\mathbb N_+$, the projected coefficients
$F_N$, $P_NB$, and $G_N$ are three times continuously
differentiable on the finite-dimensional space $H_N$. Applying \cite[Theorem 3.1]{BX-2012} successively to equation~\eqref{eq:spectralgalerkin} and to its first- and second-order variational equations yields
that the solution map $x\mapsto X^{N,(0,x)}(t)$ is three times differentiable. We denote $\zeta^N_{y,h,k,x}(t):=D^3X^{N,(0,x)}(t).(y,h,k).$ It is the unique mild solution of the following equation:
\begin{align*}
  	&	d\zeta^N_{y,h,k,x}(t)=\big(A\zeta^N_{y,h,k,x}(t)+P_NDF(X^N(t,x)).\zeta^N_{y,h,k,x}(t)+P_ND^2F(X^N(t,x)).(\xi^N_{y,h,x}(t),\eta^N_{k,x}(t))\big)dt\nonumber\\
  		&+P_ND^2F(X^N(t,x)).(\xi^N_{y,k,x}(t),\eta^N_{h,x}(t))dt+P_ND^2F(X^N(t,x)).(\xi^N_{h,k,x}(t),\eta^N_{y,x}(t))dt\nonumber\\
  		&+P_ND^3F(X^N(t,x)).(\eta^N_{y,x}(t),\eta^N_{h,x}(t),\eta^N_{k,x}(t))dt+P_N\big(DB(X^N(t,x)).\zeta^N_{y,h,k,x}(t)\nonumber\\
  		&+D^2B(X^N(t,x)).(\xi^N_{y,h,x}(t),\eta^N_{k,x}(t))+D^2B(X^N(t,x)).(\xi^N_{y,k,x}(t),\eta^N_{h,x}(t))+D^2B(X^N(t,x)).(\xi^N_{h,k,x}(t),\eta^N_{y,x}(t))\nonumber\\
  		&+D^3B(X^N(t,x)).(\eta^N_{y,x}(t),\eta^N_{h,x}(t),\eta^N_{k,x}(t))\big)dW^Q(t)+\int_\chi P_NDG(X^N(t,x),z).\zeta^N_{y,h,k,x}(t)\,\tilde{\mathbf{N}}(dz,dt)\nonumber\\
  		&+\int_\chi P_N\big(D^2G(X^N(t,x),z).(\xi^N_{y,h,x}(t),\eta^N_{k,x}(t))+D^2G(X^N(t,x),z).(\xi^N_{y,k,x}(t),\eta^N_{h,x}(t))\nonumber\\
  		&+D^2G(X^N(t,x),z).(\xi^N_{h,k,x}(t),\eta^N_{y,x}(t))+D^3G(X^N(t,x),z).(\eta^N_{y,x}(t),\eta^N_{h,x}(t),\eta^N_{k,x}(t))\big)\,\tilde{\mathbf{N}}(dz,dt)
  \end{align*}
with $
\zeta^N_{y,h,k,x}(0)=0.$ Define a  stopping time
\[
\begin{aligned}
\tau_R:=\inf\Big\{t\in[0,T]:\;&
\|X^N(t)\|
+\sum_{\ell\in\{y,h,k\}}\|\eta_{\ell,x}^N(t)\|
\\
&+\|\xi_{y,h,x}^N(t)\|
+\|\xi_{y,k,x}^N(t)\|
+\|\xi_{h,k,x}^N(t)\|
+\|\zeta_{y,h,k,x}^N(t)\|
\ge R
\Big\}\wedge T .
\end{aligned}
\]
We apply It\^o's formula to
$\|\zeta^N_{y,h,k,x}(t\wedge\tau_R)\|^2$. The stochastic Wiener and
compensated Poisson integrals have zero expectation. Moreover, $2\langle  \zeta^N_{y,h,k,x}(s),A\zeta^N_{y,h,k,x}(s)\rangle
=-2\|(-A)^{1/2}\zeta^N_{y,h,k,x}(s)\|^2.$ Using Assumptions~\ref{assum:nonlinearity-F}--\ref{assum:nonlinearity-G} and Young's inequality for the remaining drift, Gaussian quadratic-variation, and Poisson compensator terms, we obtain
\begin{align*}
&\E[\|\zeta_{y,h,k,x}^N(t\wedge\tau_R)\|^2]
+
\E\Big[\int_0^{t\wedge\tau_R}
\|(-A)^{1/2}\zeta_{y,h,k,x}^N(s)\|^2\,ds\Big]
\\
&\le
C\int_0^t
\E[\|\zeta_{y,h,k,x}^N(s\wedge\tau_R)\|^2]\,ds
+C_{\varsigma_0,\lambda_1}
\int_0^t\E\Big[
\|\eta_{y,x}^N(s)\|^2
\|\eta_{h,x}^N(s)\|^2
\|(-A)^{\varrho_0/2}\eta_{k,x}^N(s)\|^2
\Big]\,ds
\\
&\quad
+C_{\varsigma_0,\lambda_1}\int_0^t
\E\Big[\|\xi_{y,h,x}^N(s)\|^2\|\eta_{k,x}^N(s)\|^2
+
\|\xi_{y,k,x}^N(s)\|^2\|\eta_{h,x}^N(s)\|^2 +
\|\xi_{h,k,x}^N(s)\|^2\|\eta_{y,x}^N(s)\|^2
\Big]\,ds .
\end{align*}
where $\varrho_0\in(\frac12,1)$ is given by Assumption~\ref{assum:nonlinearity-F}. Since $\varrho_0\in(1/2,1)$, we fix $q\in\big(2,\frac{2}{\varrho_0}\big),$ so that $\varrho_0q<2$. 

Similar to the estimates of the moment bounds in \eqref{estimate:eta-kappa} and \eqref{estimate:xi-4-}, we derive that, for all $p\ge 2$ and $\tilde{\delta}\in(0,\frac12)$,
\begin{align}
\sup_{t\in[0,T]}
\|\eta_{k,x}^N(t)\|_{L^p(\Omega;H)}
&\le C_{p,T}\|k\|,\label{ineq:eta-N}\\
\sup_{t\in[0,T]}
\|\xi_{h_1,h_2,x}^N(t)\|_
{L^{\frac{4}{1+\tilde{\delta}}}(\Omega;H)}
&\le C_{\tilde{\delta},T}\|h_1\|\|h_2\|.\label{ineq:xi-N}
\end{align}
Then by the H\"older and interpolation inequalities and the fact that $\int_0^t\E\big[\|(-A)^{\frac{1}{2}}\eta^N_{y,x}(s)\|^2\big]\,ds\le C_T\|y\|^2$ for all $y\in H_N$ under Assumptions~\ref{assum:nonlinearity-F}--\ref{assum:nonlinearity-G} (cf. \cite[Lemma 4.1]{CDHZ-2023}), it holds that  
\begin{align*}
	&~~~\int_0^t\E\big[\|\eta^N_{y,x}(s)\|^2\|\eta^N_{h,x}(s)\|^2\|(-A)^{\frac{\varrho_0}{2}}\eta^N_{k,x}(s)\|^2\big]\,ds\\
	&\le \int_0^t\big|\|\eta^N_{y,x}(s)\|\|\eta^N_{h,x}(s)\|\big|_{L^{\frac{2q}{q-2}}(\Omega;\R)}^2\|(-A)^{\frac{\varrho_0}{2}}\eta^N_{k,x}(s)\|_{L^q(\Omega;H)}^2\,ds\\
	&\le  C_{T}\|y\|^2\|h\|^2\int_0^t\Big(\E\Big[\|\eta^N_{k,x}(s)\|^{q(1-\varrho_0)}\|(-A)^{\frac{1}{2}}\eta^N_{k,x}(s)\|^{q\varrho_0}\Big]\Big)^{\frac{2}{q}}\,ds \\
	&\le C_{T}\|y\|^2\|h\|^2\int_0^t\Big(\E\Big[\|\eta^N_{k,x}(s)\|^{\frac{2q(1-\varrho_0)}{2-q\varrho_0}}\Big]\Big)^{\frac{2-q\varrho_0}{q}}\Big(\E\Big[\|(-A)^{\frac{1}{2}}\eta^N_{k,x}(s)\|^2\Big]\Big)^{\varrho_0}\,ds\\
	&\le C_{T}\|y\|^2\|h\|^2\|k\|^{2(1-\varrho_0)}\Big(\int_0^t\E\big[\|(-A)^{\frac{1}{2}}\eta^N_{k,x}(s)\|^2\big]\,ds\Big)^{\varrho_0}\le C_{T}\|y\|^2\|h\|^2\|k\|^2.
\end{align*}
Combining \eqref{ineq:eta-N} and \eqref{ineq:xi-N} with the above estimates, we obtain
\[
\E\|\zeta_{y,h,k,x}^N(t\wedge\tau_R)\|^2
\le
C\int_0^t
\E\|\zeta_{y,h,k,x}^N(s\wedge\tau_R)\|^2\,ds
+
C_T\|y\|^2\|h\|^2\|k\|^2 ,
\]
where the constants are independent of $R$.
Hence, Gr\"onwall's inequality yields $\sup_{t\in[0,T]}
\E\|\zeta_{y,h,k,x}^N(t\wedge\tau_R)\|^2
\le
C_T\|y\|^2\|h\|^2\|k\|^2 .$ Since, for fixed $N$, all the variational processes involved have c\`adl\`ag paths on $[0,T]$, we have $\tau_R\uparrow T~\text{a.s. as }R\to\infty .$ Therefore, by Fatou's lemma, 
\begin{equation}\label{ineq:xi-t}
\sup_{t\in[0,T]}
\E\|\zeta_{y,h,k,x}^N(t)\|^2
\le
C_T\|y\|^2\|h\|^2\|k\|^2 .
\tag{5.11}
\end{equation}

Finally, since $\varrho_1\in(\frac12,1)$ is given by \eqref{setting-varrho1}, we fix an auxiliary exponent $q_1\in(2,\frac{2}{\varrho_1}),$ so that $q_1\varrho_1<2.$ Using the interpolation inequality $\|(-A)^{\varrho_1/2}\eta_{k,x}^N(t)\|
\le
\|\eta_{k,x}^N(t)\|^{1-\varrho_1}
\|(-A)^{1/2}\eta_{k,x}^N(t)\|^{\varrho_1},$ together with H\"older's inequality, estimates \eqref{ineq:eta-N}--\eqref{ineq:xi-t}, and the same argument as in the proof of \eqref{eq:regularity-eta}, we obtain
\[
\begin{aligned}
|D^3u^N(t,x).(y,h,k)|
&=
\big|
D^3\E[\varphi(X^{N,(0,x)}(t))].(y,h,k)
\big|
\\
&\le
C_p\E\Big[
(1+\|X^N(t)\|^p)
\big(
\|\zeta_{y,h,k,x}^N(t)\|
+\|\xi_{y,h,x}^N(t)\|\|\eta_{k,x}^N(t)\|
\\
&\qquad
+\|\xi_{y,k,x}^N(t)\|\|\eta_{h,x}^N(t)\|
+\|\xi_{h,k,x}^N(t)\|\|\eta_{y,x}^N(t)\|
+\|\eta_{y,x}^N(t)\|
 \|\eta_{h,x}^N(t)\|
 \|(-A)^{\varrho_1/2}\eta_{k,x}^N(t)\|
\big)
\Big]
\\
&\le
C_{L,\epsilon,T}
\exp\Big(
\frac{\|(-A)^{\frac{1-\epsilon}{4}}x\|}
     {4e^{LT}}
\Big)
\left(1+t^{-\frac{\varrho_1}{2}}\right)
\|y\|\,\|h\|\,\|k\|.
\end{aligned}
\]
Here $q_1$ is chosen independently of the exponent $q$ associated with $\varrho_0$ above. The proof is completed.

\section{Appendix}\label{sec:appendix}

\subsection{Proof of Lemma~\ref{lem:momentbounds-fulldiscre}}

Recall that $\tilde{X}^{N }(t_m)=\bar{X}^{N }_{m}$ for all $m\in\{0,1,\ldots,M-1\}$ and for all  $t\in[t_m,t_{m+1}]$,
\begin{align*}
	\!d\tilde{X}^{N }(t)\!=\!A\tilde{X}^{N }(t)dt+\!E(t\!-\!t_m)F_N(\bar{X}^{N }_{m})dt\!+\!E(t\!-\!t_m)P_NB(\bar{X}^{N }_{m})dW^Q(t)\!+\!E(t\!-\!t_m)\int_\chi \!G_N(\bar{X}^{N }_{m},z)\tilde{\mathbf{N}}(dz,dt).
\end{align*}
Set $\mu(\tilde{X}^{N }(t)):=A\tilde{X}^{N }(t)+E(t-t_m)F_N(\bar{X}_m^{N })$, $\sigma(\tilde{X}^{N }(t)):=E(t-t_{m})P_NB(\bar{X}^{N }_{m})Q^{\frac12}$, and $\mathbb g(\tilde{X}^{N }(t),z):=E(t-t_{m})G_N(\bar{X}^{N }_{m},z)$. We claim that there exists a constant $\mathcal{K}>1$ independent of $r_0$ such that
\begin{align}\label{ineq:estimate-62}
&~~~DU_{r_0}(\tilde{X}^{N }(t)).\mu(\tilde{X}^{N }(t))+\frac12\mbox{Tr}\big(D^2U_{r_0}(\tilde{X}^{N }(t)).\big(\sigma(\tilde{X}^{N }(t)),\sigma(\tilde{X}^{N }(t))\big)\big)\nonumber\\
&~~~+\frac{\|(\sigma(\tilde{X}^{N }(t)))^{*}DU_{r_0}(\tilde{X}^{N }(t))\|^2}{2}-\int_\chi DU_{r_0}(\tilde{X}^{N }(t)). \mathbb g(\tilde{X}^{N }(t),z) \nu(dz)\nonumber\\
&\le \frac{r_0}{2}\big\| (-A)^{\frac{1-\epsilon}{4}}\tilde{X}^{N }(t)\big\|+\frac{r_0(L_F+L_0)^2(1+\|\bar{X}^{N }_{m}\|^2)}{2\sqrt{1+\|(-A)^{\frac{1-\epsilon}{4}}\tilde{X}^{N }(t)\|^2}}+\mathcal{K}(1+r_0)U_{r_0}(\bar{X}^N_m).	
\end{align}
The proof of \eqref{ineq:estimate-62} is postponed to Subsection~\ref{Subsection:6.1.1}.  Now for  arbitrary positive constants $\varsigma_4$ and $\vartheta_2$, by It\^o's formula, integrating on both sides and taking the conditional expectation,  we obtain that, for all $m\in\{0,1,\ldots,M-1\}$, $\mathbb P$-a.s.,
\begin{align*}
&\Theta_{m,m+1}:=\E\Big[\exp\Big(\frac{U_{r_0}( \tilde{X}^{N }(t_{m+1}))}{e^{\varsigma_4\tau}}-\int_{t_m}^{t_{m+1}}\frac{\vartheta_2}{e^{\varsigma_4(s-t_m)}}ds\Big)\Big|\mathcal F_{t_m}\Big]-\exp\big(U_{r_0}( \tilde{X}^{N }(t_m))\big)\\
	&\le \E\Big[\int_{t_m}^{t_{m+1}}\!\!\!\!\exp\Big(\frac{U_{r_0}( \tilde{X}^{N }(s))}{e^{\varsigma_4(s-t_m)}}-\int_{t_m}^{s}\!\!\frac{\vartheta_2}{e^{\varsigma_4(r-t_m)}}dr\Big) e^{-\varsigma_4(s-t_m)}\times \Big(\!-\varsigma_4U_{r_0}( \tilde{X}^{N }(s))\!-\!\vartheta_2+\!DU_{r_0}( \tilde{X}^{N }(s)).\mu(\tilde{X}^{N }(s))\\
	&+\frac12\mbox{Tr}\big(D^2U_{r_0}(\tilde{X}^{N }(s)).\big(\sigma(\tilde{X}^{N }(s)),\sigma(\tilde{X}^{N }(s))\big)\big)+\frac12\|(\sigma(\tilde{X}^{N }(s)))^* DU_{r_0}(\tilde{X}^{N }(s))\|^2\Big)ds\\
	&+\int_{t_m}^{t_{m+1}}\exp\Big(\frac{U_{r_0}( \tilde{X}^{N }(s))}{e^{\varsigma_4(s-t_m)}}-\int_{t_m}^{s}\frac{\vartheta_2}{e^{\varsigma_4(r-t_m)}}dr\Big) e^{-\varsigma_4(s-t_m)}\times DU_{r_0}(\tilde{X}^{N }(s)).E(s-t_{m})P_NB(\bar{X}^{N }_{m})dW^Q(s)\\
  &+\!\!\int_{t_m}^{t_{m+1}}\!\!\!\int_\chi \!\!\Big[\exp\Big(\frac{U_{r_0}\big(\tilde{X}^{N }(s)\!+\!\mathbb g(\tilde{X}^{N }(s),z)\big)}{e^{\varsigma_4(s-t_m)}}\!-\!\int_{t_m}^{s}\!\!\!\frac{\vartheta_2}{e^{\varsigma_4(r-t_m)}}dr\Big)\!-\!\exp\Big(\frac{U_{r_0}(\tilde{X}^{N }(s))}{e^{\varsigma_4(s-t_m)}}\!-\!\int_{t_m}^{s}\!\!\frac{\vartheta_2}{e^{\varsigma_4(r-t_m)}}dr\Big) \!\Big]\tilde{\mathbf{N}}(dz,ds)\\
	&+\!\int_{t_m}^{t_{m+1}}\!\!\int_\chi\int_0^1 \Big\{\exp\Big(\frac{U_{r_0}\big(\tilde{X}^{N }(s)+\theta\mathbb g(\tilde{X}^{N }(s),z)\big)}{e^{\varsigma_4(s-t_m)}}-\int_{t_m}^{s}\frac{\vartheta_2}{e^{\varsigma_4(r-t_m)}}dr\Big)DU_{r_0}\big(\tilde{X}^{N }(s)\!+\!\theta\mathbb g(\tilde{X}^{N }(s),z)\big)\\
	&~~-\exp\Big(\frac{U_{r_0}(\tilde{X}^{N }(s))}{e^{\varsigma_4(s-t_m)}}-\int_{t_m}^{s}\frac{\vartheta_2}{e^{\varsigma_4(r-t_m)}}dr\Big) DU_{r_0}\big(\tilde{X}^{N }(s)\big)\Big\}.\mathbb g(\tilde{X}^{N }(s),z)\times e^{-\varsigma_4(s-t_m)}\,d\theta\nu(dz)ds~\Big| \mathcal F_{t_m}\Big].
\end{align*}
Then substituting \eqref{ineq:estimate-62} into the above inequality yields 
\begin{align}\label{estimate-Ito}
		\Theta_{m,m+1}&\le \E\bigg[\int_{t_m}^{t_{m+1}}\exp\Big(\frac{U_{r_0}(\tilde{X}^{N }(s))}{e^{\varsigma_4(s-t_m)}}-\int_{t_m}^{s}\frac{\vartheta_2}{e^{\varsigma_4(r-t_m)}}dr\Big) e^{-\varsigma_4(s-t_m)}\times \Big(-\varsigma_4U_{r_0}(\tilde{X}^N(s))-\vartheta_2\nonumber\\
		&~~~\quad +\frac{r_0}{2}\big\| (-A)^{\frac{1-\epsilon}{4}}\tilde{X}^{N }(s)\big\|+\frac{r_0(L_F+L_0)^2(1+\|\bar{X}^{N }_{m}\|^2)}{2\sqrt{1+\|(-A)^{\frac{1-\epsilon}{4}}\tilde{X}^{N }(s)\|^2}}+\mathcal{K}(1+r_0)U_{r_0}(\bar{X}^N_m)\Big)\,ds~\Big| \mathcal F_{t_m}\bigg]\nonumber\\
	&~~~+\E\bigg[\int_{t_m}^{t_{m+1}}\!\int_\chi\int_0^1\exp\Big(\frac{U_{r_0}\big(\tilde{X}^{N }(s)+\theta\mathbb g(\tilde{X}^{N }(s),z)\big)}{e^{\varsigma_4(s-t_m)}}-\int_{t_m}^{s}\frac{\vartheta_2}{e^{\varsigma_4(r-t_m)}}dr\Big)  e^{-\varsigma_4(s-t_m)}\nonumber\\
	&~~~\quad \times DU_{r_0}\big(\tilde{X}^{N }(s)+\theta\mathbb g(\tilde{X}^{N }(s),z)\big).\mathbb g(\tilde{X}^{N }(s),z)\,d\theta\nu(dz)ds~\bigg| \mathcal F_{t_m}\bigg].
\end{align}
Let $\mathcal J_m(s)$ denote the sum of the two $s$-integrands on the right-hand side of \eqref{estimate-Ito}. Hence, by the linearity of conditional expectation,
\[
\Theta_{m,m+1}
\le \E\Big[
\int_{t_m}^{t_{m+1}}\mathcal J_m(s)\,ds
\,\Big|\,\mathcal F_{t_m}
\Big].
\]

For later use, define the $\mathcal F_{t_m}$-measurable process $\{\mathrm{I}_{1,m}(s)\}_{s\in[t_m,t_{m+1}]}$ by
\[
\mathrm{I}_{1,m}(s)
:=
E(s-t_m)(-A)^{\frac{1-\epsilon}{4}}\bar X_m^N,
\qquad
s\in[t_m,t_{m+1}].
\]
Under the step-size restriction $\tau\in(0,\tau_N)$ with $\tau_N$ given in \eqref{def:tau_0}, we have
\begin{align}\label{ineq:compare-s,t_m}
	\|\mathrm I_{1,m}(s)\|=\big\|E(s-t_m)(-A)^{\frac{1-\epsilon}{4}}\bar{X}_m^{N}\big\|\ge e^{-\lambda_N\tau}\|(-A)^{\frac{1-\epsilon}{4}}\bar{X}^{N}_{m}\|\ge e^{-\frac{\sqrt{\ln 2}}{4}}\|(-A)^{\frac{1-\epsilon}{4}}\bar{X}^{N}_{m}\|.
\end{align} 

For convenience, set $\kappa_\epsilon:=\lambda_1^{-\epsilon/4}.$ Then, by Assumption~\ref{assum:nonlinearity-G} (i), $\|(-A)^{\frac{1-\epsilon}{4}}G(x,z)\|
\le
\kappa_\epsilon\bar c(z).$ For a fixed constant $k_1\ge 2$, define 
\begin{align}
\vartheta_2&:=k_1r_0+\int_\chi e^{r_0\kappa_\epsilon\bar c(z)}\,r_0\kappa_\epsilon\bar c(z)\,\nu(dz),\label{def:vartheta2}\\
 \varsigma_4&:=\frac12+\mathcal{K}(1+r_0)\times 2e^{\frac{\sqrt{\ln 2}}{4}}+2(L_F+L_0)^2(1+\lambda_1^{-\frac{1-\epsilon}{2}}e^{\frac{\sqrt{\ln 2}}{2}}),\label{def:varsigma4}
\end{align}
where  $\mathcal{K}$ is the constant appearing in \eqref{ineq:estimate-62}.

For $s\in[t_m,t_{m+1}]$, introduce the events
\[
\mathcal A_m(s)
:=
\Big\{
\|(-A)^{\frac{1-\epsilon}{4}}\tilde X^N(s)\|
\ge
\frac12\|I_{1,m}(s)\|
\Big\},
\qquad
\mathcal B_m(s):=\mathcal A_m(s)^c .
\]
By the definition of $U_{r_0}$ and the choices
\eqref{def:vartheta2}--\eqref{def:varsigma4}, the following pointwise estimates hold.

First, for $(ds\otimes d\mathbb P)$-a.e.
$(s,\omega)\in[t_m,t_{m+1}]\times\Omega$,
\begin{align}\label{ineq:estimate-64}
\mathcal J_m(s)\mathbf 1_{\mathcal A_m(s)}
    \le 0.
	\end{align}

Second, there exists a constant $\mathcal K_1>0$, independent of
$r_0$, such that, with
\[
\begin{aligned}
\widehat{\mathcal J}_m(s)
:={}&
\exp\Big(
\frac{
r_0\sqrt{1+\bigl(\frac12\|\mathrm I_{1,m}(s)\|\bigr)^2}
}{
e^{\varsigma_4(s-t_m)}
}
-
\int_{t_m}^{s}
\frac{\vartheta_2}{e^{\varsigma_4(r-t_m)}}\,dr
\Big) \times
e^{-\varsigma_4(s-t_m)}
\mathcal K_1r_0(1+r_0)
\bigl(1+\|\mathrm I_{1,m}(s)\|^2\bigr),
\end{aligned}
\]
we have
\begin{align}\label{ineq:estimate-65}
\mathcal J_m(s)\mathbf 1_{\mathcal B_m(s)}
\le
\widehat{\mathcal J}_m(s)
\mathbf 1_{\mathcal B_m(s)},
\qquad
ds\otimes d\mathbb P\text{-a.e.}
\end{align}
The proofs of \eqref{ineq:estimate-64} and \eqref{ineq:estimate-65} are postponed to Subsection~\ref{Subsection:6.1.2}.  

Combining \eqref{ineq:estimate-64} and \eqref{ineq:estimate-65}, we obtain
\[
\Theta_{m,m+1}
\le
\E\Big[ \int_{t_m}^{t_{m+1}}
\widehat{\mathcal J}_m(s)
\mathbf1_{\mathcal B_m(s)}\,ds
\,\Big|\,\mathcal F_{t_m}\Big].
\]
Choosing $C_L\ge7$, we further decompose
\begin{align}\label{ineq:estimate-Theta}
\Theta_{m,m+1}&\le
\E\Big[ \int_{t_m}^{t_{m+1}}
\widehat{\mathcal J}_m(s)
\mathbf1_{
\mathcal B_m(s)
\cap\{\|I_{1,m}(s)\|>C_L\}
}
\,ds
\,\Big|\,\mathcal F_{t_m}\Big]+
\int_{t_m}^{t_{m+1}}
\widehat{\mathcal J}_m(s)
\mathbf1_{\{\|I_{1,m}(s)\|\le C_L\}}
\,ds\nonumber\\
&=:\Theta_{m,m+1}^{1}+\Theta_{m,m+1}^{2},
\end{align}
where $\mathrm{I}_{1,m}(s)$ is $\mathcal F_{t_m}$-measurable for every fixed $s\in[t_m,t_{m+1}]$, so is $\widehat{\mathcal J}_m(s)$, and $\widehat{\mathcal J}_m(s)\ge0$.  Notice that $\mathcal B_m(s)$ depends on $\tilde X^N(s)$ and is in general not $\mathcal F_{t_m}$-measurable, which is why the conditional expectation is retained in the first term. 

For the term $\Theta_{m,m+1}^2$, there exists a constant $C>0$ dependent on $r_0$ such that 
\begin{align*}
\Theta_{m,m+1}^2&\le	 C\int_{t_m}^{t_{m+1}}\!\!\exp\Big(\frac{r_0\sqrt{1+(\frac12\|\mathrm I_{1,m}(s)\|)^2}}{e^{\varsigma_4(s-t_m)}}-\int_{t_m}^{s}\frac{\vartheta_2}{e^{\varsigma_4(r-t_m)}}dr\Big)\,ds.
\end{align*}
For the term $\Theta_{m,m+1}^1$, we claim that
 \begin{align}\label{estimate-Theta-1}
	\Theta_{m,m+1}^1&\le \int_{t_m}^{t_{m+1}} \exp\Big(\frac{r_0\sqrt{1+\big(\frac12\|(-A)^{\frac{1-\epsilon}{4}}\bar{X}^{N}_{m}\|\big)^2}}{e^{\varsigma_4(s-t_m)}}-\int_{t_m}^{s}\frac{\vartheta_2}{e^{\varsigma_4(r-t_m)}}dr\Big)\times e^{-\varsigma_4(s-t_m)}\nonumber\\
	&\quad\quad\times Cr_0(1+r_0)(s-t_m)\big(1+\|(-A)^{\frac{1-\epsilon}{4}}\bar{X}^{N}_{m}\|^{2\gamma}+\|(-A)^{\frac{1-\epsilon}{4}}\bar{X}^{N}_{m}\|^2\big)\,ds. 
\end{align}
The proof of inequality \eqref{estimate-Theta-1} is postponed to Subsection~\ref{Subsection:6.1.3}. Hence, for $\tau<\frac{1}{4}$, there exists a constant $C>0$ dependent on $r_0$ such that
\begin{align*}
	\Theta_{m,m+1}&
	\le \Theta_{m,m+1}^1+\Theta_{m,m+1}^2\\
	&\le C \int_{t_m}^{t_{m+1}}\!\!\exp\Big(\frac{r_0\sqrt{1+\frac14\|(-A)^{\frac{1-\epsilon}{4}}\bar{X}^{N}_{m}\|^2}}{e^{\varsigma_4(s-t_m)}}\!-\!\int_{t_m}^{s}\!\!\frac{\vartheta_2}{e^{\varsigma_4(r-t_m)}}dr\Big)ds\times\big(1\!+\! \tau r_0^2\|(-A)^{\frac{1-\epsilon}{4}}\bar{X}^{N}_{m}\|^2\big)\\
	&\le C\tau\times\exp\Big(\frac{r_0}{2}\|(-A)^{\frac{1-\epsilon}{4}}\bar{X}^{N}_{m}\|+\frac{r_0}{2}\|(-A)^{\frac{1-\epsilon}{4}}\bar{X}^{N}_{m}\|\Big)\le C\tau\times \exp\big(r_0\|(-A)^{\frac{1-\epsilon}{4}}\bar{X}^{N}_{m}\|\big)\\
	&\le C\tau\times\exp\big(r_0\sqrt{1+\|(-A)^{\frac{1-\epsilon}{4}}\tilde{X}^{N}(t_m)\|^2}\big)=C\tau\times\exp\big(U_{r_0}(\tilde{X}^{N}(t_m))\big),\end{align*}
which implies that 
\begin{align*}
& \E\Big[\exp\big(\frac{U_{r_0}( \tilde{X}^{N}(t_{m+1}))}{e^{\varsigma_4\tau}}\big)\Big|\mathcal F_{t_m}\Big]= \E\Big[\exp\Big(\frac{U_{r_0}( \tilde{X}^{N}(t_{m+1}))}{e^{\varsigma_4\tau}}-\int_{t_m}^{t_{m+1}}\frac{\vartheta_2}{e^{\varsigma_4(s-t_m)}}ds\Big)\Big|\mathcal F_{t_m}\Big]\times e^{\int_{t_m}^{t_{m+1}}\frac{\vartheta_2}{e^{\varsigma_4(s-t_m)}}ds}\\
	&\le (1+C\tau)e^{\vartheta_2\tau}\times \exp\big(U_{r_0}(\tilde{X}^{N}(t_m))\big). 
\end{align*}
Considering $e^{-\varsigma_4t_m}U_{r_0}$ instead of $U_{r_0}$ gives  $\E\Big[\exp\big(\frac{U_{r_0}( \tilde{X}^{N}(t_{m+1}))}{e^{\varsigma_4t_{m+1}}}\big)\Big]\le (1+C\tau)e^{\vartheta_2\tau}\E\Big[\exp\big(\frac{U_{r_0}(\tilde{X}^{N}(t_m))}{e^{\varsigma_4t_{m}}}\big)\Big]$. Then repeating this procedure yields that for all $m\in\{0,1,\ldots,M-1\}$,
\begin{align*}
 \E\Big[\exp\big(\frac{U_{r_0}( \tilde{X}^{N}(t_{m+1}))}{e^{\varsigma_4t_{m+1}}}\big)\Big]\le (1+C\tau)^{m+1}e^{\vartheta_2(m+1)\tau}\times \E\Big[\exp\big(U_{r_0}(\tilde{X}^{N}(0))\big)\Big]\le e^{CT}e^{U_{r_0}(x_0)}.
\end{align*}
For any \(t\in(t_m,t_{m+1}]\), the same argument, with \(t_{m+1}\) replaced by \(t\), yields
\[
\mathbb E\big[\exp\!\Big(
\frac{U_{r_0}(\tilde X^N(t))}{e^{\varsigma_4 t}}
\Big)\Big]\le
(1+C\tau)e^{\vartheta_2(t-t_m)}
\mathbb E\Big[\exp\!\Big(
\frac{U_{r_0}(\tilde X^N(t_m))}{e^{\varsigma_4 t_m}}
\Big)\Big].
\]
Since \(t-t_m\le\tau\), the preceding discrete-time estimates therefore imply
\[
\sup_{t\in[0,T]}
\mathbb E\Big[
\exp\Big(
\frac{U_{r_0}(\tilde X^N(t))}{e^{\varsigma_4 t}}
\Big)
\Big]
\le
e^{CT} e^{U_{r_0}(x_0)}.
\]
This completes the proof.

\subsubsection{Proof of the inequality \eqref{ineq:estimate-62}.}\label{Subsection:6.1.1}

 From the definition of mapping $U_{r_0}$, it follows that
\begin{align*}
&~~~DU_{r_0}(\tilde{X}^{N }(t)).\mu(\tilde{X}^{N }(t))+\frac12\mbox{Tr}\big(D^2U_{r_0}(\tilde{X}^{N }(t)).\big(\sigma(\tilde{X}^{N }(t)),\sigma(\tilde{X}^{N }(t))\big)\big)\\
	&~~~+\frac{\|(\sigma(\tilde{X}^{N }(t)))^{*}DU_{r_0}(\tilde{X}^{N }(t))\|^2}{2}-\int_\chi DU_{r_0}(\tilde{X}^{N }(t)). \mathbb g(\tilde{X}^{N }(t),z) \nu(dz) \\
	&=-\frac{r_0\big\langle (-A)^{\frac{1-\epsilon}{4}}\tilde{X}^{N }(t),(-A)^{1+\frac{1-\epsilon}{4}}\tilde{X}^{N }(t)\big\rangle}{\sqrt{1+\|(-A)^{\frac{1-\epsilon}{4}}\tilde{X}^{N }(t)\|^2}}+\frac{r_0\big\langle  (-A)^{\frac{1-\epsilon}{4}}\tilde{X}^{N }(t), E(t-t_m)(-A)^{\frac{1-\epsilon}{4}}F_N(\bar{X}^{N }_m)\big\rangle}{\sqrt{1+\|(-A)^{\frac{1-\epsilon}{4}}\tilde{X}^{N }(t)\|^2}}\\
	&~-\int_\chi \!\!\frac{r_0 \big\langle  (-A)^{\frac{1-\epsilon}{4}}\tilde{X}^{N }(t), E(t\!-\!t_m)(-A)^{\frac{1-\epsilon}{4}}G_N(\bar{X}^{N }_m,z)\big\rangle}{\sqrt{1+\big\|(-A)^{\frac{1-\epsilon}{4}}\tilde{X}^{N }(t)\big\|^2}}\nu(dz)+\frac{r_0\big\|E(t-t_{m})(-A)^{\frac{1-\epsilon}{4}}P_NB(\bar{X}^{N }_{m})\big\|_{\IL_2^0(H)}^2 }{2\sqrt{1+\|(-A)^{\frac{1-\epsilon}{4}}\tilde{X}^{N }(t)\|^2}}\\
	&~-\frac{r_0\big\|Q^{\frac12}(P_NB(\bar{X}^{N }_{m}))^{*}(-A)^{\frac{1-\epsilon}{4}}E(t-t_{m})(-A)^{\frac{1-\epsilon}{4}}\tilde{X}^{N }(t)\big\|^2}{2(1+\|(-A)^{\frac{1-\epsilon}{4}}\tilde{X}^{N }(t)\|^2)^{3/2}}\\
	&~+\frac{r_0^2\big\|Q^{\frac12}(P_NB(\bar{X}^{N }_{m}))^{*}(-A)^{\frac{1-\epsilon}{4}}E(t-t_{m})(-A)^{\frac{1-\epsilon}{4}}\tilde{X}^{N }(t)\big\|^2}{2(1+\|(-A)^{\frac{1-\epsilon}{4}}\tilde{X}^{N}(t)\|^2)}.
\end{align*}
Under Assumptions~\ref{assum:nonlinearity-F}--\ref{assum:nonlinearity-G} with $\gamma\in[0, \frac12)$, using the continuous embedding $\dot{H}^{\alpha_2}\subset \dot{H}^{\alpha_1}$ for $\alpha_1<\alpha_2$, which satisfies  $\|x\|_{\dot{H}^{\alpha_1}}\le \lambda_1^{-\frac{\alpha_2-\alpha_1}{2}}\|x\|_{\dot{H}^{\alpha_2}}$ for all $x\in \dot{H}^{\alpha_2}$, together with the  interpolation inequality $$\|x\|_{\dot{H}^{\alpha}}\le \|x\|^{\frac{\alpha_2-\alpha}{\alpha_2-\alpha_1}}_{\dot{H}^{\alpha_1}}\|x\|^{\frac{\alpha-\alpha_1}{\alpha_2-\alpha_1}}_{\dot{H}^{\alpha_2}},\quad \forall\, x\in\dot{H}^{\alpha_2},~\alpha\in[\alpha_1,\alpha_2],$$
 and the Young inequality, we derive that there exists a constant $\mathcal{K}>1$ independent of $r_0$ such that
\begin{align*}
&~~~DU_{r_0}(\tilde{X}^{N }(t)).\mu(\tilde{X}^{N }(t))+\frac12\mbox{Tr}\big(D^2U_{r_0}(\tilde{X}^{N }(t)).\big(\sigma(\tilde{X}^{N }(t)),\sigma(\tilde{X}^{N }(t))\big)\big)\nonumber\\
&~~~+\frac{\|(\sigma(\tilde{X}^{N }(t)))^{*}DU_{r_0}(\tilde{X}^{N }(t))\|^2}{2}-\int_\chi DU_{r_0}(\tilde{X}^{N }(t)). \mathbb g(\tilde{X}^{N }(t),z) \nu(dz)\nonumber \\
&\le \frac{r_0\big(\big\| (-A)^{\frac{1-\epsilon}{4}}\tilde{X}^{N }(t)\big\|^{\frac{1+\epsilon}{2}}\big\| (-A)^{\frac12+\frac{1-\epsilon}{4}}\tilde{X}^{N }(t)\big\|^{\frac{1-\epsilon}{2}}\|F_N(\bar{X}^{N }_m)\|-\|(-A)^{\frac12+\frac{1-\epsilon}{4}}\tilde{X}^{N }(t)\|^2\big)}{\sqrt{1+\|(-A)^{\frac{1-\epsilon}{4}}\tilde{X}^{N }(t)\|^2}}\nonumber\\
&~~~+r_0\int_{\chi}\|(-A)^{\frac{1-\epsilon}{4}}G(\bar{X}^N_m,z)\|\,\nu(dz)+r_0\big\|(-A)^{\frac{1-\epsilon}{4}}B(\bar{X}^{N }_{m})\big\|_{\IL_2^0(H)}^2+r_0^2\big\|Q^{\frac12}(P_NB(\bar{X}^{N }_{m}))^*(-A)^{\frac{1-\epsilon}{4}}\big\|_{\IL(H)}^2\nonumber\\
&\le  \frac{r_0\|(-A)^{\frac{1-\epsilon}{4}}\tilde{X}^{N }(t)\|^{1+\epsilon}\big\| (-A)^{\frac12+\frac{1-\epsilon}{4}}\tilde{X}^{N }(t)\big\|^{1-\epsilon}}{2\sqrt{1+\|(-A)^{\frac{1-\epsilon}{4}}\tilde{X}^{N }(t)\|^2}}- \frac{r_0\|(-A)^{\frac12+\frac{1-\epsilon}{4}}\tilde{X}^{N }(t)\|^2}{\sqrt{1+\|(-A)^{\frac{1-\epsilon}{4}}\tilde{X}^{N }(t)\|^2}}\nonumber\\
&\quad+\frac{r_0\|F(\bar{X}^{N }_{m})\|^2}{2\sqrt{1+\|(-A)^{\frac{1-\epsilon}{4}}\tilde{X}^{N }(t)\|^2}}+C_Gr_0\int_{\chi}\bar{c}(z)\,\nu(dz)+C_Br_0(1+r_0)(1+\|(-A)^{\frac{1-\epsilon}{4}}\bar{X}^{N }_{m}\|^{2\gamma})\nonumber\\
&\le \frac{r_0}{2}\Big(\big\| (-A)^{\frac{1-\epsilon}{4}}\tilde{X}^{N }(t)\big\|^{\epsilon}\big\| (-A)^{\frac12+\frac{1-\epsilon}{4}}\tilde{X}^{N }(t)\big\|^{1-\epsilon}-\frac{\|(-A)^{\frac12+\frac{1-\epsilon}{4}}\tilde{X}^{N }(t)\|^2}{\|(-A)^{\frac{1-\epsilon}{4}}\tilde{X}^{N }(t)\|}\Big)\times \mathbf 1_{\big\{\|(-A)^{\frac{1-\epsilon}{4}}\tilde{X}^{N }(t)\|>1\big\}}\nonumber\\
&\quad+r_0\Big(\frac12\big\| (-A)^{\frac12+\frac{1-\epsilon}{4}}\tilde{X}^{N }(t)\big\|^{1-\epsilon}- \frac{1}{\sqrt{2}}\big\| (-A)^{\frac12+\frac{1-\epsilon}{4}}\tilde{X}^{N }(t)\big\|^2  \Big)\times \mathbf 1_{\big\{\|(-A)^{\frac{1-\epsilon}{4}}\tilde{X}^{N }(t)\|\le 1\big\}}\nonumber\\
&\quad+\frac{r_0(L_F+L_0)^2(1+\|\bar{X}^{N }_{m}\|^2)}{2\sqrt{1+\|(-A)^{\frac{1-\epsilon}{4}}\tilde{X}^{N }(t)\|^2}}+C_{B,G}r_0(1+r_0)(1+\|(-A)^{\frac{1-\epsilon}{4}}\bar{X}^{N }_{m}\|^2)^{\frac12}\\
&\le  \frac{r_0}{2}\big\| (-A)^{\frac{1-\epsilon}{4}}\tilde{X}^{N }(t)\big\|+\frac{r_0(L_F+L_0)^2(1+\|\bar{X}^{N }_{m}\|^2)}{2\sqrt{1+\|(-A)^{\frac{1-\epsilon}{4}}\tilde{X}^{N }(t)\|^2}}+\mathcal{K}(1+r_0)U_{r_0}(\bar{X}^N_m).
\end{align*}


\subsubsection{Proof of the inequalities \eqref{ineq:estimate-64} and \eqref{ineq:estimate-65}.}\label{Subsection:6.1.2}

Recall that $\mathcal J_m(s)$ denotes the sum of the two
$s$-integrands on the right-hand side of \eqref{estimate-Ito}. Set
\[
\Lambda_m(s)
:=
e^{-\varsigma_4(s-t_m)}
\exp\Big(
\frac{U_{r_0}(\tilde X^N(s))}
{e^{\varsigma_4(s-t_m)}}
-
\int_{t_m}^{s}
\frac{\vartheta_2}{e^{\varsigma_4(r-t_m)}}\,dr
\Big)>0.
\]
By Assumption~\ref{assum:nonlinearity-G} (i), together with the commutativity and
contractivity of $E(t)$ and $P_N$, $\|(-A)^{\frac{1-\epsilon}{4}}E(s-t_m)G_N(\bar X_m^N,z)\|
\le\kappa_\epsilon\bar c(z).$ Moreover,  $U_{r_0}(x+g)-U_{r_0}(x)
\le r_0\|(-A)^{\frac{1-\epsilon}{4}}g\|$ and $|DU_{r_0}(x+g).g|
\le r_0\|(-A)^{\frac{1-\epsilon}{4}}g\|,$ the Poisson remainder in $\mathcal J_m(s)$ is bounded by
\[
\Lambda_m(s)
\int_\chi
e^{r_0\kappa_\epsilon\bar c(z)}
r_0\kappa_\epsilon\bar c(z)\,\nu(dz)
=
\Lambda_m(s)(\vartheta_2-k_1r_0).
\]
Therefore, by the definition of $\varsigma_4$ and the fact that
\[
-\frac{r_0}{2}
\sqrt{
1+\|(-A)^{\frac{1-\epsilon}{4}}\tilde X^N(s)\|^2
}
+
\frac{r_0}{2}
\|(-A)^{\frac{1-\epsilon}{4}}\tilde X^N(s)\|
\le0,
\]
we obtain, for $(ds\otimes d\mathbb P)$-a.e.
$(s,\omega)\in[t_m,t_{m+1}]\times\Omega$,
\begin{align}\label{6.10}
&\mathcal J_m(s)
\le
\Lambda_m(s)
\bigg[
-k_1r_0
-2r_0(L_F+L_0)^2
\left(
1+\lambda_1^{-\frac{1-\epsilon}{2}}
e^{\frac{\sqrt{\ln2}}2}
\right)\sqrt{
1+\|(-A)^{\frac{1-\epsilon}{4}}\tilde X^N(s)\|^2}\nonumber\\
&+\frac{
r_0(L_F+L_0)^2(1+\|\bar X_m^N\|^2)
}{
2\sqrt{
1+\|(-A)^{\frac{1-\epsilon}{4}}\tilde X^N(s)\|^2
}
}+
\mathcal K r_0(1+r_0)
\left(
\sqrt{
1+\|(-A)^{\frac{1-\epsilon}{4}}\bar X_m^N\|^2
}
-
2e^{\frac{\sqrt{\ln2}}4}
\sqrt{
1+\|(-A)^{\frac{1-\epsilon}{4}}\tilde X^N(s)\|^2
}
\right)
\bigg].
\end{align}

We first prove \eqref{ineq:estimate-64}.  On the event $\mathcal A_m(s)$,
\begin{equation}\label{eq:A-event-bound}
\|\mathrm{I}_{1,m}(s)\|
\le
2\|(-A)^{\frac{1-\epsilon}{4}}\tilde X^N(s)\|.
\end{equation}
Moreover, by \eqref{ineq:compare-s,t_m}, $\|(-A)^{\frac{1-\epsilon}{4}}\bar X_m^N\|
\le
e^{\frac{\sqrt{\ln2}}4}\|\mathrm{I}_{1,m}(s)\|.$ By the continuous embedding
$\dot H^{\frac{1-\epsilon}{2}}\hookrightarrow H$, it holds that
\begin{align*}
\frac{1+\|\bar X_m^N\|^2}
{2\sqrt{1+\|(-A)^{\frac{1-\epsilon}{4}}
\tilde X^N(s)\|^2}}
&\le
\frac{
1+\lambda_1^{-\frac{1-\epsilon}{2}}
e^{\frac{\sqrt{\ln2}}2}
\|\mathrm I_{1,m}(s)\|^2
}{
\sqrt{1+\|\mathrm I_{1,m}(s)\|^2}
}\le
\left(
1+\lambda_1^{-\frac{1-\epsilon}{2}}
e^{\frac{\sqrt{\ln2}}2}
\right)
\sqrt{1+\|\mathrm I_{1,m}(s)\|^2}.
\end{align*}
Using \eqref{ineq:compare-s,t_m}, \eqref{eq:A-event-bound}, and the preceding estimate in \eqref{6.10}, and discarding the non-positive term $-k_1r_0$, 
\begin{align*}
\mathcal J_m(s)
&\le
\Lambda_m(s)\,
r_0(L_F+L_0)^2\left(1+\lambda_1^{-\frac{1-\epsilon}{2}}
e^{\frac{\sqrt{\ln2}}2}\right)
\Bigl(\sqrt{1+\|\mathrm I_{1,m}(s)\|^2}-2\sqrt{1+\|(-A)^{\frac{1-\epsilon}{4}}\tilde X^N(s)\|^2}\Bigr)
\\&~\quad +
\Lambda_m(s)\,
\mathcal K r_0(1+r_0)
\Bigl(
\sqrt{
1+\|(-A)^{\frac{1-\epsilon}{4}}
\bar X_m^N\|^2
}
-
e^{\frac{\sqrt{\ln2}}4}
\sqrt{1+\|\mathrm I_{1,m}(s)\|^2}
\Bigr)
\\
&\le0 .
\end{align*}
Indeed, the first bracket is non-positive by \eqref{eq:A-event-bound}, whereas the second one is non-positive by
\eqref{ineq:compare-s,t_m}. Therefore,
\[
\mathcal J_m(s)\mathbf 1_{\mathcal A_m(s)}
\le 0,
\qquad ds\otimes d\mathbb P\text{-a.e.}
\]
which proves \eqref{ineq:estimate-64}.

We next prove \eqref{ineq:estimate-65}. On the event $\mathcal B_m(s)$, $\|(-A)^{\frac{1-\epsilon}{4}}\tilde X^N(s)\|
<
\frac12\|I_{1,m}(s)\|.$ 
Consequently,
\[
\Lambda_m(s)
\le
\widehat\Lambda_m(s)
:=
e^{-\varsigma_4(s-t_m)}
\exp\Big(
\frac{
r_0\sqrt{1+\frac14\|\mathrm I_{1,m}(s)\|^2}
}{
e^{\varsigma_4(s-t_m)}
}
-
\int_{t_m}^{s}
\frac{\vartheta_2}
{e^{\varsigma_4(r-t_m)}}\,dr\Big).
\]
Moreover, by \eqref{ineq:compare-s,t_m} and the continuous embedding $\dot H^{\frac{1-\epsilon}{2}}\hookrightarrow H$, we obtain $\|(-A)^{\frac{1-\epsilon}{4}}\bar X_m^N\|
\le
e^{\frac{\sqrt{\ln2}}4}\|\mathrm{I}_{1,m}(s)\|$ and $\|\bar X_m^N\|^2
\le
\lambda_1^{-\frac{1-\epsilon}{2}}
e^{\frac{\sqrt{\ln2}}2}
\|\mathrm I_{1,m}(s)\|^2.$ Therefore, by \eqref{6.10}, after discarding the non-positive terms,
\begin{align*}
\mathcal J_m(s)\mathbf 1_{\mathcal B_m(s)}
&\le
\widehat\Lambda_m(s)
\Bigg[
\frac{r_0(L_F+L_0)^2}{2}
\left(
1+
\lambda_1^{-\frac{1-\epsilon}{2}}
e^{\frac{\sqrt{\ln2}}2}
\|\mathrm I_{1,m}(s)\|^2
\right)+
\mathcal K r_0(1+r_0)
e^{\frac{\sqrt{\ln2}}4}
\sqrt{1+\|\mathrm I_{1,m}(s)\|^2}
\Bigg]
\mathbf 1_{\mathcal B_m(s)}
\\
&\le
\widehat\Lambda_m(s)
\mathcal K_1 r_0(1+r_0)
\bigl(1+\|\mathrm I_{1,m}(s)\|^2\bigr)
\mathbf 1_{\mathcal B_m(s)},
\end{align*}
where $\mathcal K_1>0$ is independent of $r_0$.
Hence, by $\widehat{\mathcal J}_m(s)
=
\widehat\Lambda_m(s)
\mathcal K_1 r_0(1+r_0)
\bigl(1+\|\mathrm I_{1,m}(s)\|^2\bigr),$ we obtain
\begin{equation*}
\mathcal J_m(s)\mathbf 1_{\mathcal B_m(s)}
\le
\widehat{\mathcal J}_m(s)
\mathbf 1_{\mathcal B_m(s)},
\qquad
ds\otimes d\mathbb P\text{-a.e.},
\end{equation*}
which proves \eqref{ineq:estimate-65}.

\subsubsection{Proof of the inequality \eqref{estimate-Theta-1}.} \label{Subsection:6.1.3}

For $s\in[t_m,t_{m+1}]$ set
\begin{align}\label{set:WY}
	\!\!\!\!\!\!\IW_{m}(\tilde{X}^{N })(s):=E(s-t_m)P_NB(\bar{X}^{N }_{m})(W^Q(s)-W^Q(t_m)),\quad Y_m(\tilde{X}^{N })(s):=\tilde{X}^{N }(s)-\IW_m(\tilde{X}^{N})(s),
\end{align}
and  introduce the event
\begin{align}\label{ineq:estimate-A}
\mathcal C_m:=\Big\{ &\| \mathrm I_{1,m}(s)\|>C_L, ~\big\|\int_{t_m}^s\!\int_\chi (-A)^{\frac{1-\epsilon}{4}}G_N(\bar{X}^N_m,z) \tilde{\mathbf{N}}(dz,dr)\big\|\le \frac{1}{k_1},\nonumber\\
 &~~~\|(-A)^{\frac{1-\epsilon}{4}}Y_m(\tilde{X}^{N})(s)\|\le \Big(\frac34-\frac{\tau^{\frac12}}{4\lambda_1^{\frac{1-\epsilon}{4}}}\Big)\|\mathrm I_{1,m}(s)\|+1\Big\}.	
\end{align}
Then we claim that $\E\big[\mathbf{1}_{\mathcal C_m}\big| \mathcal F_{t_m}\big]\!=\!0$.
Indeed, on the event \(\mathcal C_m\),  the first two conditions in its definition, together with \eqref{set:WY}, yield
\begin{align}\label{ineq:C_F}
	&\quad~~~ \quad \big\|(-A)^{\frac{1-\epsilon}{4}}Y_m(\tilde{X}^{N })(s)-E(s-t_m)(-A)^{\frac{1-\epsilon}{4}}\bar{X}_m^{N }\big\|\nonumber\\
	&=\Big\|(s-t_m)E(s-t_m)(-A)^{\frac{1-\epsilon}{4}}F_N(\bar{X}^{N }_{m})+\int_{t_m}^s\int_\chi E(s-t_m)(-A)^{\frac{1-\epsilon}{4}}G_N(\bar{X}^N_m,z)\,\tilde{\mathbf{N}}(dz,dr)\Big\|\nonumber\\
	&\le M_F\tau^{\frac{3+\epsilon}{4}}(1\!+\!\|\bar{X}^{N }_{m}\|)+\frac{1}{k_1}
	\le \Big(\frac{e^{-\frac{\sqrt{\ln 2}}{4}}}{4}\tau^{\frac12}+\frac12\Big)+\frac{e^{-\frac{\sqrt{\ln 2}}{4}}}{4}\tau^{\frac12}\lambda_1^{-\frac{1-\epsilon}{4}}\big\|(-A)^{\frac{1-\epsilon}{4}}\bar{X}_m^{N }\big\|\nonumber\\
	&\le \frac{C_L}{2}e^{-\frac{\sqrt{\ln 2}}{4}}+\frac{\|E(s-t_m)(-A)^{\frac{1-\epsilon}{4}}\bar{X}_m^{N }\|}{2}\times e^{-\frac{\sqrt{\ln 2}}{4}}\le \big\|E(s-t_m)(-A)^{\frac{1-\epsilon}{4}}\bar{X}_m^{N}\big\|,
\end{align}
according to the fact $k_1\ge2$, the expression of $\tau_N$ given in  \eqref{def:tau_0}, and  \eqref{ineq:compare-s,t_m}. From  $C_L\ge  7 $, it follows that 
\begin{align*}
	&\|(-A)^{\frac{1-\epsilon}{4}}Y_m(\tilde{X}^{N })(s)\|=\Big\|E(s-t_m)(-A)^{\frac{1-\epsilon}{4}}\bar{X}_m^{N }+(-A)^{\frac{1-\epsilon}{4}}Y_m(\tilde{X}^{N })(s)-E(s-t_m)(-A)^{\frac{1-\epsilon}{4}}\bar{X}_m^{N }\Big\|\\
	&\ge \big\|E(s-t_m)(-A)^{\frac{1-\epsilon}{4}}\bar{X}_m^{N }\big\|-\big\|(-A)^{\frac{1-\epsilon}{4}}Y_m(\tilde{X}^{N })(s)-E(s-t_m)(-A)^{\frac{1-\epsilon}{4}}\bar{X}_m^{N }\big\|\\
	&\ge \big\|E(s\!-\!t_m)(-A)^{\frac{1-\epsilon}{4}}\bar{X}_m^{N }\big\|-M_F\tau^{\frac{3+\epsilon}{4}}(1+\|\bar{X}^{N }_{m}\|)-\frac{1}{k_1}\ge \Big(1\!-\!\frac{\tau^{\frac12}}{4\lambda_1^{\frac{1-\epsilon}{4}}}\Big) \big\|E(s\!-\!t_m)(-A)^{\frac{1-\epsilon}{4}}\bar{X}_m^{N }\big\|-\frac{\tau^{\frac12}+2}{4} \\
	&> \Big(1-\frac{\tau^{\frac12}}{4\lambda_1^{\frac{1-\epsilon}{4}}}-\frac14\Big)\|\mathrm I_{1,m}(s)\|+ \frac{C_L-3}{4}\ge \Big(\frac34-\frac{\tau^{\frac12}}{4\lambda_1^{\frac{1-\epsilon}{4}}}\Big) \|\mathrm I_{1,m}(s)\|+1,
\end{align*}
 which contradicts the third condition in the definition of $\mathcal C_m$. Hence $\mathbb P(\mathcal C_m)=0$. Then $\mathbb E[\mathbf 1_{\mathcal C_m}|\mathcal F_{t_m}]=0,
~\mathbb P\text{-a.s.}$

In addition, under Assumption~\ref{assum:nonlinearity-G} (i), by Chebyshev's inequality for the conditional expectation and the independent increments property of the compensated Poisson random measure, we obtain
\begin{align}\label{ineq:estimate-SI-N}
&\quad  \E \Big[ \mathbf{1}_{\big\{ \|\int_{t_m}^s\!\int_\chi (-A)^{\frac{1-\epsilon}{4}}G_N(\bar{X}^N_m,z) \tilde{\mathbf{N}}(dz,dr)\|> \frac{1}{k_1}\big\}} ~\big| \mathcal F_{t_m}\Big]\nonumber\\
&\le k_1^2\times\E\Big[\Big\|\int_{t_m}^s\!\int_\chi (-A)^{\frac{1-\epsilon}{4}}G(u,z) \,\tilde{\mathbf{N}}(dz,dr)\Big\|^2\Big]~\Big|_{u=\bar{X}^N_m}\le C(s-t_m).
\end{align}

Now using the relationship \eqref{set:WY} between $\tilde{X}^{N }$ and $Y_m(\tilde{X}^N)$, we decompose the estimate of $\Theta_{m,m+1}^1$ as follows
\begin{align*}
	\Theta_{m,m+1}^1&\le  \E\bigg[ \int_{t_m}^{t_{m+1}} \!\!\!\exp\Big(\frac{r_0\sqrt{1+\big(\frac12\|\mathrm I_{1,m}(s)\|\big)^2}}{e^{\varsigma_4(s-t_m)}}-\int_{t_m}^{s}\frac{\vartheta_2}{e^{\varsigma_4(r-t_m)}}dr\Big)e^{-\varsigma_4(s-t_m)}  \mathcal{K}_1r_0(1\!+\!r_0)(1+\|\mathrm I_{1,m}(s)\|^2)\nonumber\\
	&\quad\quad\times\mathbf 1_{\big\{ \|(-A)^{\frac{1-\epsilon}{4}}Y_m(\tilde{X}^{N})(s)\|\le \big(\frac34-\frac{\tau^{\frac12}}{4\lambda_1^{\frac{1-\epsilon}{4}}}\big)\|\mathrm I_{1,m}(s)\|+1,~  \| \mathrm I_{1,m}(s)\|> C_L\big\}}\,ds~\Big| \mathcal F_{t_m}\bigg]\nonumber\\
	&+ \E\bigg[ \int_{t_m}^{t_{m+1}}\!\!\!\exp\Big(\frac{r_0\sqrt{1+\big(\frac12\|\mathrm I_{1,m}(s)\|\big)^2}}{e^{\varsigma_4(s-t_m)}}-\int_{t_m}^{s}\frac{\vartheta_2}{e^{\varsigma_4(r-t_m)}}dr\Big)e^{-\varsigma_4(s-t_m)} \mathcal{K}_1r_0(1\!+\!r_0)(1+\|\mathrm I_{1,m}(s)\|^2)\nonumber\\
	&\quad\quad\times \mathbf 1_{\big\{ \|(-A)^{\frac{1-\epsilon}{4}}\IW_{m}(\tilde{X}^{N})(s)\|\ge 1+\big(\frac14-\frac{\tau^{\frac12}}{4\lambda_1^{\frac{1-\epsilon}{4}}}\big) \|\mathrm I_{1,m}(s)\|,~  \| \mathrm I_{1,m}(s)\|> C_L\big\}}\,ds~\Big| \mathcal F_{t_m}\bigg]\nonumber.
\end{align*}
Moreover, dividing
 $\bar{B}:=\big\{ \|(-A)^{\frac{1-\epsilon}{4}}Y_m(\tilde{X}^{N})(s)\|\le \big(\frac34-\frac{\tau^{\frac12}}{4\lambda_1^{\frac{1-\epsilon}{4}}}\big)\|\mathrm I_{1,m}(s)\|+1,~  \| \mathrm I_{1,m}(s)\|> C_L\big\}$
  into $\mathcal C_m$ and $\bar{B}\setminus \mathcal C_m$, and combining the results \eqref{ineq:estimate-A} and \eqref{ineq:estimate-SI-N}, we derive
\begin{align*}
	\Theta_{m,m+1}^1&\le \int_{t_m}^{t_{m+1}} \exp\Big(\frac{r_0\sqrt{1+\big(\frac12\|\mathrm I_{1,m}(s)\|\big)^2}}{e^{\varsigma_4(s-t_m)}}-\int_{t_m}^{s}\frac{\vartheta_2}{e^{\varsigma_4(r-t_m)}}dr\Big)e^{-\varsigma_4(s-t_m)}\nonumber\\
&\quad\times 2\mathcal{K}_1r_0(1+r_0) \Big(C(s-t_m)\|\mathrm I_{1,m}(s)\|^2+16\big(1-\frac{\tau^{\frac12}}{\lambda_1^{\frac{1-\epsilon}{4}}}\big)^{-2}\times\E\big[ \|(-A)^{\frac{1-\epsilon}{4}}\IW_{m}(\tilde{X}^{N})(s)\|^{2}\big|\mathcal{F}_{t_m}\big]\Big)\,ds\nonumber\\
	&\le \int_{t_m}^{t_{m+1}} \exp\Big(\frac{r_0\sqrt{1+\big(\frac12\|(-A)^{\frac{1-\epsilon}{4}}\bar{X}^{N}_{m}\|\big)^2}}{e^{\varsigma_4(s-t_m)}}-\int_{t_m}^{s}\frac{\vartheta_2}{e^{\varsigma_4(r-t_m)}}dr\Big)e^{-\varsigma_4(s-t_m)}\\
	&\quad \times Cr_0(1+r_0)(s-t_m)\big(1+\|(-A)^{\frac{1-\epsilon}{4}}\bar{X}^{N}_{m}\|^{2\gamma}+\|(-A)^{\frac{1-\epsilon}{4}}\bar{X}^{N}_{m}\|^2\big)\,ds.\nonumber
\end{align*}
The proof of inequality \eqref{estimate-Theta-1} is completed.

\subsection{Proof of Lemma~\ref{lem:estimate-eta}}

Define $\tilde{\eta}^{N }_{y,x,m}:=\eta^{N }_{y,x,m}-E(m\tau)y$, $m\in\{1,2,\ldots,M\}$ such that $\tilde{\eta}^{N }_{y,x,0}=0$. 
Introduce a process $\{\bar{\eta}^{N}_{y,x}(t)\}_{t\in[0,T]}$ that satisfies $\bar{\eta}^{N}_{y,x}(t_m)=\tilde{\eta}^{N}_{y,x,m}, \bar{\eta}^{N}_{y,x}(t_{m-1})=\tilde{\eta}^{N}_{y,x,m-1}$, and 
\begin{align*}
	d\bar{\eta}^{N}_{y,x}(t)&=A\bar{\eta}^{N}_{y,x}(t)dt+E(t-t_{m-1})P_NDF(\bar{X}^{N}_{m-1}).(\bar{\eta}^{N}_{y,x,m-1}\!+\!E(t_{m-1})y)\,dt\\
	&\quad+E(t-t_{m-1})P_NDB(\bar{X}^{N}_{m-1}).(\bar{\eta}^{N}_{y,x,m-1}\!+\!E(t_{m-1})y)\,dW^Q(t)\\
	&\quad+E(t-t_{m-1})\int_\chi P_NDG(\bar{X}^{N}_{m-1},z).\big(\bar{\eta}^{N}_{y,x,m-1}\!+\!E(t_{m-1})y\big)\,\tilde{\mathbf{N}}(dz,dt),\quad t\in[t_{m-1},t_{m}].
\end{align*}
For all $t\in[t_{m-1}, t_m]$, define $Y^{N}_{y,x,m-1}(t):=\bar{\eta}^{N}_{y,x}(t)-\vartheta^{N}_{y,x,m-1}(t)-\mathcal N^{N}_{y,x,m-1}(t)$ with
\begin{align*}
	&\vartheta^{N}_{y,x,m-1}(t):=\int_{t_{m-1}}^t\!\!E(t-t_{m-1})P_NDB(\bar{X}^{N}_{m-1}).(\bar{\eta}^{N}_{y,x,m-1}+E(t_{m-1})y)\,dW^Q(s),\\
	&\mathcal N^{N}_{y,x,m-1}(t):=\int_{t_{m-1}}^t\!\!E(t-t_{m-1})\int_\chi P_NDG(\bar{X}^{N}_{m-1},z).\big(\bar{\eta}^{N}_{y,x,m-1}\!+E(t_{m-1})y\big)\,\tilde{\mathbf{N}}(dz,ds).
\end{align*}
They are solutions of the following equations, respectively, 
\begin{align*}
	\!\!d\vartheta^{N}_{y,x,m-1}(t)\!&=\!A \vartheta^{N}_{y,x,m-1}(t)dt+E(t\!-\!t_{m-1})P_NDB(\bar{X}^{N}_{m-1}).(\bar{\eta}^{N}_{y,x,m-1}+E(t_{m-1})y)dW^Q(t),\quad \vartheta^{N}_{y,x,m-1}(t_{m\!-\!1})=0;\\
	\!\!d\mathcal N^{N}_{y,x,m\!-\!1}(t)\!&=\!A\mathcal N^{N}_{y,x,m\!-\!1}(t)dt\!+\!E(t\!-\!t_{m-1})\int_\chi \!\! P_NDG(\bar{X}^{N}_{m-1},z).\big(\bar{\eta}^{N}_{y,x,m-1}\!+\!E(t_{m-1})y\big)\tilde{\mathbf{N}}(dz,dt),~\mathcal N^{N}_{y,x,m-1}(t_{m\!-\!1})=0;\\
	\!\!dY^{N}_{y,x,m-1}(t)\!&=\!AY^{N}_{y,x,m-1}(t)dt+E(t\!-\!t_{m-1})P_NDF(\bar{X}^{N}_{m-1}).\big(\bar{\eta}^{N}_{y,x,m-1}+E(t_{m-1})y\big)dt,\quad Y^{N}_{y,x,m\!-\!1}(t_{m\!-\!1})=\bar{\eta}^{N}_{y,x,m\!-\!1}.
\end{align*}
Under Assumption~\ref{assum:nonlinearity-B} (ii), applying It\^o's formula to $\|(-A)^{-\frac{1-\epsilon}{4}}\vartheta^{N}_{y,x,m-1}(t)\|^{p}$, $p\ge 2$ yields
\begin{align*}
&\big\|(-A)^{-\frac{1-\epsilon}{4}}\vartheta^{N}_{y,x,m-1}(t_m)\big\|^p\le \int_{t_{m-1}}^{t_m}	p\|(-A)^{-\frac{1-\epsilon}{4}}\vartheta^{N}_{y,x,m-1}(t)\|^{p-2}\big(-\lambda_1\big\| (-A)^{-\frac{1-\epsilon}{4}}\vartheta^{N}_{y,x,m-1}(t)\big\|^2\big)dt\\
&+\frac{p(p-1)}{2}\int_{t_{m-1}}^{t_m}	\|(-A)^{-\frac{1-\epsilon}{4}}\vartheta^{N}_{y,x,m-1}(t)\|^{p-2}dt\times \big\|(-A)^{-\frac{1-\epsilon}{4}}DB(\bar{X}^{N}_{m-1}).\big(\bar{\eta}^N_{y,x,m-1}+E(t_{m-1})y\big)\big\|_{\IL_2^0(H)}^2\\
&+p\int_{t_{m-1}}^{t_m}\!\!\!\|(-A)^{-\frac{1-\epsilon}{4}}\vartheta^{N}_{y,x,m-1}(t)\|^{p-2}\big\langle (-A)^{-\frac{1-\epsilon}{4}}\vartheta^{N}_{y,x,m-1}(t),\\
&\quad E(t\!-\!t_{m-1})(-A)^{-\frac{1-\epsilon}{4}}P_NDB(\bar{X}^{N}_{m-1}).(\bar{\eta}_{y,x,m-1}^N\!+\!E(t_{m-1})y)dW^Q(t)\big\rangle\\
&\le C_{p} \tau  \big(\big\|(-A)^{-\frac{1-\epsilon}{4}}\bar{\eta}_{y,x,m-1}^N\big\|^p+\big\|(-A)^{-\frac{1-\epsilon}{4}}E(t_{m-1})y\big\|^p\big)+p\int_{t_{m-1}}^{t_m}\!\!\!\|(-A)^{-\frac{1-\epsilon}{4}}\vartheta^{N}_{y,x,m-1}(t)\|^{p-2}\\
&\quad \times\big\langle (-A)^{-\frac{1-\epsilon}{4}}\vartheta^{N}_{y,x,m-1}(t), E(t-t_{m-1})(-A)^{-\frac{1-\epsilon}{4}}P_NDB(\bar{X}^{N}_{m-1}).(\bar{\eta}_{y,x,m-1}^N\!+E(t_{m-1})y)dW^Q(t)\big\rangle,
\end{align*}
where we used the Young inequality in the last step. 
Using Taylor's formula, one has
\begin{align}\label{eq:Taylor}
	\|x+y\|^p-\|x\|^p-p\|x\|^{p-2}\langle x,y\rangle \le K_p\big(\|x\|^{p-2}\|y\|^2+\|y\|^p\big) 
\end{align}
with $K_p:=2^{p-4}p(p-1)$ for $p\ge 3$ and $K_p:=\frac{p(p-1)}{2}$ for $2\le p<3$. A similar calculation to $\|(-A)^{-\frac{1\!-\!\epsilon}{4}}\mathcal N^{N}_{y,x,m\!-\!1}(t)\|^{p}$ yields that there exists a positive constant $C_p$ such that
\begin{align*}
&\|(-A)^{-\frac{1-\epsilon}{4}}\mathcal N^{N}_{y,x,m-1}(t_m)\|^{p}\le \int_{t_{m-1}}^{t_m}\Big((-p\lambda_1)\times\|(-A)^{-\frac{1-\epsilon}{4}}\mathcal N^{N}_{y,x,m-1}(t)\|^{p-2}\big\|(-A)^{-\frac{1-\epsilon}{4}}\mathcal N^{N}_{y,x,m-1}(t)\big\|^2\\
	&+K_p\!\int_\chi\!\!\|(-A)^{-\frac{1-\epsilon}{4}}\mathcal N^{N}_{y,x,m-1}(t)\|^{p-2}\big\|(-A)^{-\frac{1-\epsilon}{4}}DG(\bar{X}^{N}_{m-1},z).(\bar{\eta}_{y,x,m-1}^N\!+\!E(t_{m-1})y)\big\|^2 \Big)\nu(dz)dt\\
	&+K_p\int_{t_{m-1}}^{t_m}\int_\chi\big\|(-A)^{-\frac{1-\epsilon}{4}}DG(\bar{X}^{N}_{m-1},z).\big(\bar{\eta}_{y,x,m-1}^N+E(t_{m-1})y\big)\big\|^{p}\,\nu(dz)dt\\
	&+\int_{t_{m-1}}^{t_m}\!\!\int_\chi\!\! \Big(\big\|(-A)^{-\frac{1-\epsilon}{4}}\mathcal N^{N}_{y,x,m-1}(t)+E(t-t_{m-1})(-A)^{-\frac{1-\epsilon}{4}}P_NDG(\bar{X}^{N}_{m-1},z).(\bar{\eta}_{y,x,m-1}^N\!+\!E(t_{m-1})y)\big\|^{p}\\
	&\quad -\!\big\|(-A)^{-\frac{1-\epsilon}{4}}\mathcal N^{N}_{y,x,m-1}(t)\big\|^{p}\Big)\,\tilde{\mathbf{N}}(dz,dt)\\
	&\le C_{p}\tau  \big(\big\|(-A)^{-\frac{1-\epsilon}{4}}\bar{\eta}_{y,x,m-1}^N\big\|^p+\|E(t_{m-1})(-A)^{-\frac{1-\epsilon}{4}}y\|^p\big)+\int_{t_{m-1}}^{t_m}\!\!\int_\chi\! \Big(\big\|(-A)^{-\frac{1-\epsilon}{4}}\mathcal N^{N}_{y,x,m-1}(t)\\
	&\!\!+E(t\!-\!t_{m-1})(-A)^{-\frac{1-\epsilon}{4}}P_NDG(\bar{X}^{N}_{m-1},z).(\bar{\eta}_{y,x,m-1}^N\!+\!E(t_{m-1})y)\big\|^{p}\!-\!\big\|(-A)^{-\frac{1-\epsilon}{4}}\mathcal N^{N}_{y,x,m-1}(t)\big\|^{p}\Big)\,\tilde{\mathbf{N}}(dz,dt).
\end{align*}
Moreover, for $\alpha_0\in(0,1)$ given by Assumption~\ref{assum:nonlinearity-F} (ii), we claim that there exists a positive constant $C_{\epsilon}$ dependent on $p$ such that
\begin{align}\label{estimate-Y}
	&\|(-A)^{-\frac{1-\epsilon}{4}}Y^{N}_{y,x,m-1}(t_m)\|^{p}-\|(-A)^{-\frac{1-\epsilon}{4}}\bar{\eta}_{y,x,m-1}^N\|^{p}
	\le \big( \|(-A)^{-\frac{1-\epsilon}{4}}\bar{\eta}_{y,x,m-1}^N\|^{p}+\big\|E(t_{m-1})(-A)^{-\frac{1-\epsilon}{4}}y\big\|^p\big)\nonumber\\
	&\quad\quad\times\Big\{C_{\epsilon}\tau\big(1+\|(-A)^{\frac{1-\epsilon}{4}}\bar{X}^{N}_{m-1}\|\big)^{\alpha_0(1+\epsilon)} +C_{\epsilon}\big(\tau\big(1+\|(-A)^{\frac{1-\epsilon}{4}}\bar{X}^{N}_{m-1}\|\big)^{\alpha_0(1+\epsilon)}\big)^{\frac{p}{1+\epsilon}}\Big\}.
\end{align}
The proof of \eqref{estimate-Y} is postponed to Subsection~\ref{Subsection:6.2.1}. For later use, set
\[
\mathcal{Q}_{m-1}
:=
\tau
\big(
1+\big\|(-A)^{\frac{1-\epsilon}{4}}
\bar X_{m-1}^N\big\|
\big)^{\alpha_0(1+\epsilon)},
\]
and, for every \(r\ge0\), define
\begin{align}\label{eq:notation-G}
	\mathcal G_{m-1,r}:=\mathcal{Q}_{m-1}+\mathcal{Q}_{m-1}^{\frac{r}{1+\epsilon}}.
\end{align}
Then it holds that 
\begin{align*}
\|(-A)^{-\frac{1-\epsilon}{4}}Y^{N}_{y,x,m-1}(t_m)\|^{p}\le \|(-A)^{-\frac{1-\epsilon}{4}}\bar{\eta}_{y,x,m-1}^N\|^{p} (1+C_{\epsilon}\mathcal G_{m-1,p})+\big\|E(t_{m-1})(-A)^{-\frac{1-\epsilon}{4}}y\big\|^p \times C_{\epsilon}\mathcal G_{m-1,p}.	
\end{align*}
Using \eqref{eq:Taylor} and combining estimates obtained for $Y^{N}_{y,x,m-1}(t_m), \vartheta^{N}_{y,x,m-1}(t_m)$ and $\mathcal N^{N}_{y,x,m-1}(t_m)$ above yield 
\begin{align*}
&\|(-A)^{-\frac{1-\epsilon}{4}}\bar{\eta}^{N}_{y,x,m}\|^{p}=\big\|(-A)^{-\frac{1-\epsilon}{4}}Y^{N}_{y,x,m-1}(t_m)+(-A)^{-\frac{1-\epsilon}{4}}\vartheta^{N}_{y,x,m-1}(t_m)+(-A)^{-\frac{1-\epsilon}{4}}\mathcal N^{N}_{y,x,m-1}(t_m)\big\|^p\\
&\le p\|(-A)^{-\frac{1-\epsilon}{4}}Y^{N}_{y,x,m-1}(t_m)\|^{p-2}\big\langle(-A)^{-\frac{1-\epsilon}{4}}Y^{N}_{y,x,m-1}(t_m),  (-A)^{-\frac{1-\epsilon}{4}}\vartheta^{N}_{y,x,m-1}(t_m)+(-A)^{-\frac{1-\epsilon}{4}}\mathcal N^{N}_{y,x,m-1}(t_m)\big\rangle \\
&~~~+K_p\|(-A)^{-\frac{1-\epsilon}{4}}Y^{N}_{y,x,m-1}(t_m)\|^{p-2}\big\|(-A)^{-\frac{1-\epsilon}{4}}\vartheta^{N}_{y,x,m-1}(t_m)+(-A)^{-\frac{1-\epsilon}{4}}\mathcal N^{N}_{y,x,m-1}(t_m)\big\|^2\\
&~~~+\|(-A)^{-\frac{1-\epsilon}{4}}Y^{N}_{y,x,m-1}(t_m)\|^p+K_p\big\|(-A)^{-\frac{1-\epsilon}{4}}\vartheta^{N}_{y,x,m-1}(t_m)+(-A)^{-\frac{1-\epsilon}{4}}\mathcal N^{N}_{y,x,m-1}(t_m)\big\|^p\\
&\le \|(-A)^{-\frac{1-\epsilon}{4}}\bar{\eta}_{y,x,m-1}^N\|^{p}\times \Big\{1+C_{\epsilon}\big(\mathcal G_{m-1,p}+\tau \mathcal G_{m-1,p-2}+\tau\big)\Big\}\\
&+\big\|E(t_{m-1})(-A)^{-\frac{1-\epsilon}{4}}y\big\|^p\times C_{\epsilon}\big(\mathcal G_{m-1,p}+\tau \mathcal G_{m-1,p-2}+\tau\big)\\
&+\bigg\{p\|(-A)^{-\frac{1-\epsilon}{4}}Y^{N}_{y,x,m-1}(t_m)\|^{p-2}\big\langle(-A)^{-\frac{1-\epsilon}{4}}Y^{N}_{y,x,m-1}(t_m),  (-A)^{-\frac{1-\epsilon}{4}}\vartheta^{N}_{y,x,m-1}(t_m)+(-A)^{-\frac{1-\epsilon}{4}}\mathcal N^{N}_{y,x,m-1}(t_m)\big\rangle\\
&+C_p
\|(-A)^{-\frac{1-\epsilon}{4}}Y^{N}_{y,x,m-1}(t_m)\|^{p-2}\bigg(\int_{t_{m-1}}^{t_m}\!\!\!\!\big\langle (-A)^{-\frac{1-\epsilon}{4}}\vartheta^{N}_{y,x,m-1}(t), E(t\!-\!t_{m-1})(-A)^{-\frac{1-\epsilon}{4}}P_NDB(\bar{X}^{N}_{m-1}).(\bar{\eta}_{y,x,m-1}^N\\
&+E(t_{m-1})y)dW^Q(t)\big\rangle+\int_{t_{m-1}}^{t_m}\int_\chi \Big(\big\|(-A)^{-\frac{1-\epsilon}{4}}\mathcal N^{N}_{y,x,m-1}(t)\\
&+E(t\!-\!t_{m-1})(-A)^{-\frac{1-\epsilon}{4}}P_NDG(\bar{X}^{N}_{m-1},z).(\bar{\eta}_{y,x,m-1}^N+E(t_{m-1})y)\big\|^{2}-\big\|(-A)^{-\frac{1-\epsilon}{4}}\mathcal N^{N}_{y,x,m-1}(t)\big\|^{2}\Big)\tilde{\mathbf{N}}(dz,dt)\bigg)\\
&+C_p\int_{t_{m-1}}^{t_m}\!\!\!\|(-A)^{-\frac{1-\epsilon}{4}}\vartheta^{N}_{y,x,m-1}(t)\|^{p-2}\big\langle (-A)^{-\frac{1-\epsilon}{4}}\vartheta^{N}_{y,x,m-1}(t), \\
&\quad E(t\!-\!t_{m-1})(-A)^{-\frac{1-\epsilon}{4}}P_NDB(\bar{X}^{N}_{m-1}).(\bar{\eta}_{y,x,m-1}^N\!\!+\!E(t_{m-1})y)dW^Q(t)\big\rangle+C_p\int_{t_{m-1}}^{t_m}\!\!\int_\chi \!\Big(\big\|(-A)^{-\frac{1-\epsilon}{4}}\mathcal N^{N}_{y,x,m-1}(t)\\
&+E(t\!-\!t_{m-1})(-A)^{-\frac{1-\epsilon}{4}}P_NDG(\bar{X}^{N}_{m-1},z).(\bar{\eta}_{y,x,m-1}^N\!+\!E(t_{m-1})y)\big\|^{p}\!-\!\big\|(-A)^{-\frac{1-\epsilon}{4}}\mathcal N^{N}_{y,x,m-1}(t)\big\|^{p}\Big)\tilde{\mathbf{N}}(dz,dt)\bigg\}\\
&=:\|(-A)^{-\frac{1-\epsilon}{4}}\bar{\eta}_{y,x,m-1}^N\|^{p}\times (1+\mathcal G^1_{m-1,p})+\big\|E(t_{m-1})(-A)^{-\frac{1-\epsilon}{4}}y\big\|^p\times \mathcal G^1_{m-1,p}+\mathcal M_{m-1,p}(t_m),
\end{align*}
where
\[
\mathcal G^1_{m-1,p}
:=
C_\epsilon
\bigl(
\mathcal G_{m-1,p}
+\tau\mathcal G_{m-1,p-2}
+\tau
\bigr),
\]
with $\mathcal G_{m-1,\cdot}$ given by \eqref{eq:notation-G}. Define
\[
\Delta M_{m,p}
:=
\mathcal M_{m-1,p}(t_m),
\qquad
m\in\{1,\ldots,M\}.
\]
Since the frozen coefficients and the factors outside the stochastic
integrals are $\mathcal F_{t_{m-1}}$-measurable, while all stochastic
integrands are predictable, the Wiener and compensated Poisson
stochastic integrals appearing in $\Delta M_{m,p}$ have zero
conditional expectation. Hence,
\[
\mathbb E
\bigl[
\Delta M_{m,p}
\mid
\mathcal F_{t_{m-1}}
\bigr]
=0,
\qquad
m=1,\ldots,M.
\]
Thus, $\{\Delta M_{m,p}\}_{m\in\{1,\ldots,M\}}$ is a martingale difference
sequence with respect to $\{\mathcal F_{t_m}\}_{m\in\{0,\ldots,M\}}$. Set $\mathcal H_{m-1,p}:=\|E(t_{m-1})(-A)^{-\frac{1-\epsilon}{4}}y\|^p\times \mathcal G^1_{m-1,p}.$ By iteration, 
\begin{align*}
\|(-A)^{-\frac{1-\epsilon}{4}}\bar{\eta}^{N}_{y,x,m}\|^{p}&\le \|(-A)^{-\frac{1-\epsilon}{4}}\bar{\eta}^{N}_{y,x,m-1}\|^{p}\times (1+\mathcal G^1_{m-1,p})+\mathcal H_{m-1,p}+\mathcal M_{m-1,p}(t_m)\\
&\le \|(-A)^{-\frac{1-\epsilon}{4}}\bar{\eta}^{N}_{y,x,m-2}\|^{p}\times (1+\mathcal G^1_{m-2,p})+\mathcal H_{m-2,p}+\mathcal M_{m-2,p}(t_{m-1})\\
&~~~+\mathcal G^1_{m-1,p} \|(-A)^{-\frac{1-\epsilon}{4}}\bar{\eta}^{N}_{y,x,m-1}\|^{p}+\mathcal H_{m-1,p}+\mathcal M_{m-1,p}(t_m)\\
&\le \sum_{i=0}^{m-1}\mathcal G^1_{i,p}\|(-A)^{-\frac{1-\epsilon}{4}}\bar{\eta}^{N}_{y,x,i}\|^{p}+\sum_{i=0}^{m-1}\mathcal H_{i,p}+\sum_{i=0}^{m-1}\mathcal M_{i,p}(t_{i+1}).
\end{align*}
Note that $\bar{\eta}^{N}_{y,x,0}=0$. Denote $\chi_{m}:=\|(-A)^{-\frac{1-\epsilon}{4}}\bar{\eta}^{N}_{y,x,m}\|^{p}$ for $p\ge 2$, and $F_m:=\sum_{i=0}^{m-1}\mathcal H_{i,p}$,  $M_m:=\sum_{i=0}^{m-1}\mathcal M_{i,p}(t_{i+1})$, and $M_0=0$. Then we have $\chi_m\le F_m+M_m+\sum_{i=0}^{m-1}\mathcal G^1_{i,p}\chi_i,~ \forall m\in\{1,\ldots,M\}$, where $\{M_m\}_{m\in\{1,\ldots,M\}}$ is an $\{\mathcal F_{t_m}\}_{m\in\{1,\ldots,M\}}$-martingale.  By the time discrete version of the stochastic Gr\"onwall Lemma (cf. \cite[Theorem 1]{RM-2018}),  for any $q\in(0,1)$, and $\mu_1,\nu_1\in[1,\infty]$ with $\frac{1}{\mu_1}+\frac{1}{\nu_1}=1$ and $q\nu_1<1$, it holds that
 \begin{align*}
	\E[\sup_{0\le k\le m}\chi_k^{q}]\le \Big(1+\frac{1}{1-\nu_1q}\Big)^{1/\nu_1}\Big\|\prod_{k=0}^{m-1}(1+\mathcal G^1_{k,p})^q\Big\|_{L^{\mu_1}(\Omega;\R)}\big(\E[\sup_{0\le k\le m}F_k]\big)^q.
\end{align*}
In the following, we choose $q=\frac{1}{1+2\epsilon_1}$, $\epsilon_1\in(0,\epsilon]$ and $p\ge \frac{2}{q}$ such that $\tilde{p}:=pq\ge 2$. Hence, for any $\alpha_2\in[0,\frac{1}{\tilde{p}(1+2\epsilon_1)})$ such that $\alpha_2p<1$, there exist constants $L_3>0$ and $C_{L_3,\tilde{p},\epsilon,T}>0$ such that
\begin{align*}
&\E[\sup_{0\le k\le m}F_k]=\E\big[\sup_{0\le k\le m}\sum_{i=0}^{k-1}\mathcal H_{i,p}\big]=  \E\Big[\sup_{0\le k\le m}\sum_{i=0}^{k-1}\mathcal G^1_{i,p}\big\|E(t_{i})(-A)^{-\frac{1-\epsilon}{4}}y\big\|^p\Big]\\
&\le C_{\epsilon} \E\Big[\sum_{i=0}^{m-1}\big\|E(t_{i})(-A)^{-\frac{1-\epsilon}{4}}y\big\|^p\\
&\quad \times\Big(\tau+\tau\|(-A)^{\frac{1-\epsilon}{4}}\bar{X}^{N}_{i}\|^{\alpha_0(1+\epsilon)}+\tau^{\frac{p}{1+\epsilon}}\|(-A)^{\frac{1-\epsilon}{4}}\bar{X}^{N}_{i}\|^{\alpha_0p}+\tau^{\frac{p-1+\epsilon}{1+\epsilon}}\|(-A)^{\frac{1-\epsilon}{4}}\bar{X}^{N}_{i}\|^{\alpha_0(p-2)}\Big)\Big]\\
&\le C_{\epsilon}\tau\big\|(-A)^{-\frac{1-\epsilon}{4}}y\big\|^p\big(1+\|(-A)^{\frac{1-\epsilon}{4}}x\|^p\big)+C_{\epsilon}\big\|(-A)^{-\frac{1-\epsilon}{4}-\alpha_2}y\big\|^p\int_{t_1}^{t_m}\lfloor t\rfloor^{-\alpha_2p}\Big(1+\E[\|(-A)^{\frac{1-\epsilon}{4}}\bar{X}^{N}_{\lfloor t\rfloor/\tau}\|^p]\Big)\,dt\\
&\le C_{L_3,\tilde{p},\epsilon,T}\exp\Big(\frac{\|(-A)^{\frac{1-\epsilon}{4}}x\|}{4e^{L_3T}}\Big)\times \big(\tau\big\|(-A)^{-\frac{1-\epsilon}{4}}y\big\|^p+ \big\|(-A)^{-\frac{1-\epsilon}{4}-\alpha_2}y\big\|^p\big),
\end{align*}
where we used the result from Corollary \ref{cor:exponentialintegrability} in the last step. Next, using the convexity of $e^a, a\in\R$ yields that for $\epsilon\in(0,\frac12\wedge (\frac{1}{\alpha_0}-1))$
\begin{align*}
&\Big\|\prod_{k=0}^{m-1}(1+\mathcal G^1_{k,p})^q\Big\|_{L^{\mu_1}(\Omega;\R)}\le \Big\|\prod_{k=0}^{m-1}\Big\{1+C_{\epsilon}\tau\big(1+\|(-A)^{\frac{1-\epsilon}{4}}\bar{X}^{N}_{k}\|\big)^{\alpha_0(1+\epsilon)}\!+\!C_{\epsilon}\big(\tau\big(1+\|(-A)^{\frac{1-\epsilon}{4}}\bar{X}^{N}_{k}\|\big)^{\alpha_0(1+\epsilon)}\big)^{\frac{p}{1+\epsilon}}\\
&\quad~~~+C_{\epsilon}\tau \big(\tau\big(1+\|(-A)^{\frac{1-\epsilon}{4}}\bar{X}^{N}_{k}\|\big)^{\alpha_0(1+\epsilon)}\big)^{\frac{p-2}{1+\epsilon}}\Big\}^q\Big\|_{L^{\mu_1}(\Omega;\R)}\\
&\le  \Big\|\prod_{k=0}^{m-1}	\exp\Big(C_{\tilde{p},\epsilon}\tau \big(1+\|(-A)^{\frac{1-\epsilon}{4}}\bar{X}^{N}_{k}\|\big)^{\alpha_0(1+\epsilon)}\Big)\Big\|_{L^{\mu_1}(\Omega;\R)} =\! \Big\|\exp\Big(C_{\tilde{p},\epsilon}\!\sum_{k=0}^{m-1}\!\tau \big(1\!+\!\|(-A)^{\frac{1-\epsilon}{4}}\bar{X}^{N}_{k}\|\big)^{\alpha_0(1+\epsilon)}\Big)\Big\|_{L^{\mu_1}(\Omega;\R)}\\
&\le e^{C_{\tilde{p},\epsilon}t_m} \Big\|\exp\Big(\frac{1}{t_m}\!\!\sum_{k=0}^{m-1}\tau t_mC_{\tilde{p},\epsilon}\|(-A)^{\frac{1-\epsilon}{4}}\bar{X}^{N}_{k}\|^{\alpha_0(1+\epsilon)}\Big)\Big\|_{L^{\mu_1}(\Omega;\R)}\\
&\le e^{C_{\tilde{p},\epsilon}t_m} \frac{1}{t_m}\sum_{k=0}^{m-1}\tau \Big\|\exp\Big(t_mC_{\tilde{p},\epsilon} \|(-A)^{\frac{1-\epsilon}{4}}\bar{X}^{N}_{k}\|^{\alpha_0(1+\epsilon)} \Big)\Big\|_{L^{\mu_1}(\Omega;\R)}\\
&\le e^{C_{\tilde{p},\epsilon}t_m} \sup_{0\le k\le m}\Big(\E\Big[ \exp\Big(\mu_1t_mC_{\tilde{p},\epsilon} \|(-A)^{\frac{1-\epsilon}{4}}\bar{X}^{N}_{k}\|^{\alpha_0(1+\epsilon)}\Big]\Big)^{1/\mu_1}\le e^{C_{L_3,\tilde{p},\epsilon,T}}\exp\big(\frac{1}{4\mu_1e^{L_3T}}\|(-A)^{\frac{1-\epsilon}{4}}x\|\big).
\end{align*}
Consequently, we obtain that for all $\tilde{p}=pq\ge 2$ and all $\alpha_2\in[0,\frac{1}{\tilde{p}(1+2\epsilon_1)})$, 
\begin{align*}
	\E\big[\sup_{0\le k\le m}\|(-A)^{-\frac{1-\epsilon}{4}}\bar{\eta}^{N}_{y,x,k}\|^{\tilde{p}}\big]&=\E[\sup_{0\le k\le m}\chi_k^{q}]\le \Big(1+\frac{1}{1-\nu_1q}\Big)^{1/\nu_1}\Big\|\prod_{k=0}^{m-1}(1+\mathcal G^1_{k,p})^q\Big\|_{L^{\mu_1}(\Omega;\R)}\big(\E[\sup_{0\le k\le m}F_k]\big)^q\\
	&\le C_{L_3,\tilde{p},\epsilon,T}\exp\big(\frac{\|(-A)^{\frac{1-\epsilon}{4}}x\|}{2e^{L_3T}}\big) \times\big(\tau^{\frac{1}{1+2\epsilon_1}}\big\|(-A)^{-\frac{1-\epsilon}{4}}y\big\|^{\tilde{p}}+ \big\|(-A)^{-\frac{1-\epsilon}{4}-\alpha_2}y\big\|^{\tilde{p}}\big),
\end{align*}
and for all $m\in\{1,2,\ldots, M\}$,
\begin{align}\label{estimate-eta1}
	&\E\big[\|(-A)^{-\frac{1-\epsilon}{4}}\eta^{N}_{y,x,m}\|^{\tilde{p}}\big]\le 2^{\tilde{p}-1}\E\big[\|(-A)^{-\frac{1-\epsilon}{4}}\bar{\eta}^{N}_{y,x,m}\|^{\tilde{p}}\big]+2^{\tilde{p}-1}\|E(t_m)(-A)^{-\frac{1-\epsilon}{4}}y\|^{\tilde{p}}\nonumber\\
	&\le C_{L_3,\tilde{p},\epsilon,T}\big(\exp\big(\frac{\|(-A)^{\frac{1-\epsilon}{4}}x\|}{2e^{L_3T}}\big) +t_m^{-\alpha_2\tilde{p}}\big) \times \big(\tau^{\frac{1}{1+2\epsilon_1}}\big\|(-A)^{-\frac{1-\epsilon}{4}}y\big\|^{\tilde{p}}+ \big\|(-A)^{-\frac{1-\epsilon}{4}-\alpha_2}y\big\|^{\tilde{p}}\big).
\end{align}
In particular, recall that $\eta^{N}_{y,x,m}=E(m\tau)y+\tilde{\eta}^{N}_{y,x,m}$,  $m\in\{1,2,\ldots,M\}$. 
Then taking the moment estimate based on expression \eqref{expression:eta1} gives
\begin{align*}
&\|(-A)^{\frac{1-\epsilon}{4}}\tilde{\eta}^{N}_{y,x,m}\|_{L^{2}(\Omega;H)}
     \le C\int_0^{t_m}\big(1+(t_{m}-s)^{-\frac{1}{2}}\big) \Big|\big(1+\|(-A)^{\frac{1-\epsilon}{4}}\bar{X}^{N}_{\lfloor s\rfloor/\tau}\|\big)\times\big\|(-A)^{-\frac{1-\epsilon}{4}}\eta^{N}_{y,x,\lfloor s\rfloor/\tau}\big\|\Big|_{L^{2}(\Omega;\R)}\,ds\\
		&~~~+C\Big(\int_0^{t_m}\big(t_m-\lfloor s\rfloor\big)^{-1+\epsilon}\times\E\big[\big\|(-A)^{-\frac{1-\epsilon}{4}}(P_NDB(\bar{X}^{N}_{\lfloor s\rfloor/\tau}).\eta^{N}_{y,x,\lfloor s\rfloor/\tau})\big\|^2_{\IL_2^0(H)}\big]\,ds\Big)^{\frac12}\\
		&~~~+C\Big(\int_0^{t_m} \big(t_m-\lfloor s\rfloor\big)^{-1+\epsilon}\times \int_\chi\E\big[\big\|(-A)^{-\frac{1-\epsilon}{4}}(P_NDG(\bar{X}^{N}_{\lfloor s\rfloor/\tau},z).\eta^{N}_{y,x,\lfloor s\rfloor/\tau})\big\|^2\big]\,\nu(dz)ds \Big)^{\frac12}\\
		&\le C\int_0^{t_m}\big(1+(t_{m}-s)^{-\frac{1}{2}}\big)\times\big(1+\|(-A)^{\frac{1-\epsilon}{4}}\bar{X}^{N}_{\lfloor s\rfloor/\tau}\|_{L^{2+2/\epsilon_1}(\Omega;H)}\big)\times \big\|(-A)^{-\frac{1-\epsilon}{4}}\eta^{N}_{y,x,\lfloor s\rfloor/\tau}\big\|_{L^{2(1+\epsilon_1)}(\Omega;H)}\,ds\\
		&+C_{L_3,\epsilon,T}\Big(\int_{t_1}^{t_m}(t_m-\lfloor s\rfloor)^{-1+\epsilon}\big(\exp\big(\frac{\|(-A)^{\frac{1-\epsilon}{4}}x\|}{2e^{L_3T}}\big) +\lfloor s\rfloor^{-2\iota}\big)\big(\tau^{\frac{1}{1+2\epsilon_1}}\big\|(-A)^{-\frac{1-\epsilon}{4}}y\big\|^2+ \big\|(-A)^{-\frac{1-\epsilon}{4}-\iota}y\big\|^2\big)ds\Big)^{\frac12} \\
		&\times\mathbf 1_{\{m>1\}}+C_{\epsilon,T}\tau^{\frac12}\|E(t_m)(-A)^{\frac{1-\epsilon}{2}+\iota}\|_{\IL(H)}\|(-A)^{-\frac{1-\epsilon}{4}}y\|\\
		&\le C_{L_3,\epsilon,T}\big(\exp\big(\frac{\|(-A)^{\frac{1-\epsilon}{4}}x\|}{4e^{L_3T}}\big)+t_m^{-\frac{1-\epsilon}{2}-\iota}+t_{m-1}^{-\iota+\frac{\epsilon}{2}}\mathbf 1_{\{m>1\}}\big) \big(\tau^{\frac{1}{2(1+\epsilon_1)(1+2\epsilon_1)}}\big\|(-A)^{-\frac{1-\epsilon}{4}}y\big\|+ \big\|(-A)^{-\frac{1-\epsilon}{4}-\iota}y\big\|\big)
\end{align*}
with $\iota\in[0,\frac{1}{2(1+\epsilon_1)(1+2\epsilon_1)})$ and we used the condition of Assumption~\ref{assum:nonlinearity-B} in the second inequality. Hence we obtain
\begin{align}\label{estimate-eta2}
	&~~~\big\|(-A)^{\frac{1-\epsilon}{4}}\eta^{N}_{y,x,m}\big\|_{L^{2}(\Omega;H)}\le \big\|(-A)^{\frac{1-\epsilon}{4}}\tilde{\eta}^{N}_{y,x,m}\big\|_{L^{2}(\Omega;H)}+\|(-A)^{\frac{1-\epsilon}{4}}E(m\tau)y\|\\
	&\!\!\le C_{L_3,\epsilon,T}\big(\exp\big(\frac{\|(-A)^{\frac{1-\epsilon}{4}}x\|}{4e^{L_3T}}\big)\!+\!t_m^{-\frac{1-\epsilon}{2}-\iota}\!+\!t_{m-1}^{-\iota+\frac{\epsilon}{2}}\mathbf 1_{\{m>1\}}\big)\times \big(\tau^{\frac{1}{2(1+\epsilon_1)(1+2\epsilon_1)}}\big\|(-A)^{-\frac{1-\epsilon}{4}}y\big\|+ \big\|(-A)^{-\frac{1-\epsilon}{4}-\iota}y\big\|\big).\nonumber
\end{align}
In addition, for all $\epsilon_1\in(0,\epsilon]$ and $\iota_1\in[0,\frac{1-\epsilon_1}{2(1+\epsilon_1)(1+2\epsilon_1)})$,  using the interpolation and H\"older  inequalities and combining   \eqref{estimate-eta1} with $\tilde{p}=\frac{2(1+\epsilon_1)}{1-\epsilon_1}$, $\alpha_2=\iota_1$, and \eqref{estimate-eta2} with $\iota=\iota_1$ yield 
\begin{align*}
&\big\|\eta^{N}_{y,x,m}\big\|_{L^{2(1+\epsilon_1)}(\Omega;H)}\le \big\|(-A)^{\frac{1-\epsilon}{4}}\eta^{N}_{y,x,m}\big\|_{L^{2}(\Omega;H)}^{\frac{1}{2}}\big\|(-A)^{-\frac{1-\epsilon}{4}}\eta^{N}_{y,x,m}\big\|_{L^{\frac{2(1+\epsilon_1)}{1-\epsilon_1}}(\Omega;H)}^{\frac{1}{2}}\nonumber\\
&\le C_{L_3,\epsilon,T}\exp\big(\frac{\|(-A)^{\frac{1-\epsilon}{4}}x\|}{4(1+\epsilon_1)e^{L_3T}}\big)\times \big(1+t_m^{-\frac{1-\epsilon}{4}-\iota_1}+t_{m-1}^{-\frac{\iota_1}{2}+\frac{\epsilon}{4}}(1+t_m^{-\frac{\iota_1}{2}})\times \mathbf 1_{\{m>1\}}\big)\nonumber\\
&~~~\times\big(\tau^{\frac{1}{2(1+\epsilon_1)(1+2\epsilon_1)}}\big\|(-A)^{-\frac{1-\epsilon}{4}}y\big\|\!+\! \big\|(-A)^{-\frac{1-\epsilon}{4}-\iota_1}y\big\|\big)^{\frac12}\big(\tau^{\frac{1-\epsilon_1}{2(1+\epsilon_1)(1+2\epsilon_1)}}\big\|(-A)^{-\frac{1-\epsilon}{4}}y\big\|\!+\! \big\|(-A)^{-\frac{1-\epsilon}{4}-\iota_1}y\big\|\big)^{\frac12}\nonumber\\
&\le C_{L_3,\epsilon,T}\exp\big(\frac{\|(-A)^{\frac{1-\epsilon}{4}}x\|}{4(1+\epsilon_1)e^{L_3T}}\big)\times \big(1+t_m^{-\frac{1-\epsilon}{4}-\iota_1}+t_{m-1}^{-\frac{\iota_1}{2}+\frac{\epsilon}{4}}(1+t_m^{-\frac{\iota_1}{2}})\times \mathbf 1_{\{m>1\}}\big)\nonumber\\
&~~~\times\big(\tau^{\frac{1-\epsilon_1}{2(1+\epsilon_1)(1+2\epsilon_1)}}\big\|(-A)^{-\frac{1-\epsilon}{4}}y\big\|\!+\! \big\|(-A)^{-\frac{1-\epsilon}{4}-\iota_1}y\big\|\big).
\end{align*}
The proof is completed.

\subsubsection{Proof of the inequality \eqref{estimate-Y}.} \label{Subsection:6.2.1}

By interpolation and Young's inequalities, for any
$\sigma_1>\sigma>\sigma-\epsilon_1$ satisfying $2 - \frac{\epsilon_1}{\sigma_1-\sigma+\epsilon_1} - \frac{2}{1+\epsilon} \ge 0$,  $x\in\dot H^{\sigma_1}$,
and $y\in H$, we have, for $x\neq0$, 
\begin{align}\label{ineq:interpolation}
	\|x\|_{\dot{H}^{\sigma}}\|y\|&\le \|x\|^{1+\frac{\sigma_1-\sigma}{\sigma_1-\sigma+\epsilon_1}}_{\dot{H}^{\sigma-\epsilon_1}}\|x\|^{\frac{\epsilon_1}{\sigma_1-\sigma+\epsilon_1}}_{\dot{H}^{\sigma_1}}\frac{\|y\|}{\|x\|_{\dot{H}^{\sigma-\epsilon_1}}}=\|x\|^{\frac{\epsilon_1}{\sigma_1-\sigma+\epsilon_1}}_{\dot{H}^{\sigma_1}}\|x\|^{2-\frac{\epsilon_1}{\sigma_1-\sigma+\epsilon_1}}_{\dot{H}^{\sigma-\epsilon_1}}\frac{\|y\|}{\|x\|_{\dot{H}^{\sigma-\epsilon_1}}}\nonumber\\
	&\le \frac{\epsilon}{1+\epsilon}\Big(\|x\|^{\frac{\epsilon_1}{\sigma_1-\sigma+\epsilon_1}}_{\dot{H}^{\sigma_1}}\|x\|^{2-\frac{\epsilon_1}{\sigma_1-\sigma+\epsilon_1}-\frac{2}{1+\epsilon}}_{\dot{H}^{\sigma-\epsilon_1}}\Big)^{\frac{1+\epsilon}{\epsilon}}+\frac{1}{1+\epsilon}\Big(\|x\|^{\frac{2}{1+\epsilon}}_{\dot{H}^{\sigma-\epsilon_1}}\frac{\|y\|}{\|x\|_{\dot{H}^{\sigma-\epsilon_1}}}\Big)^{1+\epsilon}\nonumber\\
	&\le \frac{\epsilon}{1+\epsilon}\lambda_1^{-\frac{\sigma_1-\sigma+\epsilon_1}{2}(2-\frac{\epsilon_1}{\sigma_1-\sigma+\epsilon_1}-\frac{2}{1+\epsilon})\times \frac{1+\epsilon}{\epsilon}}\|x\|^{2}_{\dot{H}^{\sigma_1}}+\frac{1}{1+\epsilon}\|x\|^{1-\epsilon}_{\dot{H}^{\sigma-\epsilon_1}}\|y\|^{1+\epsilon}.
\end{align}
The case $x=0$ is trivial.
Let $\sigma:=-\frac{1-2\epsilon-\epsilon_0}{2}$, $\sigma_1:=1-\frac{1-\epsilon}{2}$ and $\epsilon_1:=\frac{\epsilon+\epsilon_0}{2}$ for $\epsilon_0\in(0,\epsilon]$ and $\epsilon\in(0,\frac12\wedge\frac{1}{(\lambda_1\wedge 1)^{-3/4}-1})$. Then we obtain
\begin{align*}
	&~~~\|(-A)^{-\frac{1-\epsilon}{4}}Y^{N}_{y,x,m-1}(t_m)\|^{p}-\|(-A)^{-\frac{1-\epsilon}{4}}\bar{\eta}_{y,x,m-1}^N\|^{p}\\
	&\le \int_{t_{m-1}}^{t_m}p\|(-A)^{-\frac{1-\epsilon}{4}}Y^{N}_{y,x,m-1}(t)\|^{p-2}\Big(-\big\|(-A)^{\frac12-\frac{1-\epsilon}{4}}Y^{N}_{y,x,m-1}(t)\big\|^2+\big\|(-A)^{-\frac{1-\epsilon}{4}+\frac{\epsilon+\epsilon_0}{4}}Y^{N}_{y,x,m-1}(t)\big\|\\
	&~~~\times\big\|(-A)^{-\frac{1+\epsilon_0}{4}}P_NDF(\bar{X}^{N}_{m-1}).(\bar{\eta}_{y,x,m-1}^N+E(t_{m-1})y)\big\|\Big)dt\\
&\le  \int_{t_{m-1}}^{t_m}\!\!\!p\|(-A)^{-\frac{1-\epsilon}{4}}Y^{N}_{y,x,m-1}(t)\|^{p-2}\Big(\!-\big\|(-A)^{\frac12-\frac{1-\epsilon}{4}}Y^{N}_{y,x,m-1}(t)\big\|^2+C_{\epsilon,\epsilon_0}\big\|(-A)^{\frac12-\frac{1-\epsilon}{4}}Y^{N}_{y,x,m-1}(t)\big\|^{\frac{\epsilon+\epsilon_0}{2}}\times\\
	& \|(-A)^{-\frac{1-\epsilon}{4}}Y^{N}_{y,x,m-1}(t)\big\|^{1-\frac{\epsilon+\epsilon_0}{2}}\big(1+\|(-A)^{\frac{1-\epsilon}{4}}\bar{X}^{N}_{m-1}\|\big)^{\alpha_0}\big[\big\|(-A)^{-\frac{1-\epsilon}{4}}\bar{\eta}_{y,x,m-1}^N\|+\|E(t_{m-1})(-A)^{-\frac{1-\epsilon}{4}}y\|\big]\Big)dt\\
&\le C_{p}\int_{t_{m-1}}^{t_m}\!\!\!\!\|(-A)^{-\frac{1-\epsilon}{4}}Y^{N}_{y,x,m-1}(t)\|^{p-1-\epsilon} \Big(\|E(t_{m-1})(-A)^{-\frac{1-\epsilon}{4}}y\|^{1+\epsilon}\!+\!\big\|(-A)^{-\frac{1-\epsilon}{4}}\bar{\eta}_{y,x,m-1}^N\big\|^{1+\epsilon}\Big)dt\\
&\quad \times\big(1\!+\!\|(-A)^{\frac{1-\epsilon}{4}}\bar{X}^{N}_{m-1}\|\big)^{\alpha_0(1+\epsilon)}\\
&\le C_{p}\int_{t_{m-1}}^{t_m}\!\!\!\|(-A)^{-\frac{1-\epsilon}{4}}\bar{\eta}_{y,x,m-1}^N\|^{p-1-\epsilon} \Big(\|E(t_{m-1})(-A)^{-\frac{1-\epsilon}{4}}y\|^{1+\epsilon}\!+\!\big\|(-A)^{-\frac{1-\epsilon}{4}}\bar{\eta}_{y,x,m-1}^N\big\|^{1+\epsilon}\Big)dt\\
&\quad\times\big(1\!+\!\|(-A)^{\frac{1-\epsilon}{4}}\bar{X}^{N}_{m-1}\|\big)^{\alpha_0(1+\epsilon)}+C_{p}\int_{t_{m-1}}^{t_m}\big\|(-A)^{-\frac{1-\epsilon}{4}}(Y^{N}_{y,x,m-1}(t)-\bar{\eta}_{y,x,m-1}^N)\big\|^{(1+\epsilon)(\frac{p}{1+\epsilon}-1)}\\
&\quad\times \Big(\|E(t_{m-1})(-A)^{-\frac{1-\epsilon}{4}}y\|^{1+\epsilon}+\big\|(-A)^{-\frac{1-\epsilon}{4}}\bar{\eta}_{y,x,m-1}^N\big\|^{1+\epsilon}\Big)dt\times\big(1+\|(-A)^{\frac{1-\epsilon}{4}}\bar{X}^{N}_{m-1}\|\big)^{\alpha_0(1+\epsilon)},
\end{align*}
 where constant $C_{\epsilon,\epsilon_0}>0$ is given by the condition on $DF$ in Assumption \ref{assum:nonlinearity-F} (ii) with $\varsigma_0=\epsilon_0$ and $\varsigma=\epsilon$. Indeed, using this condition with $\epsilon_0\le\epsilon$, from the expression of $Y^{N}_{y,x,m-1}(t)$, it follows that
\begin{align*}
&\big\|(-A)^{-\frac{1-\epsilon}{4}}(Y^{N}_{y,x,m-1}(t)-\bar{\eta}_{y,x,m-1}^N)\big\|^{1+\epsilon}	\le 2\|(-A)^{-\frac{1-\epsilon}{4}}\bar{\eta}_{y,x,m-1}^N\|^{1+\epsilon}+C_{\epsilon}(t-t_{m-1})^{(1-\frac{\epsilon+\epsilon_0}{4})(1+\epsilon)} \\
&\times\big(1+\|(-A)^{\frac{1-\epsilon}{4}}\bar{X}^{N}_{m-1}\|\big)^{\alpha_0(1+\epsilon)}\times \big(\big\|E(t_{m-1})(-A)^{-\frac{1-\epsilon}{4}}y\big\|^{1+\epsilon}+\big\|(-A)^{-\frac{1-\epsilon}{4}}\bar{\eta}_{y,x,m-1}^N\big\|^{1+\epsilon}\big).
\end{align*}
Here $(t-t_{m-1})^{(1-\frac{\epsilon+\epsilon_0}{4})(1+\epsilon)}\le \tau$ for all $t\in[t_{m-1},t_m]$ due to $\epsilon_0\le\epsilon$. Then combining this result and applying Young's inequality once more  gives the desired assertion.

\bibliographystyle{abbrv} 
  \bibliography{references.bib}

\end{document}